\documentclass[12pt, eqno]{article}
\usepackage{amssymb,latexsym}
\usepackage{amsmath,latexsym}
\usepackage{amscd,latexsym} 
\usepackage[all]{xy}

\newcommand{\N}{{\mathbb N}}
\newcommand{\C}{{\mathbb C}}
\newcommand{\R}{{\mathbb R}}

\newcommand{\Z}{{\mathbb Z}}
\newcommand{\CP}{{\mathbb CP}}
\newcommand{\RP}{{\mathbb RP}}

\newtheorem{theorem}{Theorem}[section]

\newtheorem{corollary}[theorem]{Corollary}

\newtheorem{remark}[theorem]{Remark}

\newtheorem{lemma}[theorem]{Lemma}

\newtheorem{proposition}[theorem]{Proposition}
\newtheorem{claim}[theorem]{Claim}

\begin{document}

\title{Deforming symplectomorphism of certain irreducible Hermitian symmetric spaces of compact type
\\ by mean curvature flow}

\author{Guangcun Lu
\thanks{Partially supported by the NNSF  10671017 and 10971014 of
China, and PCSIRT and Research Fund for the Doctoral
Program Higher Education of China (Grant No. 200800270003).}\hspace{10mm}Bang Xiao\\
{\normalsize School of Mathematical Sciences, Beijing Normal University},\\
{\normalsize Laboratory of Mathematics
 and Complex Systems,  Ministry of
  Education},\\
  {\normalsize Beijing 100875, The People's Republic
 of China}\\
{\normalsize (gclu@bnu.edu.cn\hspace{5mm}bangxiao@mail.bnu.edu.cn)}}
\date{Preliminary version October 20, 2010,\\ The revision July 5, 2011}
\maketitle \vspace{-0.1in}

\abstract{In this paper, we generalize Medos-Wang's  arguments and
results on the mean curvature flow  deformations of
symplectomorphisms of $\CP^n$ in \cite{MeWa} to complex Grassmann
manifold $G(n, n+m;\C)$ and compact totally geodesic
K\"ahler-Einstein submanifolds of $G(n, 2n;\C)$ such as irreducible
Hermitian symmetric spaces  $SO(2n)/U(n)$ and $Sp(n)/U(n)$ (in the
terminology of \cite[p. 518]{He}).  We also give an abstract result
and discuss the
case of complex tori. }\\



\tableofcontents

\section{Introduction}
\setcounter{equation}{0}

 Recall that a symplectic
manifold $(M, \omega)$ is said to be K\"{a}hler if there exists an
integrable almost complex structure $J$ on $M$ such that the
bilinear form
$$
g(X,Y)=\langle X, Y\rangle:=\omega(X, JY)
$$
defines a Riemannian metric on $M$. The triple $(\omega, J, g)$ is
called a K\"{a}hler structure on $M$,  $g$ and $\omega$ are called a
K\"{a}hler metric and a K\"ahler form, respectively.  Such a
K\"{a}hler manifold is called a {\it K\"{a}hler-Einstein manifold}
if the Ricci form $\rho_\omega\equiv \rho_g$ of $g$ satisfies
$\rho_\omega= c \omega$ for some constant $c\in\R$. For a K\"{a}hler
manifold  $(M, J, g, \omega)$ let ${\rm Symp}(M,\omega)$ and ${\rm
Aut}(M, J)$ denote the group of symplectomorphisms of the symplectic
manifold $(M,\omega)$ and the  group of biholomorphisms of the
complex manifold $(M, J)$,
 respectively. Their intersection is equal to
 the group of isometries of the K\"ahler manifold $(M, J, g, \omega)$,
 ${\rm I}(M, J, g):=\{\phi\in{\rm Aut}(M, J)\,|\,\phi^\ast g=g\}$.

Without special statements we always assume that $M$ is closed (i.e.
compact and boundaryless) and connected throughout this paper. It is
well-known that ${\rm Symp}(M,\omega)$ is an infinite dimensional
Lie group whose Lie algebra is the space of symplectic vector
fields.  A lot of symplectic topology information of $(M, \omega)$
is contained in ${\rm Symp}(M,\omega)$. (See beautiful books
\cite{Ban, HoZe, McSa, Po} for detailed study). On the other hand
${\rm I}(M, J, g)$ is a finite dimensional Lie subgroup of ${\rm
Symp}(M,\omega)$. Hence in order to understand topology of ${\rm
Symp}(M,\omega)$, e.g. its homotopy groups, it is helpful to study
the topology properties of the inclusion ${\rm I}(M, J,
g)\hookrightarrow{\rm Symp}(M,\omega)$. To the author's knowledge,
the first result in this direction was obtained by Smale \cite{Sma},
who proved that there exists a continuous strong deformation
retraction from ${\rm Symp}(S^2,\omega_{\rm FS}^{(1)})$ to
$SO(3)={\rm I}(S^2, i, g_{\rm FS}^{(1)})$. Hereafter $g_{\rm
FS}^{(n)}$ and $\omega_{\rm FS}^{(n)}$ denote, up to multiplying  a
positive number, the Fubini-Study metric and the associated K\"ahler
form on the complex projective spaces $\CP^n$ respectively, and $i$
is the standard complex structure on $\CP^n$.  Recently, Jiayong Li
and Jordan Alan Watts \cite{LiWat} strengthened this result. They
constructed a strong deformation retraction from ${\rm
Symp}(S^2,\omega_{\rm FS}^{(1)})$ to $SO(3)$ which is
diffeologically smooth.

In his famous paper \cite{Gr} Gromov invented a powerful
pseudo-holomorphic curve theory to study symplectic topology and got
the following important
results: \\
$\bullet$ For any two area forms $\omega_1$ and $\omega_2$ on
$\CP^1$ with $\int_{\CP^1}\omega_1=\int_{\CP^1}\omega_2$,  ${\rm
Symp}(\CP^1\times\CP^1, \omega_1\oplus\omega_2)$ contracts onto
\begin{eqnarray*}
{\rm I}(\CP^1\times\CP^1, i\times i, g_{\rm FS}^{(1)}\oplus g_{\rm
FS}^{(1)})=\Z/2\Z\; \hbox{\bf extension of}\; SO(3)\times SO(3)
\end{eqnarray*}
(\cite[\S 2.4.$A_1$]{Gr}),   and ${\rm Symp}(\CP^1\times\CP^1,
\omega_1\oplus\omega_2)$ cannot contract onto $SO(3)\times SO(3)$ if
$\int_{\CP^1}\omega_1\ne\int_{\CP^1}\omega_2$ (\cite[\S
2.4.$C_2$]{Gr}). (A simple application of Moser theorem can reduce
these to the case $\omega_1=a\omega_{\rm FS}^{(1)}$ and
$\omega_2=b\omega_{\rm FS}^{(1)}$ for nonzero
$a,b\in\R$).\\
$\bullet$ ${\rm Symp}(\CP^2, \omega_{\rm FS}^{(2)})$ contracts onto
 ${\bf I}(\CP^2, i, g_{\rm FS}^{(2)})$ (\cite[\S 2.4.$B'_3$]{Gr}).

Since ${\rm Symp}(S^2,\omega_{\rm FS}^{(1)})={\rm Diff}_+(S^2)$,
Gromov's results may be viewed as generalizations of Smale's theorem
above in a direction. For ${\rm Symp}(S^2\times S^2, \omega_{\rm
FS}^{(1)}\oplus\lambda\omega_{\rm FS}^{(1)})$ with
$\int_{S^2}\omega_{\rm FS}^{(1)}=1$ and $\lambda\ne 1$,  so far some
deep results were made by Abreu \cite{Ab}, Abreu and McDuff
\cite{AbMc}, Anjos and Granja \cite{AnGr} and others following an
approach suggested by Gromov \cite[\S 2.4.$C_2$]{Gr}.  A different
direction generalization of Smale's theorem is to study the topology
properties of the inclusion ${\rm Symp}(M,\omega)\hookrightarrow{\rm
Diff}(M,\omega)$. Using the parameterized Gromov-Witten invariant
theory, L\^{e} and  Ono \cite{LeO}, and Seidel \cite{Sei}
 got a few of interesting results in this direction.  For example,
  with $P_{mn}:=\CP^m\times \CP^n$ and
$\int_{\CP^1}\omega_{\rm FS}^{(m)}=\int_{\CP^1}\omega_{\rm
FS}^{(n)}=1$, it was showed in \cite{Sei} that the homomorphisms
$$\beta_{k}\colon \pi_{k}({\rm Symp}(\CP^m\times \CP^n, \omega_{\rm
FS}^{(m)}\oplus\omega_{\rm FS}^{(n)}))\to \pi_{k}({\rm
Diff}(P_{mn}))$$ induced by the natural inclusion ${\rm
Symp}(\CP^m\times \CP^n, \omega_{\rm FS}^{(m)}\oplus\omega_{\rm
FS}^{(n)})\subset{\rm Diff}(P_{mn})$ are not surjective for odd
numbers $k\le\max\{2m-1,2n-1\}$. This gave the first examples of
symplectic manifolds of dimension $> 4$ for which the map
$\pi_1({\rm Symp})\to \pi_1({\rm Diff})$ is not surjective. In
particular, this result implies that ${\rm rank}({\rm
coker}\beta_1)\ge 2$ for $m = n =1$, which can also be derived from
the above Gromov's first result. (See McDuff's survey \cite{Mc} for
recent developments).

In past ten years  a new method ({\it mean curvature flow} (MCF)
method) to the above question was developed by  Mu-Tao Wang
\cite{Wa1, Wa2, Wa3, Wa4, Wa5, TsWa, MeWa} and  Smoczyk \cite{Smo2}.
For compact Riemann surfaces they obtained the desired results (cf.
\cite{Wa4, Wa5, Smo2}). Recently Ivana Medos and Mu-Tao Wang
\cite{MeWa} applied the MCF  to deform symplectomorphisms of
$\CP^{n}$ for each dimension $n$, and obtained a constant
$\Lambda_0(n)\in (1, +\infty]$ only depending on $n\in\N$, (see
(\ref{e:3.7}) for its definition),
 such that any $\Lambda$-pinched symplectomorphism of
$\CP^{n}$ with

 {\scriptsize
 \begin{equation}\label{e:1.1}
1\le
\Lambda\le\Lambda_1(n):=\biggl[\frac{1}{2}\biggl(\Lambda_0(n)+
\frac{1}{\Lambda_0(n)}\biggr)\biggl]^{\frac{1}{n}}+ \sqrt{
\biggl[\frac{1}{2}\biggl(\Lambda_0(n)+
\frac{1}{\Lambda_0(n)}\biggr)\biggl]^{\frac{2}{n}}-1}
\end{equation}}
 is symplectically
isotopic to a biholomorphic isometry (cf.\cite[Corollay 5]{MeWa}).
 Here a symplectomorphism
$\varphi$ of the K\"{a}hler manifold $(M, \omega, J, g)$ is called
{\it $\Lambda$-pinched} if
$$
\frac{1}{\Lambda^2}g\leq \varphi^{\ast}g\leq\Lambda^2 g
$$
(cf. \cite[Def.1]{MeWa}).  The constant $\Lambda_0(n)$  was
introduced above Remark 2 of \cite[p.322]{MeWa}, and it was shown
that $\Lambda_0(1)=\infty$ there. So \cite{MeWa} provides a new
proof for the Smale's result above. For $n\in\N$ we define an
increasing function $[1, \infty)\ni\Lambda\mapsto\Lambda'_n$ by

{\scriptsize
 \begin{equation}\label{e:1.2}
\Lambda'_n:=\biggl[\frac{1}{2}\Bigl(\Lambda+
\frac{1}{\Lambda}\Bigr)\biggr]^n+
\sqrt{\biggl[\frac{1}{2}\Bigl(\Lambda+
\frac{1}{\Lambda}\Bigr)\biggr]^{2n}-1}.
\end{equation}}
(This is obtained from \cite[(3.11)]{MeWa} when $\Lambda_1$ in
\cite[(3.10)]{MeWa} is replaced by $\Lambda$.) Then
$\Lambda'_n=\Lambda_0(n)$ if $\Lambda=\Lambda_1(n)$ by the proof of
\cite[Cor.5]{MeWa}.

 By Cartan's classification, in addition to two exceptional spaces
$E_6/({\rm Spin}(10)\times{\rm SO}(n+2))$ and $E_7/(E_6\times{\rm
SO}(2))$, all irreducible Hermitian symmetric spaces of compact type
have the following form of four types (in the terminology of
\cite[p. 518]{He}):
\begin{eqnarray*}
&&U(n+m)/U(n)\times U(m),\;n, m\ge 1, \quad SO(2n)/U(n), \;n\ge 2,\\
&& Sp(n)/U(n)\; n\ge 2, \qquad SO(n+2)/SO(n)\times SO(2),\;n\ge 3.
\end{eqnarray*}
 They are, respectively,  holomorphically equivalent to  (cf. \cite{CaVe}):\\
$\bullet$ $G^{\rm I}(n, n+m)=G(n, n+m;\C)$ the complex Grassmann
manifold  consisting of the $(n-1)$-dimensional, complex projective
linear subspaces $\CP^{n-1}$ of the complex projective $(n+
m-1)$-space
$\CP^{n+m-1}$;\\
$\bullet$ $G^{\rm II}(n, 2n)$ the complex analytic submanifold of
$G(n, 2n;\C)$ consisting of all the $(n-1)$-dimensional, complex
projective linear subspaces $\CP^{n-1}$ that lie in a non- singular
quadric
hypersurface $Q^{2n-2}(\C)$ in $\CP^{2n-1}$;\\
$\bullet$ $G^{\rm III}(n, 2n)$ the complex subvariety of $G(n,
2n,\C)$ consisting of all the $(n-1)$-dimensional, complex
projective linear subspaces $\CP^{n-1}\subset\CP^{2n-1}$, all of
whose projective lines are contained in a general (i.e.,
nonsingular) linear complex. In terms of homogeneous coordinates in
$\CP^{2n-1}$, it is represented by the m-dimensional, complex vector
subspaces of the vector space $\C^{2n}$, that are totally isotropic
with respect to a nonsingular, alternating bilinear form
on $\C^{2n}$; \\
$\bullet$ $G^{\rm IV}(1, n+1)$ the nonsingular, $n$-dimensional,
complex quadric hypersurface $Q_n\subset\CP^{n+1}$, and
real-analytically isomorphic to the real Grassmann manifold $G^+(2,
n,\R)$ of oriented, real projective lines in $\RP^{n+1}$.

They have complex dimensions $mn, n(n-1)/2, n(n+1)/2, n$,
respectively. Let $h$ and $h_{\rm I}$ be the canonical K\"ahler
metrics on $G(n, n+m;\C)$ and  $G^{\rm I}(n, 2n)$, respectively.
Both $G^{\rm II}(n, 2n)$ and $G^{\rm III}(n, 2n)$ are totally
geodesic submanifolds of $(G^{\rm I}(n, 2n), h_I)$. Denote by
$h_{\rm II}$ and $h_{\rm III}$ the induced metrics on $G^{\rm II}(n,
2n)$ and $G^{\rm III}(n, 2n)$, respectively.

\begin{theorem}\label{th:1.1}
Let $\omega$ be the K\"ahler form corresponding with the canonical
metric $h$ on $G(n, n+m;\C)$, $g={\rm Re}(h)$ and $J$ the standard
complex structure.  Then for every $\Lambda$-pinched
symplectomorphism $\varphi\in{\rm Symp}(G(n, n+m;\C),\omega)$ with
$\Lambda\in [1, \Lambda_1(mn)]\setminus\{\infty\}$
 the following holds:
\begin{description}
\item[(i)] The mean curvature flow $\Sigma_t$ of the graph of $\varphi$ in
$G(n, n+m;\C)\times G(n, n+m;\C)$ exists for all $t>0$.

\item[(ii)] $\Sigma_{t}$ is the graph of a symplectomorphism $\varphi_t$ for each
$t>0$, and $\varphi_{t}$ is $\Lambda'_{mn}$-pinched along the mean
curvature flow, where $\Lambda'_{mn}$ is defined by (\ref{e:1.2}).

\item[(iii)] $\varphi_{t}$ converges smoothly to a biholomorphic isometry of
$(G(n, n+m;\C), J, g)$ as $t\rightarrow\infty$.
\end{description}
Consequently, each such $\Lambda$-pinched symplectomorphism
$\varphi\in{\rm Symp}(G(n, n+m;\C),\omega)$
 is symplectically isotopic to a biholomorphic isometry of $(G(n, n+m;\C), J, g)$.
\end{theorem}

\begin{theorem}\label{th:1.2}
Let $(M, \omega, J, g)$ be a compact K\"ahler-Einstein submanifold
of $(G(n, n+m;\C), h)$ which is totally geodesic (e.g. $(G^{\rm
II}(n, 2n), h_{\rm II})$ and $(G^{\rm III}(n, 2n), h_{\rm III})$ are
such submanifolds of $(G^{\rm I}(n, 2n), h_{\rm I})$). Set $\dim
M=2N$. Then for every $\Lambda$-pinched symplectomorphism
$\varphi\in{\rm Symp}(M,\omega)$ with $\Lambda\in [1,
\Lambda_1(N)]\setminus\{\infty\}$
 the following holds:
\begin{description}
\item[(i)] The mean curvature flow $\Sigma_t$ of the graph of $\varphi$ in
$M\times M$ exists for all $t>0$.

\item[(ii)] $\Sigma_{t}$ is the graph of a symplectomorphism $\varphi_t$ for each
$t>0$, and $\varphi_{t}$ is  $\Lambda'_N$-pinched along the mean
curvature flow, where $\Lambda'_{N}$ is defined by (\ref{e:1.2}).

\item[(iii)] $\varphi_{t}$ converges smoothly to a biholomorphic isometry of
$(M, J, g)$ as $t\rightarrow\infty$.
\end{description}
Consequently, each such $\Lambda$-pinched symplectomorphism
$\varphi: (M,\omega)\rightarrow (M,\omega)$
 is symplectically isotopic to a biholomorphic isometry of $(M, J, g)$.
\end{theorem}

Recall that a complex torus of complex dimension $n$ is the quotient
space $T^n=\C^n/\Gamma$, where $\Gamma$ is a lattice in $\C^n$
generated by  $2n$ vectors $\{u_1,\cdots,u_{2n}\}$ in $\C^n$ which
are linearly independent over $\R$. It has a natural flat K\"ahler
metric induced from the flat metric of $\C^n$. By Bieberbach theorem
(\cite[page 65]{Ch}), any compact flat K\"ahler manifold is
holomorphically covered by a complex torus (\cite[Example
2.60]{Be}). From this and Calabi-Yau theorem it follows that any
compact K\"ahler manifold $M$ with the first and the second (real)
Chern class vanishing must be (holomorphically) covered by a complex
torus (\cite[Cor. 11.27]{Be}). Unfortunately, for complex tori  we cannot obtain the
corresponding result with (iii) of Theorems~\ref{th:1.1} and \ref{th:1.2} yet though other conclusions are proved
under the weaker pinching
condition.

\begin{theorem}\label{th:1.3}
Let $(M, \omega, J, g)$ and $(\widetilde{M}, \tilde\omega, \tilde J,
\tilde g)$ be two real $2n$-dimensional compact K\"{a}hler-Einstein
manifolds of {\rm constant zero holomorphic sectional curvature}.
Then for every $\Lambda$-pinched symplectomorphism $\varphi:
M\rightarrow\widetilde{M}$ with $\Lambda\in (1, \Lambda_0(n))$ there
hold:
\begin{description}
\item[(i)] The mean curvature flow $\Sigma_{t}$ of the graph of $\varphi$ in
$M\times\widetilde{M}$ exists smoothly for all $t>0$;

\item[(ii)] $\Sigma_{t}$ is the graph of a symplectomorphism $\varphi_{t}$ for
each $t>0$, and $\varphi_{t}$ is still $\Lambda_{0}(n)$-pinched
along the mean curvature flow.

\item[(iii)] If $\Lambda<\widehat\Lambda_1$ for some $\Lambda_1\in
(\Lambda, \Lambda_0(n))$, where $\widehat\Lambda_1>1$ is a constant
determined by $\Lambda_1$ and\, $n$ (see Lemma~\ref{lem:5.2}), then
 the flow converges to a totally
geodesic submanifold of $M\times\widetilde{M}$ as
$t\rightarrow\infty$. (In addition we have
$$
\widehat\Lambda_1\ge
\left(2\exp\Bigl(\frac{0.141446\delta_{\Lambda_1}}{5n} \Bigr)+
2\exp\Bigl(\frac{0.141446\delta_{\Lambda_1}}{10n}
\Bigr)\sqrt{\exp\Bigl(\frac{0.141446\delta_{\Lambda_1}}{5n} \Bigr)
-1}-1\right)^{1/2},
$$
where $\delta_{\Lambda_1}$ is defined by (\ref{e:3.6})).
\end{description}
\end{theorem}

It is easily seen that the convergence assertion in
Theorem~\ref{th:1.3} cannot be derived from \cite[Theorem B]{Wa3}.
Moreover, it was pointed out in \cite[Remark 8.1]{Wa2} that when $M$
is locally a product of two Riemannian surfaces of nonpositive
curvature the uniform convergence of the flow can also be proved
with the method in \cite{Wa4}.

It is possible to generalize the above three theorems  to a larger
class of manifolds---- compact homogeneous K\"{a}hler-Einstein
manifolds. (See Theorem~\ref{th:6.1}). Recall that a K\"{a}hler
manifold $(M, \omega, J, g)$ is called {\it homogeneous} if ${\rm
I}(M, J, g)$ acts transitively on $M$. In particular, a
simply-connected compact homogeneous K\"ahler manifold is called a
{\it K\"ahler $C$-space} in \cite{W} (or a {\it generalized flag
manifold}).  However, except the manifolds contained in the three
theorems above we do not find an example satisfying the conditions
of Theorem~\ref{th:6.1}.

In this paper  we follow \cite{KoNo} to define the curvature tensor
$R$ of a K\"ahler manifold $(M, \omega, J, g)$   by
$$
R(X, Y, Z, W)=g(R(X,Y)W, Z)=g(R(Z,W)Y, X)
$$
for $X, Y, Z, W\in\Gamma(TM)$. Then the {\it holomorphic sectional
curvature}  in the direction $X\in TM\setminus\{0\}$ is
 defined by $H(X)=R(X, JX,X, JX)/[g(X,X)]^2$.
(After extending $g$ and $R$ by $\C$-linearity to
$TM\otimes_{\R}\C$, $H(X)$ is equal to $-R(Z, \overline{Z}, Z,
\overline{Z})/[g(Z,\overline{Z})]^2$ for $Z=(X-\sqrt{-1}JX)/2\in
T^{(1,0)}M$).

The paper is organized as follows. In Section 2 we review
differential geometry of Grassmann manifolds, the key
Proposition~\ref{prop:2.3} seems to be new.  Section 3 is our
technical core, where we study evolution along the mean curvature
flow under different pinching conditions for different cases. In
Section 4 we prove Theorems~\ref{th:1.1}, \ref{th:1.2} and
\ref{th:1.3}. Finally, Section 5 gives a general result under
stronger assumptions as a concluding remark. \vspace{2mm}

\

\noindent{\bf Acknowledgments}. The authors would like to thank the
anonymous referees for pointing out errors in the arguments of
improving pinching condition, a number of typos in a previous
version, and suggestions in improving the presentation.

\section{Differential geometry of Grassmann
manifolds} \setcounter{equation}{0}

We shall review some necessary results in differential geometry of
Grassmann manifolds. Some of them are first observed.

\subsection{Grassmann manifold and curvature}

  Let $M(n+m, n;\C)=\{ A\in\C^{n\times (n+m)}\,|\, {\rm
  rank}A=n\,\}$. Then
${\rm GL}(n;\C):=\{Q\in\C^{n\times n}\,|\,{\rm det}Q\ne 0\}$ acts
freely on $M(n, n;\C)$ from the left by matrix multiplication. The
complex Grassmann manifold $G(n, n+m;\C)$ may be defined  as the
quotient $M(n, n+m;\C)/{\rm GL}(n;\C)$. For $A\in M(n, n+m;\C)$,
denote by $[A]\in G(n, n+m;\C)$ the ${\rm GL}(n;\C)$-orbit of $A$ in
$M(n, n+m;\C)$.
 Any representative matrix $B$ of $[A]$ is
called a {\it homogeneous coordinate} of the point $[A]$. For
increasing integers $1\le\alpha_1<\cdots<\alpha_n\le n+m$ let
$\{\alpha_{n+1},\cdots, \alpha_{n+m}\}$ be the complement of
$\{\alpha_1,\cdots,\alpha_n\}$ in the set $\{1,2,\dots, n+m\}$.
Write $A\in M(n, n+m;\C)$ as $A=(A_1,\cdots, A_{n+m})$ and
\begin{eqnarray*}
A_{\alpha_1\cdots\alpha_n}=(A_{\alpha_1},\cdots,
A_{\alpha_n})\in\C^{n\times n},\quad
A_{\alpha_{n+1}\cdots\alpha_{n+m}}=(A_{\alpha_{n+1}},\cdots,
A_{\alpha_{n+m}}) \in\C^{n\times m},
\end{eqnarray*}
where $A_1,\cdots, A_{n+m}$ are $n\times 1$ matrices. Define
\begin{eqnarray*}
&&U_{\alpha_1,\cdots,\alpha_n}=\{[A]\in G(n, n+m;\C)\,|\, {\rm
det}A_{\alpha_1\cdots\alpha_n}\ne 0\,\},\\
&&\Theta_{\alpha_1\cdots\alpha_n}: U_{\alpha_1\cdots \alpha_n}\to
\C^{n\times m}\equiv\C^{nm},\quad [A]\to
Z=(A_{\alpha_1\cdots\alpha_n})^{-1}A_{\alpha_{n+1}\cdots\alpha_{n+m}}.
\end{eqnarray*}
 We call $Z$ the {\it local
coordinate} of $[A]\in G(n, n+m;\C)$, and
$$
 \bigl\{\bigl(U_{\alpha_1\cdots\alpha_n},\,
\Theta_{\alpha_1\cdots\alpha_n}\bigr)\, \bigm|\,
1\le\alpha_1<\cdots<\alpha_n\le n\bigr\}
$$
the {\it canonical atlas}  on $G(n, n+m;\C)$ (\cite{Le,Lu1, Wo2}).
There exists a unique (up to multiplying a nonzero real constant)
metric $h$ on $G(n, n+m;\C)$ which is invariant under the action of
the group of motions in $G(n, n+m;\C)$ (cf.\cite{Le}). It is the
pullback of the Fubini-Study metric on ${\bf P}(\wedge^n\C^{n+m})$
by the Pl\"ucker embedding from $G(n,m;\C)$ to ${\bf
P}(\wedge^n\C^{n+m})$). It is well-known (e.g. \cite{Lu1, Lu2} )
that this metric in the local chart $(U_{1\cdots n},
Z=\Theta_{1\cdots n})$ on $G(n, n+m;\C)$ is given by
\begin{eqnarray}\label{e:2.1}
h&=&{\rm Tr}[(I_n+
Z\overline{Z}')^{-1}dZ(I_m+\overline{Z}'Z)^{-1}\overline{dZ}']\\
&=&\partial\bar \partial\log\det(I+ Z\overline{Z}'),\nonumber
\end{eqnarray}
where $\overline{Z}'$ and $\overline{dZ}'$ are the conjugate
transposes of $Z$ and $dZ$ respectively, ${\rm Tr}$ denotes the trace, and
$\partial=\sum_{i,\alpha}dZ^{i\alpha}\frac{\partial}{\partial
Z^{i\alpha}}$ and $\bar
\partial=\sum_{i,\alpha}\overline{dZ^{i\alpha}}\frac{\partial}{\partial
Z^{i\alpha}}$. It is  K\"ahler-Einstein.

 If a  (real) tangent
vector $T$ at the point $Z\in U_{1\cdots n}$ is represented by their
component matrices, i.e., we identify
\begin{eqnarray}\label{e:2.2}
T=\sum_{k,l}{\rm Re}(T^{kl})\frac{\partial}{\partial X^{kl}} +
\sum_{k,l}{\rm Im}(T^{kl})\frac{\partial}{\partial Y^{kl}}
\end{eqnarray}
with complex matrices $T=(T^{kl})\in \C^{n\times m}$, where
$Z^{kl}=X^{kl}+ iY^{kl}$, $k=1,\cdots,n$ and $l=1,\cdots, m$, then
the Riemannian metric $g:=Re(h)$ is given by
\begin{equation}\label{e:2.3}
g_Z(T_1, T_2)=Re {\rm Tr}[(I+
Z\overline{Z}')^{-1}T_1(I+\overline{Z}'Z)^{-1}\overline{T_2}']
\end{equation}
(cf. \cite[(2)]{Wo2}). By \cite[(4)]{Wo2}), the curvature tensor $R_Z$
of $g$ at $Z$ has the expression
\begin{eqnarray*}
R_Z(T_1, T_2)T\!\! &=&\!\!T\bigl[(I+
\overline{Z}'Z)^{-1}\overline{T_2}'(I+ Z\overline{Z}')^{-1}T_1 - (I+
\overline{Z}'Z)^{-1}\overline{T_1}'(I+
Z\overline{Z}')^{-1}T_2\bigr]\nonumber\\
&+&\!\!\!\bigl[T_1(I+ \overline{Z}'Z)^{-1}\overline{T_2}'(I+
Z\overline{Z}')^{-1} - T_2(I+ \overline{Z}'Z)^{-1}\overline{T_1}'(I+
Z\overline{Z}')^{-1}\bigr]T.\nonumber
\end{eqnarray*}
Here as above the left is a real tangent vector and the right  is
the corresponding complex matrix representation of it. Let $p_0\in
U_{1\cdots n}$ has coordinate $Z(p_0)=0$. Then
\begin{eqnarray}
&&R_{p_0}(T_1, T_2, T_3, T_4):=g_p(R_p(T_3, T_4)T_2, T_1)\nonumber\\
&=&Re {\rm
Tr}\bigl[(T_2\overline{T}'_4T_3\overline{T}_1-T_2\overline{T}'_3T_4\overline{T}'_1
+T_3\overline{T}_4T_2\overline{T}_1-T_4\overline{T}_3T_2\overline{T}_1)\bigr]\label{e:2.4}
\end{eqnarray}
for any tangent vectors in $T_{p_0}G(n, n+m;\C)$ as in
(\ref{e:2.2}), $T_i$, $i=1,2,3,4$, which are identified with
complex matrices $(T^{kl}_i)\in\C^{n\times m}$, $i=1,2,3,4$. It
follows that
 the sectional
curvature sits between $0$ and $4$, and that the
 holomorphic sectional curvature of $G(n, n+m;\C)$
at the point $p_0\in U_{1\cdots n}$ in the direction $T$ is given by
\begin{equation}\label{e:2.5}
H(0, T)=\frac{2{\rm Tr}(T\overline{T}'T\overline{T}')}{[{\rm
Tr}(T\overline{T}')]^2}\in [4/\min(n,m), 4]
\end{equation}
(cf. \cite[(2.11)]{Lu1} and  \cite[page 77]{Wo2}).

\begin{proposition}\label{prop:2.1}
For the metric $h$ in (\ref{e:2.1}) let $R$ be the Riemannian
curvature tensor $R$ of the Riemannian metric $g=Re(h)$ (extended to
$TG(n, n+m;\C)\otimes_{\R}\C$ in a $\C$-linear way). For $1\le
i,j,k,h\le n$ and $1\le\alpha,\beta,\gamma,\delta\le m$ let
\begin{eqnarray*}
R_{i\alpha, \overline{j\beta}, k\gamma,
\overline{h\delta}}&=&R\Bigl(\frac{\partial}{\partial
Z^{i\alpha}}\Bigm|_0 , \frac{\partial}{\partial
\overline{Z}^{j\beta}}\Bigm|_0 , \frac{\partial}{\partial
Z^{k\gamma}}\Bigm|_0 , \frac{\partial}{\partial
\overline{Z}^{h\delta}}\Bigm|_0 \Bigr)\\
&=&g\biggl(R\Bigl(\frac{\partial}{\partial Z^{i\alpha}}\Bigm|_0 ,
\frac{\partial}{\partial \overline{Z}^{j\beta}}\Bigm|_0\Bigr)
\frac{\partial}{\partial \overline{Z}^{h\delta}}\Bigm|_0 ,
\frac{\partial}{\partial {Z}^{k\gamma}}\Bigm|_0 \biggr)
\end{eqnarray*}
and others be defined similarly. Then
$$
R_{i\alpha, \overline{j\beta}, k\gamma,
\overline{h\delta}}=R_{i\alpha, \overline{h\delta}, k\gamma,
\overline{j\beta}} =-R_{i\alpha, \overline{h\delta},
\overline{j\beta},k\gamma}=\frac{1}{2}(-\delta_{ij}\delta_{kh}\delta_{\alpha\delta}\delta_{\beta\gamma}-
\delta_{ih}\delta_{k j}\delta_{\alpha\beta}\delta_{\gamma\delta})
$$
for all $1\le i, j, k, l\le n$ and
$1\le\alpha,\beta,\gamma,\delta\le m$. These and their complex
conjugates are all component types different from zero.
\end{proposition}

\noindent{\bf Proof}. By (\ref{e:2.1}), for $h=2\partial\bar
\partial \Phi(Z)$, where $\Phi(Z)=\frac{1}{2}\ln\det(I+ Z\overline{Z}')$, from
the well-known formula $\det A=\exp\bigl\{{\rm Tr}\ln A\bigr\}$ we
have
\begin{eqnarray*}
2\Phi(Z)&=&{\rm Tr}\ln(I+ Z\overline{Z}')={\rm
Tr}\biggl(\sum^\infty_{q=1}\frac{(-1)^{q+1}}{q} (Z\overline{Z}')^q
\biggr)\nonumber\\
&=&\sum_{i,\alpha}|Z^{i\alpha}|^2-\frac{1}{2}\sum_{i,j,\alpha,\beta}\overline{
Z}^{i\alpha}Z^{i\beta}\overline{Z}^{j\beta}Z^{j\alpha}+
\hbox{(higher order terms)}
\end{eqnarray*}
for $\|Z\overline{Z}'\|<1$. (See also \cite[page 493]{CaVe}). From
this and the arguments on the pages 155-159 of \cite{KoNo}, it
follows that the curvature tensor at $Z=0$ is given by
$$
R_{i\alpha, \overline{j\beta}, k\gamma,
\overline{h\delta}}=\frac{\partial^4\Phi}{\partial
Z^{i\alpha}\partial \overline{Z}^{j\beta}\partial
Z^{k\gamma}\partial\overline {Z}^{h\delta}
}\Bigm|_{Z=0}=\frac{1}{2}(-\delta_{ij}\delta_{kh}\delta_{\alpha\delta}\delta_{\beta\gamma}-
\delta_{ih}\delta_{k j}\delta_{\alpha\beta}\delta_{\gamma\delta})
$$
for all $1\le i, j, k, l\le n$ and
$1\le\alpha,\beta,\gamma,\delta\le m$. Moreover, from the Bianchi
identity and the fact that the curvature tensor $R$ of K\"ahler
manifold is of type $(2,2)$ it is not hard to derive that
$$
R_{i\alpha, \overline{j\beta}, k\gamma,
\overline{h\delta}}=R_{i\alpha, \overline{h\delta}, k\gamma,
\overline{j\beta}} =-R_{i\alpha, \overline{h\delta},
\overline{j\beta},k\gamma}
$$
for all $1\le i, j, k, l\le n$ and
$1\le\alpha,\beta,\gamma,\delta\le m$. These and their complex
conjugates are all component types different from zero.
$\Box$\vspace{2mm}

Consider the non-degenerate alternating (resp. symmetric) bilinear
form $J_n$ (resp. $\Sigma_n$) on $\C^{2n}$ represented by the matrix
$J_n=\left(
     \begin{array}{ccccc}
       0 & I_n  \\
       -I_n & 0
     \end{array}
   \right)$ (resp. $\Sigma_n=\left(
     \begin{array}{ccccc}
       0 & I_n  \\
       I_n & 0
     \end{array}
   \right)$) in
Euclidean coordinate on $\C^{2n}$. Then $G^{\rm II}(n, 2n)$ (resp.
$G^{\rm III}(n, 2n)$) is the set of complex $n$-planes $V$ in
$\C^{2n}$ such that $\Sigma_n|_V=0$ (resp. $J_n|_V=0$). Following
\cite{Lu1, Lu2, Lu3} we have also more convenient matrix
representations of the Hermitian symmetric spaces $G^{\rm II}(n,
2n)$, $G^{\rm III}(n, 2n)$ and $G^{\rm IV}(1, n+1)$ corresponding
with the matrix definition of the complex Grassmann manifolds:
\begin{eqnarray*}
 &&G^{\rm II}(n, 2n)=\biggl\{[A]\in G(n, 2n)\,\biggm|\, \exists A\in
[A] \;\hbox{such that}\; A \left(
     \begin{array}{ccccc}
       0 & I_n  \\
       I_n & 0
     \end{array}
   \right)A'=0\,\biggr\},\\
&&G^{\rm III}(n, 2n)=\biggl\{[A]\in G(n, 2n)\,\biggm|\, \exists A\in
[A] \;\hbox{such that}\; A \left(
     \begin{array}{ccccc}
       0 & I_n  \\
       -I_n & 0
     \end{array}
   \right)A'=0\,\biggr\},\\
&&G^{\rm IV}(1, n+1)=\Bigl\{[(z_1,\cdots,
z_{n+2})]\in\CP^{n+1}\,\bigm|\,
\sum^n_{j=1}z_j^2-z_{n+1}^2-z_{n+2}^2=0\Bigr\}.
\end{eqnarray*}
They are  the compact duals (or extended spaces) of the classical
domains  ${\rm D}^{\rm II}_n$, ${\rm D}^{\rm III}_n$ and ${\rm
D}^{\rm IV}_n$ as $G(n, n+m;\C)$ is that of ${\rm D}^{\rm I}_{n,m}$
(cf. \cite{Lu1, Lu2, Lu3}).


Let $h_{\rm I}$ be the canonical K\"ahler metric on $G^{\rm I}(n,
2n)$, which in the coordinate chart $U_{\alpha_1\cdots\alpha_n}$ is
given by $\partial\bar
\partial\ln\det(I+ Z\overline{Z}')$. $h_{\rm I}$ induces a K\"ahler
metric $h_{\rm II}$ on $G^{\rm II}(n, 2n)$ which in the induced
coordinate system
\begin{equation}\label{e:2.6}
 G^{\rm II}(n, 2n)\cap
U_{\alpha_1\cdots \alpha_n}\ni [A]\mapsto
\Bigl(Z^{kl}([A])\Bigr)_{k<l}
\end{equation}
is given by
\begin{equation}\label{e:2.7}
h_{\rm II}=\partial\bar \partial\ln\det(I- Z\overline{Z})
\end{equation}
 with $Z\in\C^{n\times n}$ and $Z=-Z'$; moreover $h_{\rm I}$ induces a
K\"ahler metric $h_{\rm III}$ on $G^{\rm III}(n, 2n)$ which in the
induced coordinate system
\begin{equation}\label{e:2.8}
 G^{\rm III}(n, 2n)\cap
U_{\alpha_1\cdots \alpha_n}\ni [A]\mapsto
\Bigl(Z^{kl}([A])\Bigr)_{k\le l}
\end{equation}
is given by
\begin{equation}\label{e:2.9}
h_{\rm III}=\partial\bar \partial\ln\det(I+ Z\overline{Z})
\end{equation}
 with $Z\in\C^{n\times n}$ and $Z=Z'$.

Let $h_{\rm FS}$ be the Fubini-Study metric on $\CP^{n+1}$, which is
given by
\begin{equation}\label{e:2.10}
h_{\rm FS}=\partial\bar\partial\ln(1+ |\xi_1|^2+\cdots+
|\xi_{n+1}|^2)
\end{equation}
with $\xi_k=\xi_k([z])=z_k/z_{n+2}$, $k=1,\cdots, n+1$, $[z]\in
U_{n+2}=\{[z_1,\cdots, z_{n+2}]\in\CP^{n+1}\,|\, z_{n+2}\ne 0\}$.
Then $G^{\rm IV}(1, n+1)$ is a K\"ahler submanifold of $\CP^{n+1}$
with the induced K\"ahler metric
\begin{equation}\label{e:2.11}
h_{\rm IV}=\partial\bar\partial\ln(1+ |\xi_1|^2+\cdots+ |\xi_n|^2+
|1-\xi^2_1-\cdots-\xi_n^2|)
\end{equation}
on $G^{\rm IV}(1, n+1)\cap U_{n+2}$ from $h_{\rm FS}$. If ${\rm
Im}\xi_{n+1}\ne 0$, in the new coordinate chart on $G^{\rm IV}(1,
n+1)$,
$$
(\xi_1,\cdots,\xi_n)\mapsto Z=(Z_1,\cdots,
Z_n)=\Bigl(\frac{\xi_1}{\xi_{n+1}+ i},\cdots,
\frac{\xi_n}{\xi_{n+1}+ i}\Bigr),
$$
the metric $h_{\rm IV}$ has the following expression (cf.\cite{Lu1})
\begin{equation}\label{e:2.12}
h_{\rm IV}=\partial\bar\partial\ln(1+ |ZZ'|^2+ 2Z\bar Z').
\end{equation}

All irreducible symmetric spaces of compact type have positive
holomorphic sectional curvatures (cf. \cite{Bo, CaVe, Lu1}). As in
(\ref{e:2.5})  the following explicit expressions come from
\cite{Lu1} too. Under the above coordinate charts their holomorphic
sectional curvatures are
$$
 H_{\rm II}(Z, T)= 2\frac{{\rm
Tr}\bigl[(I- Z\bar Z)^{-1}T(I- \bar ZZ)^{-1}\bar T (I- Z\bar
Z)^{-1}T(I-\bar ZZ)^{-1}\bar T\bigr]}{\bigl\{{\rm Tr}\bigl[(I- Z\bar
Z)^{-1}T(I- \bar ZZ)^{-1}\bar T\bigr]\}^2}
$$
where the tangent vector $ T=\sum_{k,l}{\rm
Re}(T^{kl})\frac{\partial}{\partial X^{kl}} + \sum_{k,l}{\rm
Im}(T^{kl})\frac{\partial}{\partial Y^{kl}} $ with skew-symmetric
complex matrix $T=(T^{kl})\in \C^{n\times n}$,
$$
H_{\rm III}(Z, T)= 2\frac{{\rm
Tr}\bigl[(I+ Z\bar Z)^{-1}T(I+ \bar ZZ)^{-1}\bar T (I+ Z\bar
Z)^{-1}T(I+\bar ZZ)^{-1}\bar T\bigr]}{\bigl\{{\rm Tr}\bigl[(I+ Z\bar
Z)^{-1}T(I+ \bar ZZ)^{-1}\bar T\bigr]\}^2}
$$
where the tangent vector $ T=\sum_{k,l}{\rm
Re}(T^{kl})\frac{\partial}{\partial X^{kl}} + \sum_{k,l}{\rm
Im}(T^{kl})\frac{\partial}{\partial Y^{kl}} $ with symmetric complex
matrix $T=(T^{kl})\in \C^{n\times n}$,
$$
H_{\rm IV}(0, T)= \frac{2(T\bar T')^2-|TT'|^2}{(T\bar T')^2}
$$
where the tangent vector $ T=\sum_{k}{\rm
Re}(T^{k})\frac{\partial}{\partial X^{k}} + \sum_{k}{\rm
Im}(T^{k})\frac{\partial}{\partial Y^{k}} $ with complex vector
$T=(T^{1},\cdots, T^n)\in \C^{n}$.

Let $R^{\rm I}$ denote the curvature tensor of the metric $h_{\rm
I}=\partial\bar
\partial \ln\det(I+ Z\overline{Z}')$ on $G^{\rm I}(n, 2n)$. By
Proposition~\ref{prop:2.1}, at $Z=0$ we have
\begin{eqnarray}\label{e:2.13}
R^{\rm I}_{i\alpha, \overline{j\beta}, k\gamma,
\overline{h\delta}}&=&R^{\rm I}_{i\alpha, \overline{h\delta},
k\gamma, \overline{j\beta}} =-R^{\rm I}_{i\alpha,
\overline{h\delta},
\overline{j\beta},k\gamma}\nonumber\\
&=&\frac{1}{2}(-\delta_{ij}\delta_{kh}\delta_{\alpha\delta}\delta_{\beta\gamma}-
\delta_{ih}\delta_{k j}\delta_{\alpha\beta}\delta_{\gamma\delta})
\end{eqnarray}
for all $1\le i, j, k, l, \alpha,\beta,\gamma,\delta\le n$. These
and their complex conjugates are all component types different from
zero.

 Denote the curvature tensors of
$(G^{\rm II}(n,2n), h_{\rm II})$ and $(G^{\rm III}(n, 2n), h_{\rm
III})$ by $R^{\rm II}$ and $R^{\rm III}$, respectively. Note that at
$Z=0$ the local coordinate systems $(U_{1\cdots n}, Z)$ on $G(n,
2n;\C)$ and (\ref{e:2.6})-(\ref{e:2.8})  are normal coordinates (or
complex geodesic coordinates) for the metrics $h_{\rm I}$, $h_{\rm
II}$ and $h_{\rm III}$. (In fact, $(G^{\rm II}(n,2n), h_{\rm II})$
and $(G^{\rm III}(n, 2n), h_{\rm III})$ are totally geodesic
submanifolds of $(G^{\rm I}(n, 2n), h_{\rm I})$, see the claim on
the page 136 of \cite{Mok} and the proof of Lemma 1 on the page 85
of \cite{Mok}).  By (\ref{e:2.13}) we have

\begin{proposition}\label{prop:2.2}
At $Z=0$ the curvature tensors $R^{\rm II}$ and $R^{\rm III}$ are
the restrictions of $R^{\rm I}$, that is,
$$
R^{\rm II}_{i\alpha, \overline{j\beta}, k\gamma,
\overline{h\delta}}=R^{\rm II}_{i\alpha, \overline{h\delta},
k\gamma, \overline{j\beta}} =-R^{\rm II}_{i\alpha,
\overline{h\delta},
\overline{j\beta},k\gamma}=\frac{1}{2}(-\delta_{ij}\delta_{kh}\delta_{\alpha\delta}\delta_{\beta\gamma}-
\delta_{ih}\delta_{k j}\delta_{\alpha\beta}\delta_{\gamma\delta})
$$
for all $1\le i<\alpha\le n, 1\le j<\beta\le n, 1\le k<\gamma\le n,
l\le l<\delta\le n$, and
$$
R^{\rm III}_{i\alpha, \overline{j\beta}, k\gamma,
\overline{h\delta}}=R^{\rm III}_{i\alpha, \overline{h\delta},
k\gamma, \overline{j\beta}} =-R^{\rm III}_{i\alpha,
\overline{h\delta},
\overline{j\beta},k\gamma}=\frac{1}{2}(-\delta_{ij}\delta_{kh}\delta_{\alpha\delta}\delta_{\beta\gamma}-
\delta_{ih}\delta_{k j}\delta_{\alpha\beta}\delta_{\gamma\delta})
$$
for all $1\le i\le\alpha\le n, 1\le j\le\beta\le n, 1\le
k\le\gamma\le n, 1\le l\le\delta\le n$.
\end{proposition}

Now we consider  $(G^{\rm IV}(1, n+1), h_{\rm IV})$. By
(\ref{e:2.12}) the K\"ahler potential function
$\Phi(Z)=\frac{1}{2}\ln(1+ |ZZ'|^2+ 2Z\bar Z')$ has the following
power series expansion
\begin{eqnarray*}
&&\frac{1}{2}\ln(1+ |ZZ'|^2+ 2Z\bar Z')=\frac{1}{2}\ln\bigl(1+
2\sum_{k}|z_k|^2+
|\sum^n_{k=1}z_k^2|^2\bigr)\\
&=&\sum^n_{k=1}|Z_k|^2+ \frac{1}{2}|\sum^n_{k=1}Z_k^2|^2-
(\sum^n_{k=1}|Z_k|^2)^2+ \hbox{higher order terms}
\end{eqnarray*}
near $Z=0$. Since the coordinates $Z_k$ ($1\le k\le n$) are normal
coordinates, the curvature tensor at $Z=0$ is given by
$$
R^{\rm IV}_{i\bar jk\bar l}=\frac{\partial^4\Phi}{\partial
Z_i\partial\bar Z_j\partial Z_k\partial\bar
Z_l}\Bigm|_{Z=0}=2(\delta_{ik}\delta_{jl}-\delta_{ij}\delta_{kl}-\delta_{il}\delta_{jk})
$$
for all $1\le i, j, k, l\le n$. In particular we get
\begin{equation}\label{e:2.14}
R^{\rm IV}_{i\bar ii\bar i}=-2\;\forall i   \quad\hbox{and}\quad
R^{\rm IV}_{i\bar ji\bar j}=2\;\forall i\ne j.
\end{equation}

\subsection{An expected local coordinate chart}

 Let $J$ be the standard complex structure on $G(n, n+m;\C)$. For $p\in
G(n, n+m;\C)$, recall that by  $\{a_{ij}, b_{ij},\;i=1,\cdots,
n,\;j=1,\cdots,m\}$ being a {\it unitary base} of $(T_pG(n, n+m),
J_p, g_p)$ we mean
$$
a_{ij}, b_{ij}=J_pa_{ij}\in T_pG(n, n+m;\C),\;i=1,\cdots,
n,\;j=1,\cdots,m,
$$
is a unit orthogonal base of $(T_pG(n, n+m;\C),  g_p)$. To our
knowledge the following result seems to be new. It is key for us
completing the proofs of Theorems~\ref{th:1.1},~\ref{th:1.2}.

\begin{proposition}\label{prop:2.3}
For any $p\in G(n, n+m;\C)$ and a  unitary base of $(T_pG(n,
n+m;\C), J_p, g_p)$,
$$
a_{ij}, b_{ij}:=J_pa_{ij}\in T_pG(n, n+m;\C),\;i=1,\cdots,
n,\;j=1,\cdots,m,
$$
there exists a local chart around $p$ on $G(n, n+m;\C)$,
\begin{equation}\label{e:2.15}
{\cal U}\ni q\to Z(q)=X(q)+ iY(q)\in \C^{n\times m}
\end{equation}
satisfying $Z(p)=0$, such that
\begin{description}
\item[(i)] In this chart the metric $h$ and $g=Re(h)$ are given by
(\ref{e:2.1}) and (\ref{e:2.3}), respectively;

\item[(ii)] $a_{ij}=\frac{\partial}{\partial X^{ij}}|_p,\; b_{ij}=\frac{\partial}{\partial
Y^{ij}}|_p   ,\;i=1,\cdots, n,\;j=1,\cdots,m$.
\end{description}
\end{proposition}

\noindent{\bf Proof.} Since the isometry group of the K\"ahler
manifold $(G(n, n+m;\C), h)$, ${\rm I}(G(n, n+m;\C), h)=SU(n+m)$,
acts transitively on $(G(n, n+m;\C), h)$, for any $p\in G(n,
n+m;\C)$ there exists a $\tau\in {\rm I}(G(n, n+m;\C), h)$ such
that $\tau(p_0)=p$. Clearly, we get a coordinate chart around $p$
on $G(n, n+m;\C)$,
\begin{equation}\label{e:2.16}
W=U+ iV: \tau(U_{1\cdots n})\to\C^{n\times m},\;q\mapsto
Z(\tau^{-1}(q)).
\end{equation}
Since $\tau$ is a K\"ahler isometry, using (\ref{e:2.1}) one easily
shows that the metric $h$ in this chart is given by
$$
h={\rm Tr}[(I+ W\overline{W}')^{-1}dW(I+
\overline{W}'W)^{-1}\overline{dW}'].
$$
It follows that the Riemannian metric $g=Re(h)$ is given by
$$
g_W(T_1, T_2)=Re {\rm Tr}[(I+
W\overline{W}')^{-1}T_1(I+\overline{W}'W)^{-1}\overline{T_2}']
$$
for real tangent vectors $T_1, T_2$ at $W\in\tau(U_{1\cdots n})$,
\begin{eqnarray*}
&&T_1=\sum_{k,l}{\rm Re}(T^{kl}_1)\frac{\partial}{\partial U^{kl}} +
\sum_{k,l}{\rm Im}(T^{kl}_1)\frac{\partial}{\partial
V^{kl}},\\
&&T_2=\sum_{k,l}{\rm Re}(T^{kl}_2)\frac{\partial}{\partial U^{kl}} +
\sum_{k,l}{\rm Im}(T^{kl}_2)\frac{\partial}{\partial V^{kl}},
\end{eqnarray*}
which are identified with complex matrices $(T_1^{kl}),
(T_2^{kl})\in\C^{n\times m}$, respectively.

  Define vectors
\begin{eqnarray*}
&&\overrightarrow{a}=(a_{11}, a_{12},\cdots, a_{1m}, a_{21},\cdots,
a_{2m},\cdots, a_{n1},\cdots, a_{nm}),\\
&&\overrightarrow{b}=(b_{11}, b_{12},\cdots, b_{1m}, b_{21},\cdots,
b_{2m},\cdots, b_{n1},\cdots, b_{nm}),\\
&&\overrightarrow{\frac{\partial}{\partial
U}|_p}=(\frac{\partial}{\partial
U^{11}}|_p,\cdots,\frac{\partial}{\partial U^{1m}}|_p,
\frac{\partial}{\partial U^{21}}|_p,\cdots,\frac{\partial}{\partial
U^{2m}}|_p, \cdots, \frac{\partial}{\partial
U^{n1}}|_p,\cdots,\frac{\partial}{\partial U^{nm}}|_p),\\
&&\overrightarrow{\frac{\partial}{\partial
V}|_p}=(\frac{\partial}{\partial
V^{11}}|_p,\cdots,\frac{\partial}{\partial V^{1m}}|_p,
\frac{\partial}{\partial V^{21}}|_p,\cdots,\frac{\partial}{\partial
V^{2m}}|_p, \cdots, \frac{\partial}{\partial
V^{n1}}|_p,\cdots,\frac{\partial}{\partial V^{nm}}|_p).
\end{eqnarray*}

Since
$$
\left\{\frac{\partial}{\partial U^{ij}}|_p, \frac{\partial}{\partial
V^{ij}}|_p   ,\;i=1,\cdots, n,\;j=1,\cdots,m\right\}
$$
is  a  unitary base of $(T_pG(n, n+m;\C), J_p, g_p)$,  there exists
a unique real matrix $\Theta$ such that
\begin{equation}\label{e:2.17}
(\overrightarrow{a},
\overrightarrow{b})=\left(\overrightarrow{\frac{\partial}{\partial
U}|_p}, \overrightarrow{\frac{\partial}{\partial
V}|_p}\right)\Theta.
\end{equation}
The matrix $\Theta$ must have form $\left(
     \begin{array}{ccccc}
       {\cal A} & {\cal B}  \\
       -{\cal B} & {\cal A}
     \end{array}
   \right)$, where ${\cal A}, {\cal B}\in\R^{nm\times nm}$ is such
   that ${\cal A}+ i{\cal B}$ is a unitary matrix (which is equivalent to
$$
{\cal B}'{\cal A}=({\cal A}'{\cal B})'={\cal A}'{\cal
B}\quad\hbox{and}\quad {\cal A}'{\cal A}+ {\cal B}'{\cal
B}=I_{nm\times nm}.)
$$
Note that (\ref{e:2.17}) is
   equivalent to
\begin{equation}\label{e:2.18}
\overrightarrow{a}+
i\overrightarrow{b}=\left(\overrightarrow{\frac{\partial}{\partial
U}|_p}+ i\overrightarrow{\frac{\partial}{\partial
V}|_p}\right)({\cal A}+ i{\cal B}).
\end{equation}
Recall that the tensor product or Kronecker product of matrices
$A=(a_{ij})\in\C^{n\times m}$ and $B=(b_{ij})\in\C^{p\times q}$ is
a $(np\times mq)$-matrix given by
\begin{eqnarray*}
A\otimes B=[a_{ij}B]^{n,m}_{i,j=1}=\left(
  \begin{array}{ccc}
    a_{11}B & \cdots & a_{1m}B \\
    \cdots & \cdots & \cdots \\
    a_{n1}B & \cdots & a_{nm}B \\
  \end{array}
\right).
\end{eqnarray*}
Define matrices  ${\bf a}=(a_{ij})$, ${\bf b}=(b_{ij})$ and
$\frac{\partial}{\partial U}|_p=(\frac{\partial}{\partial
U^{ij}}|_p)$, $\frac{\partial}{\partial
V}|_p=(\frac{\partial}{\partial V^{ij}}|_p)$. It follows from
(\ref{e:2.18}) that there exist unitary matrices ${\cal
R}\in\C^{n\times n}$ and ${\cal S}\in\C^{m\times m}$ such that
\begin{equation}\label{e:2.19}
{\cal A}+ i{\cal B}={\cal R}'\otimes{\cal S}\quad\hbox{and}\quad
{\bf a}+ i{\bf b}={\cal R}(\frac{\partial}{\partial U}|_p+
i\frac{\partial}{\partial V}|_p){\cal S}.
\end{equation}
Let ${\cal R}=R_1+ iR_2$ with $R_1, R_2\in\R^{n\times n}$, and
${\cal S}=S_1+ iS_2$ with $S_1, S_2\in\R^{m\times m}$. Then
$$
\left.\begin{array}{ll}
 & (R_1'R_2)'=R_1'R_2\quad\hbox{and}\quad R_1'R_1+ R_2'R_2=I_{n\times n},\\
 & (S_1'S_2)'=S_1'S_2\quad\hbox{and}\quad S_1'S_1+ S_2'S_2=I_{m\times
 m}.
\end{array}\right\}
$$
Moreover, the first equality in (\ref{e:2.19}) implies
$$
{\cal A}=R_1'\otimes S_1-R'_2\otimes S_2\quad\hbox{and}\quad {\cal
B}=R_2'\otimes S_1+ R_1'\otimes S_2.
$$

From the local chart $(\tau(U_{1\cdots n}), W)$ in (\ref{e:2.16}),
we define a new chart
\begin{equation}\label{e:2.20}
{\cal U}\to\C^{n\times m},\;q\mapsto G(q)=E(q)+ iF(q):= {\cal
R}^{-1}W(q){\cal S}^{-1}.
\end{equation}
Then $G(p)=W(p)=0$. Define vectors
\begin{eqnarray*}
&&\overrightarrow{W}=(W^{11}, W^{12},\cdots, W^{1m}, Z^{21},\cdots,
W^{2m},\cdots, W^{n1},\cdots, W^{nm}),\\
&&\overrightarrow{G}=(G^{11}, G^{12},\cdots, G^{1m}, G^{21},\cdots,
G^{2m},\cdots, G^{n1},\cdots, G^{nm}).
\end{eqnarray*}
By \cite[page 364, (6)]{Lu2} we get
\begin{equation}\label{e:2.21}
\frac{\partial G}{\partial W}=\frac{\partial
\overrightarrow{G}}{\partial \overrightarrow{W}}=({\cal
R}^{-1})'\otimes{\cal S}^{-1}=({\cal R}'\otimes{\cal S})^{-1}=({\cal
A}+ i{\cal B})^{-1}.
\end{equation}
Writing $G=\Phi(W)$ and
\begin{eqnarray*}
&&\overrightarrow{\frac{\partial}{\partial
W}|_p}=(\frac{\partial}{\partial
W^{11}}|_p,\cdots,\frac{\partial}{\partial W^{1m}}|_p,
\frac{\partial}{\partial W^{21}}|_p,\cdots,\frac{\partial}{\partial
W^{2m}}|_p, \cdots, \frac{\partial}{\partial
W^{n1}}|_p,\cdots,\frac{\partial}{\partial W^{nm}}|_p),\\
&&\Phi_\ast\Bigl(\overrightarrow{\frac{\partial}{\partial
W}|_p}\Bigr)=\Bigr(\Phi_\ast(\frac{\partial}{\partial
W^{11}}|_p),\cdots, \Phi_\ast(\frac{\partial}{\partial W^{1m}}|_p),
\Phi_\ast(\frac{\partial}{\partial W^{21}}|_p),\cdots,\\
&&\hspace{45mm}\Phi_\ast(\frac{\partial}{\partial W^{2m}}|_p),
\cdots, \Phi_\ast(\frac{\partial}{\partial
W^{n1}}|_p),\cdots,\Phi_\ast(\frac{\partial}{\partial
W^{nm}}|_p)\Bigr),
\end{eqnarray*}
since $\overrightarrow{\frac{\partial}{\partial U}|_p}+
i\overrightarrow{\frac{\partial}{\partial
V}|_p}=\overrightarrow{\frac{\partial}{\partial W}|_p}$, by
(\ref{e:2.18}) and (\ref{e:2.21}) we get
\begin{eqnarray*}
\overrightarrow{\frac{\partial}{\partial
G}|_p}&=&\Phi_\ast\Bigl(\overrightarrow{\frac{\partial}{\partial
W}|_p}\Bigr)=\overrightarrow{\frac{\partial}{\partial
W}|_p}\frac{\partial \overrightarrow{W}}{\partial
\overrightarrow{G}}\nonumber\\
&=&\overrightarrow{\frac{\partial}{\partial
W}|_p}\left(\frac{\partial \overrightarrow{G}}{\partial
\overrightarrow{W}}\right)^{-1}=\overrightarrow{\frac{\partial}{\partial
W}|_p}({\cal A}+ i{\cal B})=\overrightarrow{a}+ i\overrightarrow{b}.
\end{eqnarray*}
That is, the coordinate chart in (\ref{e:2.20}), ${\cal
U}\to\C^{n\times m},\;q\mapsto G(q)$, satisfies
$$
a_{ij}=\frac{\partial}{\partial E^{ij}}|_p,\quad
b_{ij}=\frac{\partial}{\partial F^{ij}}|_p   ,\;i=1,\cdots,
n,\;j=1,\cdots,m.
$$

 It remains to prove that the transformation
$$
\C^{n\times m}\to \C^{n\times m},\;W\mapsto G=\Phi(W)
$$
preserves the K\"ahler metric
$$
ds^2={\rm Tr}[(I+ W\overline{W}')^{-1}dW(I+
\overline{W}'W)^{-1}\overline{dW}']
$$
on $\C^{n\times m}$. In fact,
since
\begin{eqnarray*}
(I+ G\overline{G}')^{-1}dG &=&(I+ {\cal R}^{-1}W{\cal S}^{-1}
\overline{{\cal R}^{-1}W{\cal S}^{-1}}' ){\cal R}^{-1}dW{\cal
S}^{-1}\\
&=&(I+ {\cal R}^{-1}W\overline{W}' \overline{{\cal R}^{-1}}' ){\cal
R}^{-1}dW{\cal
S}^{-1}\\
&=&({\cal R}^{-1}\overline{{\cal R}^{-1}}'+ {\cal
R}^{-1}W\overline{W}' \overline{{\cal R}^{-1}}' ){\cal
R}^{-1}dW{\cal S}^{-1}\\
&=&{\cal R}^{-1}(I+ W\overline{W}' )dW{\cal S}^{-1},\\
(I+ \overline{G}'G)^{-1}d\overline{G}' &=&(I+ \overline{{\cal
R}^{-1}W{\cal S}^{-1}}'{\cal R}^{-1}W{\cal S}^{-1} )\overline{{\cal
R}^{-1}dW{\cal S}^{-1}}'\\
 &=&(I+ \overline{{\cal S}^{-1}}'\overline{W}'\overline{{\cal
R}^{-1}}'{\cal R}^{-1}W{\cal S}^{-1} )\overline{{\cal
S}^{-1}}'\overline{dW}'\overline{{\cal R}^{-1}}'\\
&=&\overline{{\cal S}^{-1}}'(I+ \overline{W}'W
)\overline{dW}'\overline{{\cal R}^{-1}}'
\end{eqnarray*}
we get
\begin{eqnarray*}
&&{\rm Tr}[(I+
G\overline{G}')^{-1}dG(I+ \overline{G}'G)^{-1}\overline{dG}']\\
&=&{\rm Tr}\bigl[(I+ \Phi(W)\overline{\Phi(W)}')^{-1}d\Phi(W)(I+
\Phi(W)\overline{\Phi(W)}'\Phi(W))^{-1}\overline{d\Phi(W)}'\bigr]\\
&=&{\rm Tr}[(I+ W\overline{W}')^{-1}dW(I+
\overline{W}'W)^{-1}\overline{dW}'].
\end{eqnarray*}
Hence the coordinate chart in (\ref{e:2.20}) satisfies the desired
requirements. $\Box$\vspace{2mm}

\begin{corollary}\label{cor:2.4}
For any $p, q\in G(n, n+m;\C)$, let
\begin{eqnarray*}
&&\{a_{ij}, b_{ij}:=J_pa_{ij},\;i=1,\cdots,
n,\;j=1,\cdots,m\}\quad\hbox{and}\\
&&\{a'_{ij}, b'_{ij}:=J_qa'_{ij},\;i=1,\cdots, n,\;j=1,\cdots,m\}
\end{eqnarray*}
be unitary bases of $(T_pG(n, n+m;\C), J_p, g_p)$ and $(T_qG(n,
n+m;\C), J_q, g_q)$, respectively. Consider the sequence
$u_1,\cdots, u_{2nm}$ whose all odd (resp. even) terms are given by
\begin{eqnarray*}
&&a_{11}, a_{12},\cdots, a_{1m}, a_{21},\cdots,
a_{2m},\cdots, a_{n1},\cdots, a_{nm},\\
\hbox{(resp.} &&b_{11}, b_{12},\cdots, b_{1m}, b_{21},\cdots,
b_{2m},\cdots, b_{n1},\cdots, b_{nm}. \hbox{)}
\end{eqnarray*}
Similarly let the sequence $u'_1,\cdots, u'_{2nm}$ be given by
$\{a'_{ij}, b'_{ij}:=J_qa'_{ij},\;i=1,\cdots, n,\;j=1,\cdots,m\}$.
Then the curvature tensor $R$ of $(G(n, n+m;\C), g)$ satisfies
\begin{equation}\label{e:2.22}
R_p(u_\alpha, u_\beta, u_\gamma, u_\delta)=R_q(u'_\alpha, u'_\beta,
u'_\gamma, u'_\delta)
\end{equation}
for any $\alpha, \beta, \gamma,\delta\in\{1,\cdots, 2nm\}$.
\end{corollary}

\noindent{\bf Proof}. This can be directly derived from
Propositions~\ref{prop:2.1},~\ref{prop:2.3}. We here give another
proof of it with (\ref{e:2.4}). Let $({\cal U}, Z)$ be a local
chart around $p$ as in (\ref{e:2.15}). Then
$a_{ij}=\frac{\partial}{\partial X^{ij}}|_p,\;
b_{ij}=\frac{\partial}{\partial Y^{ij}}|_p ,\;i=1,\cdots,
n,\;j=1,\cdots,m$. Let Let $({\cal V}, W=U+ \sqrt{-1}V)$ be a
local chart around $q$ as in (\ref{e:2.15}). Then
$a'_{ij}=\frac{\partial}{\partial U^{ij}}|_q,\;
b'_{ij}=\frac{\partial}{\partial V^{ij}}|_q   ,\;i=1,\cdots,
n,\;j=1,\cdots,m$. Note that according to the above correspondence
the tangent vectors $\frac{\partial}{\partial X^{kl}}|_p$ and
$\frac{\partial}{\partial Y^{st}}|_p$ have matrices
representations
\begin{eqnarray}\label{e:2.23}
S_{(k,l)}=\left(
      \begin{array}{ccc}
        0 & 0 & 0 \\
        0 & 1_{(k,l)} & 0 \\
        0 & 0 & 0 \\
      \end{array}
    \right)_{n\times m}\quad\hbox{and}\quad
T_{(s,t)}=\left(
      \begin{array}{ccc}
        0 & 0 & 0 \\
        0 & i_{(s,t)} & 0 \\
        0 & 0 & 0 \\
      \end{array}
    \right)_{n\times m}
\end{eqnarray}
respectively, where the first index $(k,l)$ means that $1$ is in
the $k$-th row and $l$-th array of the matrix and similarly for
other indexes in the sequel. Clearly, the tangent vectors
$\frac{\partial}{\partial U^{kl}}|_p$ and
$\frac{\partial}{\partial V^{st}}|_p$ are also represented by
these two matrices. So for any $\alpha\in\{1,\cdots, 2nm\}$ both
$u_\alpha$ and $u_\alpha'$ have the same matrix representations.
The desired conclusions follow from (\ref{e:2.4}) immediately.
$\Box$\vspace{2mm}

This corollary and Proposition~\ref{prop:2.1} immediately lead to

\begin{corollary}\label{cor:2.5}
Let $(M, \omega^M, J^M, g^M)$ be a compact K\"ahler-Einstein
submanifold of $(G(n, n+m;\C), h)$ which is totally geodesic (e.g.
$(G^{\rm II}(n, 2n), h_{\rm II})$ and $(G^{\rm III}(n, 2n), h_{\rm
III})$ are such submanifolds of $(G(n, 2n;\C), h_{\rm I})$). Set
$\dim M=2N$.
 For any $p, q\in M$, let
\begin{eqnarray*}
&&\{a_{2i-1}, a_{2i}:=J^M_pa_{2i-1},\;i=1,\cdots, N\}\quad\hbox{and}\\
&& \{a'_{2i-1}, a'_{2i}:=J^M_qa'_{2i-1},\;i=1,\cdots, N\}
\end{eqnarray*}
be unitary bases of $(T_pM, g^M_p, J^M_p)$ and $(T_qM, g^M_q,
J^M_q)$, respectively. Then the curvature tensor $R^M$ of $(M, g)$
satisfies
$$
R^M_p(a_\alpha, a_\beta, a_\gamma, a_\delta)=R^M_q(a'_\alpha,
a'_\beta, a'_\gamma, a'_\delta)
$$
for any $\alpha, \beta, \gamma,\delta\in\{1,\cdots, \dim M\}$.
\end{corollary}

\noindent{\bf Proof}. Since $(T_pM, g^M_p, J^M_p)$ and $(T_qM,
g^M_q, J^M_q)$ are Hermitian subspaces of $(T_pG(n, n+m;\C), h_p)$
and $(T_qG(n, n+m;\C), h_q)$, respectively, we may extend
$\{a_1,\cdots, a_{2N}\}$ and $\{a'_1,\cdots, a'_{2N}\}$ into unitary
bases
$$
\{a_1,\cdots, a_{2nm}\}\quad\hbox{and}\quad\{a'_1,\cdots, a'_{2nm}\}
$$
 of  $(T_pG(n, n+m;\C), h_p)$ and $(T_qG(n, n+m;\C), h_q)$,
respectively. By the assumptions $(M, \omega^M, J^M, g^M)$ is a
totally geodesic submanifold of $(G(n, n+m;\C), h)$. $R^M$ is equal
to the restriction of $R$ to $M$. Hence the desired conclusion
follows from (\ref{e:2.22}). (Of course it may also be obtained from
Proposition~\ref{prop:2.2} for $(G^{\rm II}(n, 2n), h_{\rm II})$ and
$(G^{\rm III}(n, 2n), h_{\rm III})$). $\Box$\vspace{2mm}

Let $({\cal U}, Z)$ be the local chart around $p$ on $G(n, n+m;\C)$
as in Proposition~\ref{prop:2.3}.

\begin{proposition}\label{prop:2.6}
For any $1\leq k,s,\mu\leq n, 1\leq l,t, \nu\leq m$ it holds that
$$
R\Bigl(\frac{\partial}{\partial X^{kl}}\Bigm|_p,
\frac{\partial}{\partial X^{st}}\Bigm|_p,
 \frac{\partial}{\partial X^{kl}}\Bigm|_p , \frac{\partial}{\partial Y^{\mu\nu}}\Bigm|_p\Bigr)=0,
$$

\begin{displaymath}
R\Bigl(\frac{\partial}{\partial X^{kl}}\Bigm|_p,
\frac{\partial}{\partial X^{st}}\Bigm|_p,
 \frac{\partial}{\partial X^{kl}}\Bigm|_p ,
  \frac{\partial}{\partial X^{\mu\nu}}\Bigm|_p\Bigr) = \left\{ \begin{array}{ll}
1
&\hbox{if $\mu=s\neq k, l=t=\nu $},\\

1
& \hbox{if $\mu=s=k, l\neq t=\nu$},\\

0 & \hbox{otherwise},
\end{array} \right.
\end{displaymath}

\begin{displaymath}
R\Bigl(\frac{\partial}{\partial X^{kl}}\Bigm|_p,
\frac{\partial}{\partial Y^{st}}\Bigm|_p,
 \frac{\partial}{\partial X^{kl}}\Bigm|_p ,
  \frac{\partial}{\partial Y^{\mu\nu}}\Bigm|_p\Bigr) = \left\{ \begin{array}{ll}
1
&\textrm{if $\mu=s\neq k, l=t=\nu $},\\

1
& \textrm{if $\mu=s=k, l\neq t=\nu$},\\

4
& \textrm{if $\mu=s=k, l=t=\nu$},\\

0 & \hbox{otherwise}.
\end{array} \right.
\end{displaymath}
Consequently, for $S_{(k,l)}$ and $T_{(s,t)}$ in (\ref{e:2.23}) we
get the sectional curvatures
\begin{displaymath}
K_p(S_{(k,l)}, T_{(s,t)}):=R\Bigl(\frac{\partial}{\partial
X^{kl}}\Bigm|_p, \frac{\partial}{\partial
Y^{st}}\Bigm|_p,\frac{\partial}{\partial
X^{kl}}\Bigm|_p,\frac{\partial}{\partial Y^{st}}\Bigm|_p\Bigr)
 = \left\{ \begin{array}{ll}
1
&\textrm{if $k=s, l\neq t$,}\\

1
&\textrm{if $k\neq s, l=t$,}\\

4
& \textrm{if $k=s, l=t$,}\\

0 & \textrm{if $k\neq s, l\neq t$,}
\end{array} \right.
\end{displaymath}

\begin{displaymath}
K_p(S_{(k,l)}, S_{(s,t)}):=R\Bigl(\frac{\partial}{\partial
X^{kl}}\Bigm|_p, \frac{\partial}{\partial
X^{st}}\Bigm|_p,\frac{\partial}{\partial
X^{kl}}\Bigm|_p,\frac{\partial}{\partial X^{st}}\Bigm|_p\Bigr)
 = \left\{ \begin{array}{ll}
1
&\textrm{if $k=s, l\neq t$,}\\

1
&\textrm{if $k\neq s, l=t$,}\\

0
& \textrm{if $k=s, l=t$,}\\

0 & \textrm{if $k\neq s, l\neq t$.}
\end{array} \right.
\end{displaymath}

\end{proposition}

\noindent{\bf Proof}. Since the only possible non-vanishing terms of
the curvature components are of the form $R_{i\alpha,
\overline{j\beta}, k\gamma, \overline{h\delta}}$ and those obtained
from the universal symmetries of the curvature tensor, a direct
computation leads to
\begin{eqnarray*}
&&R\Bigl(\frac{\partial}{\partial X^{kl}}\Bigm|_p,
\frac{\partial}{\partial X^{st}}\Bigm|_p,
 \frac{\partial}{\partial X^{kl}}\Bigm|_p,
  \frac{\partial}{\partial X^{\mu\nu}}\Bigm|_p\Bigr)\\
&=&R\Bigl(\frac{\partial}{\partial
Z^{kl}}\Bigm|_p+\frac{\partial}{\partial \overline{Z}^{kl}}\Bigm|_p,
\frac{\partial}{\partial Z^{st}}\Bigm|_p+\frac{\partial}{\partial
\overline{Z}^{st}}\Bigm|_p,
 \frac{\partial}{\partial Z^{kl}}\Bigm|_p+\frac{\partial}{\partial \overline{Z}^{kl}}\Bigm|_p ,
  \frac{\partial}{\partial Z^{\mu\nu}}\Bigm|_p+\frac{\partial}{\partial \overline{Z}^{\mu\nu}}\Bigm|_p\Bigr)\\
&=& R\Bigl(\frac{\partial}{\partial Z^{kl}}\Bigm|_p ,
\frac{\partial}{\partial \overline{Z}^{st}}\Bigm|_p ,
 \frac{\partial}{\partial Z^{kl}}\Bigm|_p  ,
  \frac{\partial}{\partial \overline{Z}^{\mu\nu}}\Bigm|_p\Bigr)
+ R\Bigl(\frac{\partial}{\partial Z^{kl}}\Bigm|_p ,
\frac{\partial}{\partial \overline{Z}^{st}}\Bigm|_p ,
 \frac{\partial}{\partial \overline{Z}^{kl}}\Bigm|_p  ,
  \frac{\partial}{\partial {Z}^{\mu\nu}}\Bigm|_p \Bigr)\\
&&+R\Bigl(\frac{\partial}{\partial \overline{Z}^{kl}}\Bigm|_p ,
\frac{\partial}{\partial {Z}^{st}}\Bigm|_p ,
 \frac{\partial}{\partial Z^{kl}}\Bigm|_p  ,
  \frac{\partial}{\partial \overline{Z}^{\mu\nu}}\Bigm|_p \Bigr)
+R\Bigl(\frac{\partial}{\partial \overline{Z}^{kl}}\Bigm|_p ,
\frac{\partial}{\partial {Z}^{st}}\Bigm|_p ,
 \frac{\partial}{\partial \overline{Z}^{kl}}\Bigm|_p  ,
  \frac{\partial}{\partial {Z}^{\mu\nu}}\Bigm|_p \Bigr) \\
&=&R_{kl,\overline{st},kl,\overline{\mu\nu}}-R_{kl,\overline{st},\mu\nu,\overline{kl}}
-R_{st,\overline{kl},kl,\overline{\mu\nu}}+R_{st,\overline{kl},\mu\nu,\overline{kl}}\\
&=&\frac{1}{2}[-\delta_{ks}\delta_{k\mu}\delta_{l\nu}\delta_{tl}-\delta_{k\mu}\delta_{sk}\delta_{tl}\delta_{l\nu}]
+\frac{1}{2}[\delta_{ks}\delta_{k\mu}\delta_{ll}\delta_{t\nu}+\delta_{kk}\delta_{s\mu}\delta_{tl}\delta_{l\nu}]\\
&&+\frac{1}{2}[\delta_{ks}\delta_{k\mu}\delta_{t\nu}\delta_{ll}+\delta_{s\mu}\delta_{kk}\delta_{tl}\delta_{l\nu}]
+\frac{1}{2}[-\delta_{ks}\delta_{k\mu}\delta_{tl}\delta_{l\nu}-\delta_{sk}\delta_{k\mu}\delta_{tl}\delta_{l\nu}]\\
&=&-2\delta_{sk}\delta_{k\mu}\delta_{tl}\delta_{l\nu}+
\delta_{ks}\delta_{k\mu}\delta_{t\nu}+
\delta_{s\mu}\delta_{tl}\delta_{l\nu},
\end{eqnarray*}
where the final equality comes from Proposition~\ref{prop:2.1}. So
we get
\begin{displaymath}
R\Bigl(\frac{\partial}{\partial X^{kl}}\Bigm|_p,
\frac{\partial}{\partial X^{st}}\Bigm|_p,
 \frac{\partial}{\partial X^{kl}}\Bigm|_p ,
  \frac{\partial}{\partial X^{\mu\nu}}\Bigm|_p\Bigr) = \left\{ \begin{array}{ll}
1
&\hbox{if $\mu=s\neq k, l=t=\nu $},\\

1
& \hbox{if $\mu=s=k, l\neq t=\nu$},\\

0 & \hbox{otherwise}.
\end{array} \right.
\end{displaymath}
Similarly we may obtain
\begin{eqnarray*}
&&R\Bigl(\frac{\partial}{\partial X^{kl}}\Bigm|_p ,
\frac{\partial}{\partial X^{st}}\Bigm|_p , \frac{\partial}{\partial
X^{kl}}\Bigm|_p , \frac{\partial}{\partial Y^{\mu\nu}}\Bigm|_p
\Bigr)=0,\\
&&R\Bigl(\frac{\partial}{\partial X^{kl}}\Bigm|_p ,
\frac{\partial}{\partial Y^{st}}\Bigm|_p , \frac{\partial}{\partial
X^{kl}}\Bigm|_p ,
\frac{\partial}{\partial Y^{\mu\nu}}\Bigm|_p \Bigr)\\
&=& R\Bigl(\frac{\partial}{\partial
Z^{kl}}\Bigm|_p+\frac{\partial}{\partial \overline{Z}^{kl}}\Bigm|_p
, i\frac{\partial}{\partial
Z^{st}}\Bigm|_p-i\frac{\partial}{\partial \overline{Z}^{st}}\Bigm|_p
, \frac{\partial}{\partial Z^{kl}}\Bigm|_p+\frac{\partial}{\partial
\overline{Z}^{kl}}\Bigm|_p  ,
i\frac{\partial}{\partial Z^{\mu\nu}}\Bigm|_p-i\frac{\partial}{\partial \overline{Z}^{\mu\nu}}\Bigm|_p \Bigr)\\
&=&2\delta_{sk}\delta_{k\mu}\delta_{tl}\delta_{l\nu}+\delta_{ks}\delta_{k\mu}\delta_{t\nu}+\delta_{s\mu}\delta_{tl}\delta_{l\nu}
\end{eqnarray*}
and therefore
\begin{displaymath}
R\Bigl(\frac{\partial}{\partial X^{kl}}\Bigm|_p,
\frac{\partial}{\partial Y^{st}}\Bigm|_p,
 \frac{\partial}{\partial X^{kl}}\Bigm|_p ,
  \frac{\partial}{\partial Y^{\mu\nu}}\Bigm|_p\Bigr) = \left\{ \begin{array}{ll}
1
&\textrm{if $\mu=s\neq k, l=t=\nu $},\\

1
& \textrm{if $\mu=s=k, l\neq t=\nu$},\\

4
& \textrm{if $\mu=s=k, l=t=\nu$},\\

0 & \hbox{otherwise}.
\end{array} \right.
\end{displaymath}
$\Box$\vspace{2mm}

\section{Evolution along the mean curvature flow}
\setcounter{equation}{0}

\subsection{Preliminaries}

For convenience we review results in Section 2 of \cite{MeWa}. A
real $2N$-dimensional Hermitian vector space is a real
$2N$-dimensional vector space  $V$ equipped with a Hermitian
structure, i.e.  a triple $(\omega, J, g)$ consisting of a
symplectic bilinear form $\omega:V\times V\to\R$, an inner product
$g$  and an complex structure $J$ on $V$ satisfying
$g=\omega\circ({\rm Id}\times J)$. A {\it Hermitian isomorphism}
from $(V, \omega,  J,  g)$ to another Hermitian vector space
$(\widetilde V, \tilde\omega, \tilde J, \tilde g)$ of real $2n$
dimension is  a linear isomorphism  $L: V\to\widetilde V$
satisfying: $LJ=\tilde JL$, $L^\ast\tilde\omega=\omega$ and
$L^\ast\tilde g=g$.  Proposition 1 and Corollary 2 in Section 2.1 of
\cite{MeWa} can be summarized as follows.

\begin{proposition}\label{prop:3.1}
For any  linear symplectic isomorphism  $L$ from the real
$2N$-dimensional Hermitian $(V, \omega, J,  g)$ to $(\widetilde V,
\tilde\omega, \tilde J, \tilde g)$, let $L^\star:\widetilde V\to V$
be the adjoint of $L$ determined by $g(L^\star \tilde u, v)=\tilde
g(\tilde u, Lv)$. Then $L^\star L:V\to V$ is positive definite, and
$E:=L(L^\star L)^{-1/2}$ gives rise to a Hermitian isomorphism
 from $(V, \omega,  J, g)$ to $(\widetilde
V, \tilde\omega, \tilde J, \tilde g)$. Moreover, there exists an
{\rm unitary basis} $\{v_{1},\cdots, v_{2N}\}$ of $(V, \omega, J,
g)$, i.e.,
$$
g(v_i, v_j)=\delta_{ij}\quad\hbox{and}\quad Jv_{2k-1}=v_{2k},\;
k=1,\cdots,N,
$$
 (and hence an unitary basis  of $(\widetilde V, \tilde\omega,
\tilde J, \tilde g)$,
$\{\tilde{v}_{1},\cdots,\tilde{v}_{2N}\}=\{E(v_{1}),\cdots,E(v_{2N})\}$
 ), such that
\begin{description}
\item[(i)] The matrix representations of $J$ and $\tilde J$ under them
are all
\begin{equation}\label{e:3.1}
J_0=\left(
     \begin{array}{ccccc}
       0 & 1 &  &  & \\
       -1 & 0 &  &  &  \\
        &  & \ddots &  &  \\
        &  &  & 0 & 1 \\
        &  &  & -1 & 0 \\
     \end{array}
   \right)
\end{equation}

\item[(ii)] The map $(L^{\star} L)^{1/2}$ has the matrix representation
under the basis $\{v_{1},\cdots, v_{2N}\}$,
$$(L^{\star}L)^{1/2}={\rm
Diag}(\lambda_1,\lambda_2,\cdots,\lambda_{2N-1}, \lambda_{2N}),
$$
where $\lambda_{2i-1}\lambda_{2i}=1$ and $\lambda_{2i-1}\le
1\le\lambda_{2i}$, $i=1,\cdots,N$.

\item[(iii)] Under the bases $\{v_{1},\cdots, v_{2N}\}$ and
$\{\tilde{v}_{1},\cdots,\tilde{v}_{2N}\}$ the map $L$ has the matrix
representation $L={\rm
Diag}(\lambda_1,\lambda_2,\cdots,\lambda_{2N-1}, \lambda_{2N})$.
\end{description}
\end{proposition}

\begin{remark}\label{rm:3.2}
{\rm From the arguments in \cite{MeWa} one  can also choose the
$\{v_{1},\cdots, v_{2N}\}$ such that $\lambda_k, k=1,\cdots,2N$ in
Proposition~\ref{prop:3.1}(ii) satisfy: $\lambda_{2i}\le
1\le\lambda_{2i-1}$, $i=1,\cdots,N$.}
\end{remark}

 Let $(M, \omega, J, g)$ and $(\widetilde{M}, \tilde\omega, \tilde
J, \tilde g)$ be two real $2N$-dimensional K\"{a}hler-Einstein
manifolds, and let $\pi_1:M\times\widetilde M\to M$ and
$\pi_2:M\times\widetilde M\to\widetilde M$ be two natural
projections. We have a product K\"{a}hler manifold
$(M\times\widetilde M, \pi_1^\ast\omega-\pi_2^\ast\tilde\omega,
{\cal J}, G)$, where $G=\pi_1^\ast g+ \pi_2^\ast \tilde g$ and
${\cal J}(u,v)=(Ju, -\tilde Jv)$ for $(u, v)\in T(M\times\widetilde
M)$.

For a symplectomorphism  $\varphi:(M, \omega)\rightarrow
(\widetilde{M}, \tilde{\omega})$ let
$$\Sigma={\rm
Graph}(\varphi)=\{(p, \varphi(p))\,|\, p\in M\},
$$
and let $\Sigma_t$ be the mean curvature flow of $\Sigma$ in
$M\times\widetilde{M}$.

Denote by $\Omega:=\pi_1^\ast\omega^N$, and by $\ast\Omega$ the
Hodge star of $\Omega|_{\Sigma_t}$ with respect to the induced
metric on $\Sigma_t$ by $G$. Then $\ast\Omega$ is the Jacobian of
the projection from $\Sigma_t$ onto $M$, and
$\ast\Omega(q)=\Omega(e^1,\cdots,e^{2N})$ for $q\in\Sigma_t$ and any
oriented orthogonal basis $\{e^1,\cdots,e^{2N}\}$ of $T_q\Sigma_t$.
The implicit function theorem implies that $\ast\Omega(q)>0$ if and
only if $\Sigma_t$ is locally a graph over $M$ at $q$.

Let $q=(p, \varphi_t(p))\in\Sigma_{t}\subset
M\times\widetilde{M}$. Set $L:=D_p\varphi_t:T_{p}M\rightarrow
T_{\varphi_t(p)}\widetilde{M}$ and
$E:=D_p\varphi_t[(D_p\varphi_t)^{\star}D_p\varphi_t]^{-\frac{1}{2}}:T_{p}M\rightarrow
T_{\varphi_t(p)}\widetilde{M}$. Since $L^{\ast}L$ is a positive
definite matrix, by the above arguments one can choose a
holomorphic local coordinate system $\{z^{1},\cdots,z^{N}\}$
around $p$, $z^{j}=x^{j}+iy^{j}$, $j=1,\cdots,N$, such that
\begin{description}
\item[(i)] $\{\frac{\partial}{\partial x^{1}}|_p,\cdots, \frac{\partial}{\partial
x^{N}}|_p,\frac{\partial}{\partial
y^{1}}|_p,\cdots,\frac{\partial}{\partial y^{N}}|_p\}$ is an
orthogonal basis of the real $2N$-dimensional vector space $T_{p}M$,

\item[(ii)] The complex structure $J_p$
is given by the matrix $J_0$ in (\ref{e:3.1}) with respect to the
base $\frac{\partial}{\partial x^{1}},
      \frac{\partial}{\partial y^{1}},
      \cdots,
      \frac{\partial}{\partial x^{N}},
      \frac{\partial}{\partial y^{N}}$.

\item[(iii)]  $L^{\ast}L={\rm Diag}(\lambda^2_1,\lambda^2_2,\cdots,
\lambda^2_{2N-1}, \lambda^2_{2N})$ with respect to these basis,
where $\lambda_{2i-1}\lambda_{2i}=1$, $\lambda_{2i-1}\leq
1\le\lambda_{2i}$ for $i=1,\cdots,N$. Obviously
$\frac{\partial}{\partial x^{j}}=\frac{\partial}{\partial
z^{j}}+\frac{\partial}{\partial \overline{z}^{j}}$,
$\frac{\partial}{\partial
y^{j}}=\frac{1}{\sqrt{-1}}(\frac{\partial}{\partial
\overline{z}^{j}}-\frac{\partial}{\partial z^{j}})$.

\item[(iv)] There exists a Hermitian vector space isomorphism
$$
E: (T_pM, \omega_p, J_p, g_p)\to (T_{\varphi_t(p)}\widetilde{M},
\tilde\omega_{\varphi_t(p)}, \tilde J_{\varphi_t(p)}, \tilde
g_{\varphi_t(p)})
$$
such that  under the orthogonal basis of
$(T_{\varphi_t(p)}\widetilde M, \tilde g_{\varphi_t(p)})$,
$$
\left\{E(\frac{\partial}{\partial x_1}|_p),
E(\frac{\partial}{\partial y_1}|_p),\cdots,
E(\frac{\partial}{\partial x_N}|_p), \cdots,
E(\frac{\partial}{\partial y_N}|_p)\right\},
$$
 $\tilde J_{\varphi_t(p)}$ is also given by the matrix $J_0$ in
 (\ref{e:3.1}).
\end{description}

By the choose of basis, we have
\begin{eqnarray*}
&&g(\frac{\partial}{\partial x^{i}}|_p,\frac{\partial}{\partial
x^{j}}|_p)=g(\frac{\partial}{\partial
y^{i}}|_p,\frac{\partial}{\partial
y^{j}}|_p)=\delta_{ij},\\
&&g(\frac{\partial}{\partial x^{i}}|_p, \frac{\partial}{\partial
y^{j}}|_p)=g(\frac{\partial}{\partial
y^{i}}|_p,\frac{\partial}{\partial
x^{j}}|_p)=0,\\
&&g_{l\overline{d}}=g(\frac{\partial}{\partial
z^{l}}|_p,\frac{\partial}{\partial
\overline{z}^{d}}|_p)=g_{\overline{d}l}=g_{d\overline{l}}=g_{\overline{l}d}=\frac{\delta_{ld}}{2},\\
&& g_{\overline{l}\overline{d}}=g_{ld}=0.
\end{eqnarray*}
For $j=1,\cdots,N$, set
\begin{equation}\label{e:3.2}
a^{2j-1}=\frac{\partial}{\partial x_{j}}|_p\quad\hbox{and}\quad
a^{2j}=\frac{\partial}{\partial y_{j}}|_p.
\end{equation}
  Then by ({\bf ii}) above it
holds that
$$
J_p(a^{2j-1})=a^{2j}\quad\hbox{and}\quad J_p(a^{2j})=-a^{2j-1},\quad
j=1,\cdots, N.
$$
Let $s'=s+ (-1)^{s+1}$, $s=1,\cdots,2N$, and let $J_{rs}:=g(Ja_s,
a_r)$. It follows that
$$
J_{s's}=-J_{ss'}\quad\hbox{and}\quad J_{rs}=
\left\{\begin{array}{ll}
 0   &\quad\hbox{if}\;r\ne s',\\
(-1)^{s+1} &\quad\hbox{if}\; r=s'.
\end{array}\right.
$$
For $i=1,\cdots,2N$, let
\begin{equation}\label{e:3.3}
e^{i}=\frac{1}{\sqrt{1+\lambda^{2}_{i}}}(a^{i},\lambda_{i}E(a^{i}))\quad\hbox{and}\quad
e^{2N+i}=\frac{1}{\sqrt{1+\lambda^{2}_{i}}}(J_{p}a^{i},-\lambda_{i}E(J_{p}a^{i})).
\end{equation} They form an orthogonal basis of
$T_{q}(M\times\widetilde{M})$, and
\begin{eqnarray*}
T_{q}\Sigma_t={\rm span}(\{e^1,\cdots, e^{2N})\}\quad\hbox{and}\quad
N_{q}\Sigma_t={\rm span}(\{ e^{2N+1},\cdots, e^{4N}\})
\end{eqnarray*}
and
$\ast\Omega=\Omega(e^1,\cdots,e^{2N})=1/\sqrt{\prod^{2N}_{j=1}(1+\lambda_j^2)}$.

\

\begin{proposition}\label{prop:3.3}{\rm (\cite[Prop.2]{MeWa})}
Let  $(M,g, J, \omega)$ and $(\widetilde{M},\tilde{g}, \tilde J,
\widetilde\omega)$ be two compact K$\ddot{a}$hler-Einstein manifolds
of real dimension $2N$, and let $\Sigma_t$ be the mean curvature
flow of the graph $\Sigma$ of a symplectomorphism
$\varphi:(M,\omega)\rightarrow (\widetilde{M},\widetilde{\omega})$.
Then  $\ast\Omega$ at each point $q\in\Sigma_{t}$ satisfies the
following equation:
\begin{eqnarray}\label{e:3.4}
\frac{d}{dt}\ast\Omega=\Delta\ast\Omega+\ast\Omega\left\{Q(\lambda_{i},h_{jkl})
+\sum_{k}\sum_{i\neq
k}\frac{\lambda_{i}(R_{ikik}-\lambda_{k}^{2}\widetilde{R}_{ikik})}
{(1+\lambda_{k}^{2})(\lambda_{i}+\lambda_{i'})}\right\}.
\end{eqnarray}
where
$$
Q(\lambda_{i},h_{jkl})=\sum_{i,j,k}h_{ijk}^{2}-2\sum_{k}\sum_{i<j}(-1)^{i+j}
\lambda_{i}\lambda_{j}(h_{i'ik}h_{j'jk}-h_{i'jk}h_{j'ik})
$$
with $i'=i+(-1)^{i+1}$, and  $R_{ijkl}=R(a^{i},a^{j},a^{k},a^{l})$
and
$\widetilde{R}_{ijkl}=\widetilde{R}(E(a^{i}),E(a^{j}),E(a^{k}),E(a^{l}))$
are, respectively, the coefficients of the curvature tensors $R$ and
$\widetilde{R}$ with respect to the chosen bases of $T_{p}M$ and
$T_{f(p)}\widetilde{M}$ as in Proposition~\ref{prop:3.1}.
\end{proposition}

For $\vec{\lambda}=(\lambda_1,\cdots,\lambda_{2N})\in\R^{2N}$,
according to \cite[p.322]{MeWa} let
\begin{equation}\label{e:3.5}
\delta_{\vec{\lambda}}:=\inf\Bigl\{Q(\lambda_{i},h_{jkl})\,\Bigm|\,
h_{ijk}\in\R,\, 1\le i,j,k\le 2N,\; \sum_{i,j,k}h^2_{ijk}=1\Bigr\},
\end{equation}
that is, the smallest eigenvalue of $Q$ at $\vec{\lambda}$, and for
$\Lambda\in [1, \infty)$ let
\begin{eqnarray}
&&\delta_\Lambda:=\inf\left\{\delta_{\vec{\lambda}}\,|\,\frac{1}{\Lambda}\le\lambda_i\le\Lambda\;
\hbox{for}\;i=1,\cdots,2N\right\},\label{e:3.6}\\
&&\Lambda_0(N):=\sup\{\Lambda\,|\,\Lambda\ge
1\;\hbox{and}\;\delta_\Lambda>0\}.\label{e:3.7}
\end{eqnarray}
By Remark 2 and Lemma 4 in \cite{MeWa} (or the proof of
\cite[Prop.3]{MeWa}), $\Lambda_0(1)=\infty$,
 and
 $$
 Q((1,\cdots,1),h_{ijk})\geq
\frac{3-\sqrt{5}}{6}|II|^{2}=\frac{3-\sqrt{5}}{6}\sum_{i,j,k}h_{ijk}^2.
$$
Clearly,  $\delta_{\vec{\lambda}}$ is continuous in
$\vec{\lambda}$, and $[1, \infty)\ni\Lambda\to \delta_\Lambda$ is
nonincreasing. They imply  $\Lambda_0(N)>1$. Note that
$\delta_{\Lambda'}>0$  for every $\Lambda'\in [1, \Lambda_0(N))$.
Indeed, by the definition of supremum we have a $\Lambda\in
(\Lambda', \Lambda_0(N))$ with $\delta_\Lambda>0$. So
$\delta_{\Lambda'}\ge\delta_{\Lambda}>0$. In addition,
(\ref{e:3.5}) and (\ref{e:3.6}) imply
\begin{eqnarray*}
&&\inf\Bigl\{Q(\lambda_{i},h_{jkl})\,\Bigm|\, h_{ijk}\in\R,\;
\sum_{i,j,k}h^2_{ijk}=1,\; \frac{1}{\Lambda'}\le\lambda_i\le\Lambda'\Bigr\}\\
&=&\inf\Bigl\{\delta_{\vec{\lambda}}\,\Bigm|\,
\frac{1}{\Lambda'}\le\lambda_i\le\Lambda'\Bigr\}=\delta_{\Lambda'}
\end{eqnarray*}
for every $\Lambda'\in [1, \Lambda_0(N))$. Hence we get:

\begin{proposition}\label{prop:3.4}{\rm (\cite[Prop. 3]{MeWa})}
Let $Q(\lambda_{i},h_{jkl})$ be the the quadratic form defined in
Proposition~\ref{prop:3.3}. Then for the  constant
$\Lambda_0(N)\in (1, +\infty]$ in (\ref{e:3.7}), which only
depends on $2N=\dim M$,  $Q(\lambda_{i},h_{jkl})$ is nonnegative
whenever
$\frac{1}{\Lambda_{0}(N)}\leq\lambda_{i}\leq\Lambda_{0}(N)$ for
$i=1,\cdots,2N$. Moreover, for any $\Lambda'\in
[1,\Lambda_{0}(N))$ it holds that
$$
Q(\lambda_{i},h_{jkl})\geq\delta_{\Lambda'}\sum_{ijk}h_{jkl}^{2}
$$
whenever $\frac{1}{\Lambda'}\leq\lambda_{i}\leq\Lambda'$ for
$i=1,\cdots,2N$.
\end{proposition}

\subsection{The case of Grassmann manifolds}

Let $\varphi: M=G(n, n+m;\C)\rightarrow \widetilde{M}=G(n, n+m;\C)$
be a $\Lambda$-pinched symplectomorphism and $\Sigma={\rm
Graph}(\varphi)$. For $(p, \varphi_t(p))\in\Sigma_t$, let $a^{j},
\;j=1,\cdots,nm$, be the chosen unitary base of $(T_pG(n, n+m;\C),
J_p, g_p)$ as in Proposition~\ref{prop:3.1}. Then
\begin{eqnarray}
&& R_{ijkl}=R(a^{i},a^{j},a^{k},a^{l}),\label{e:3.8}\\
&&\widetilde{R}_{ijkl}=\widetilde{R}(E(a^{i}),E(a^{j}),E(a^{k}),E(a^{l}))\label{e:3.9}
\end{eqnarray}
are, respectively, the coefficients of the curvature tensors $R$ and
$\widetilde{R}$ with respect to the  chosen unitary bases of
$T_{p}G(n, n+m;\C)$ and $T_{\varphi_t(p)}G(n, n+m;\C)$.

From Corollary~\ref{cor:2.4} it follows that
\begin{eqnarray}\label{e:3.10}
 R_{ijkl}=\widetilde{R}_{ijkl}\quad\forall 1\le i, j, k, l\le 2nm.
\end{eqnarray}
Let $({\cal U}, Z)$ be the local chart around $p$ on $G(n, n+m;\C)$
as in Proposition~\ref{prop:2.3}. The final two equalities in
Proposition~\ref{prop:2.6} show
\begin{eqnarray}\label{e:3.11}
&&R\Bigl(\frac{\partial}{\partial X^{kl}}\Bigm|_p,
\frac{\partial}{\partial Y^{st}}\Bigm|_p,\frac{\partial}{\partial
X^{kl}}\Bigm|_p,\frac{\partial}{\partial Y^{st}}\Bigm|_p\Bigr)\nonumber\\
&&-R\Bigl(\frac{\partial}{\partial X^{kl}}\Bigm|_p,
\frac{\partial}{\partial X^{st}}\Bigm|_p,\frac{\partial}{\partial
X^{kl}}\Bigm|_p,\frac{\partial}{\partial
X^{st}}\Bigm|_p\Bigr)=4\delta_{ks}\delta_{lt}.
\end{eqnarray}
 Writing
 $$Z^{11}, Z^{12},\cdots, Z^{1m}, Z^{21},\cdots,
Z^{2m},\cdots, Z^{n1},\cdots, Z^{nm}
$$
into $z^{1}, z^{2},\cdots, z^{nm}$ we have
\begin{eqnarray}\label{e:3.12}
 e_k:=a^{2k-1}=\frac{\partial}{\partial
x^k}|_p\quad\hbox{and}\quad f_k:=a^{2k}=\frac{\partial}{\partial
y^k}|_p
\end{eqnarray}
for $k=1,\cdots,nm$. Then (\ref{e:3.11}) can be written as
\begin{eqnarray}\label{e:3.13}
&&R\bigl(e_{(k-1)m+l}, f_{(s-1)m+t}, e_{(k-1)m+l},  f_{(s-1)m+t}\bigr)\nonumber\\
&&-R\bigl(e_{(k-1)m+l}, e_{(s-1)m+t}, e_{(k-1)m+l},
e_{(s-1)m+t}\bigr)=4\delta_{ks}\delta_{lt}
\end{eqnarray}
for any $1\leq k,s\leq n$ and $1\leq l,t\leq m$.
 Clearly, this is equivalent to
\begin{equation}\label{e:3.14}
R(e_i, f_j, e_i, f_j)-R(e_i, e_j, e_i, e_j)=4\delta_{ij}
\end{equation}
for $1\le i,j\le nm$.

Now for $M=\widetilde M=G(n, n+m;\C)$, by (\ref{e:3.10})  we may
rewrite the second term in the big bracket of (\ref{e:3.4})  as
follows:
\begin{eqnarray}\label{e:3.15}
&&\sum_{k}\sum_{i\neq
k}\frac{\lambda_{i}(R_{ikik}-\lambda_{k}^{2}\widetilde{R}_{ikik})}
{(1+\lambda_{k}^{2})(\lambda_{i}+\lambda_{i'})}=\sum_{k}\sum_{i\neq
k}\frac{\lambda_{i}(1-\lambda_{k}^{2})R_{ikik}}
{(1+\lambda_{k}^{2})(\lambda_{i}+\lambda_{i'})}\\
&=&\sum_{k=2r-1,i=2s-1,r\neq
s}\frac{\lambda_{2s-1}(1-\lambda_{2r-1}^{2})R(e_s, e_r, e_s,
e_r)}{(1+\lambda_{2r-1}^{2})(\lambda_{2s-1}+\lambda_{2s})}\nonumber\\
&&+\sum_{k=2r-1,i=2s}\frac{\lambda_{2s}(1-\lambda_{2r-1}^{2})R(f_s,
e_r, f_s, e_r)}{(1+\lambda_{2r-1}^{2})(\lambda_{2s-1}+\lambda_{2s})}\nonumber\\
&&+\sum_{k=2r,i=2s-1}\frac{\lambda_{2s-1}(1-\lambda_{2r}^{2})R(e_s,
f_r, e_s, f_r)}{(1+\lambda_{2r}^{2})(\lambda_{2s-1}+\lambda_{2s})}\nonumber\\
&&+\sum_{k=2r,i=2s,r\neq
s}\frac{\lambda_{2s}(1-\lambda_{2r}^{2})R(f_s, f_r, f_s,
f_r)}{(1+\lambda_{2r}^{2})(\lambda_{2s-1}+\lambda_{2s})}\nonumber\\
&=&\sum_{r\neq s}\frac{R(e_s, e_r, e_s,
e_r)}{(\lambda_{2s-1}+\lambda_{2s})}
\left[\frac{\lambda_{2s-1}(1-\lambda_{2r-1}^{2})}{(1+\lambda_{2r-1}^{2})}
+\frac{\lambda_{2s}(1-\lambda_{2r}^{2})}{(1+\lambda_{2r}^{2})}\right]\nonumber\\
&&+\sum_{r,s}\frac{R(e_s, f_r, e_s,
f_r)(\lambda_{2r}^{2}-1)(\lambda_{2s}-
\lambda_{2s-1})}{(\lambda_{2s-1}+\lambda_{2s})(1+\lambda_{2r}^{2})}\nonumber\\
&=&\sum_{r,s}\frac{(\lambda_{2r}^{2}-1)(\lambda_{2s}-\lambda_{2s-1})}{(\lambda_{2s-1}+\lambda_{2s})(1+\lambda_{2r}^{2})}
\left[R(e_s, f_r, e_s, f_r)-R(e_s, e_r, e_s,
e_r)\right]\nonumber\\
&=&\sum_{r,s}\frac{(\lambda_{2r}^{2}-1)(\lambda_{2s}^{2}-1)}{(1+\lambda_{2s}^{2})(1+\lambda_{2r}^{2})}
\left[R(e_s, f_r, e_s, f_r)-R(e_s, e_r, e_s,
e_r)\right]\nonumber\\
&=&\sum_{r=s}4\frac{(\lambda_{2r}^{2}-1)(\lambda_{2s}^{2}-1)}{(1+\lambda_{2s}^{2})(1+\lambda_{2r}^{2})}\hspace{30mm}\hbox{
(by (\ref{e:3.14})}\nonumber\\
&=&
4\sum_{s}\frac{(\lambda_{2s}^{2}-1)^{2}}{(1+\lambda_{2s}^{2})^{2}}.\nonumber
\nonumber
\end{eqnarray}
Hence in the present case (\ref{e:3.4}) becomes
\begin{eqnarray}\label{e:3.16}
\frac{d}{dt}\ast\Omega=\Delta\ast\Omega+\ast\Omega\left\{Q(\lambda_{i},h_{jkl})
+4\sum^{nm}_{s=1}\frac{(\lambda_{2s}^{2}-1)^{2}}{(1+\lambda_{2s}^{2})^{2}}\right\}.
\end{eqnarray}
This and Proposition~\ref{prop:3.4} immediately lead to the
following generalization of \cite[\S3,Cor.4]{MeWa}.

\begin{proposition}\label{prop:3.5}
Let $\Lambda_0=\Lambda_{0}(nm)>1$ be the constant defined by
(\ref{e:3.7}). For any $\Lambda\in [1, \Lambda_{0})$ it holds that
$$
(\frac{d}{dt}-\Delta)\ast\Omega\geq\delta_{\Lambda}\ast\Omega|II|^{2}
+
4\ast\Omega\sum^{nm}_{s=1}\frac{(1-\lambda_{2s}^{2})^{2}}{(1+\lambda_{2s}^{2})^{2}}
$$
whenever $\frac{1}{\Lambda}\leq\lambda_{i}\leq\Lambda$ for
$i=1,\cdots,2nm$. Here $|II|$ is the norm of the second
fundamental form of $\Sigma_{t}$.
\end{proposition}

Recall that $\ast\Omega=1/\sqrt{\prod^{2mn}_{j=1}(1+\lambda_j^2)}=
1/\prod_{i \,{\rm odd}}\frac{1}{\lambda_{i}+\lambda_{i'}}$ on
$\Sigma_t$, where $i'=i+ (-1)^{i+1}$ for $i=1,\cdots,2nm$. For
$\Lambda>1$ and $0<\epsilon<1/2^{nm}$ set
$$
\epsilon(mn,
\Lambda)=\frac{1}{2^{mn}}-\frac{1}{(\Lambda+\frac{1}{\Lambda})^{mn}},\quad
 \Lambda(mn,\epsilon)=\frac{2^{-mn}}{2^{-mn}-\epsilon}+
\sqrt{\Bigl( \frac{2^{-mn}}{2^{-mn}-\epsilon} \Bigr)^2-1 }.
$$
Then $\epsilon(mn, \Lambda)>0$ and $\Lambda(mn,\epsilon)>1$. Lemmas
5 and 6 in \cite{MeWa} showed
\begin{eqnarray*}
&&\frac{1}{\Lambda}\leq\lambda_{i}\leq\Lambda\quad\forall
i\quad\Longrightarrow \quad\frac{1}{2^{mn}}-\epsilon(mn, \Lambda)\le\ast\Omega,\\
&&\frac{1}{2^{mn}}-\epsilon\le\ast\Omega\quad \Longrightarrow\quad
\frac{1}{\Lambda(mn,\epsilon)}\leq\lambda_{i}\leq\Lambda(mn,\epsilon)\quad\forall
i.
\end{eqnarray*}
From these and Proposition~\ref{prop:3.5}  we may repeat the proofs
of Proposition 4 and Corollary 5 in \cite{MeWa} to obtain the
following generalization of them.

\begin{proposition}\label{prop:3.6}
For some $T>0$ let $[0, T)\ni t\to\Sigma_{t}$ be  the mean curvature flow of the graph $\Sigma$
of a symplectomorphism $\varphi:G(n, n+m;\C)\to G(n, n+m;\C)$, where
$G(n, n+m;\C)$ is equipped with the unique (up to $\times$ nonzero
factor) invariant K\"{a}hler-Einstein metric. Let $\ast\Omega(t)$ be
the Jacobian of the projection $\pi_{1}:\Sigma_{t}\rightarrow G(n,
n+m;\C)$.  Suppose for some $\Lambda\in (1, \Lambda_0(nm))$ that
$$
\frac{1}{2^{mn}}-\epsilon=\frac{1}{2^{mn}}-\frac{1}{2^{mn}}\left(1-\frac{2\Lambda}{\Lambda^2+1}\right)
=\frac{1}{2^{mn-1}}\frac{\Lambda}{\Lambda^2+1}\le\ast\Omega(0).
$$
Then along the mean curvature flow  $\ast\Omega$ satisfies
\begin{eqnarray*}
\left(\frac{d}{dt}-\Delta\right)\ast\Omega\geq\delta_{\Lambda}\ast\Omega|II|^{2}
+
4\ast\Omega\sum^{nm}_{s=1}\frac{(1-\lambda_{2s}^{2})^{2}}{(1+\lambda_{2s}^{2})^{2}},
\end{eqnarray*}
where $\delta_{\Lambda}$ is given in (\ref{e:3.6}), and so
$\min_{\Sigma_{t}}\ast\Omega$ is nondecreasing as a function in $t$
and $\Sigma_t$ is the graph of a symplectomorphism $\varphi_t:G(n,
n+m;\C)\to G(n, n+m;\C)$. In particular,  if $\varphi$ is
$\Lambda$-pinched for some $\Lambda\in (1,
\Lambda_1(mn)]\setminus\{\infty\}$, then each $\varphi_t$ is
$\Lambda'_{mn}$-pinched along the mean curvature flow, where
$\Lambda'_{mn}$ is defined by (\ref{e:1.2}). {\rm (Note:}
$\Lambda'_{mn}=\Lambda_0(mn)$ if $\Lambda=\Lambda_1(mn)<\infty$.{\rm
)}
\end{proposition}

\begin{remark}\label{rm:3.7}
{\rm Let $(M, \omega, J, g)$ be a compact totally geodesic
K\"ahler-Einstein submanifold of $(G(n, n+m;\C), h)$  (e.g. $(G^{\rm
II}(n, 2n), h_{\rm II})$ and $(G^{\rm III}(n, 2n), h_{\rm III})$ are
such submanifolds of $(G(n, 2n;\C), h_{\rm I})$), $\dim M=2N$. By
Corollary~\ref{cor:2.5} we immediately obtain corresponding results
with Propositions~\ref{prop:3.5} and \ref{prop:3.6}. }
\end{remark}

\subsection{The case of flat complex tori}

 The following proposition is actually contained in the proof of
Corollary 3 of \cite[p.320]{MeWa}. We still give its proof.

\begin{proposition}\label{prop:3.8}
If $M$ and $\widetilde{M}$ are real $2n$-dimensional K\"{a}hler
manifolds with constant holomorphic sectional curvature $c\ge 0$
(hence are Einstein and have the same scalar curvature), then
\begin{eqnarray*}
\frac{d}{dt}\ast\Omega=\Delta\ast\Omega+\ast\Omega\left\{Q(\lambda_{i},h_{jkl})
+ c\sum_{k\;{\rm odd}}\frac{(1-\lambda_{k}^{2})^{2}}
{(1+\lambda_{k}^{2})^{2}}\right\}.
\end{eqnarray*}
\end{proposition}

\noindent{\textbf{Proof}}. With the choice of bases of $T_{p}M$ and
$T_{f(p)}\widetilde{M}$, (we shall suppress  $|_p$ in
$\frac{\partial}{\partial x^{r}}|_p$ and $\frac{\partial}{\partial
y^{r}}|_p$, $r=1,\cdots,n$ for simplicity), it is easily computed
that
\begin{displaymath}
R_{ikik}=R(a^{i},a^{j},a^{k},a^{l}) = \left\{ \begin{array}{ll}
R(\frac{\partial}{\partial x^{r}},\frac{\partial}{\partial x^{s}},
  \frac{\partial}{\partial x^{r}},\frac{\partial}{\partial x^{s}}) & \textrm{if $i=2r-1,k=2s-1$},\\
R(\frac{\partial}{\partial x^{r}},\frac{\partial}{\partial y^{s}},
 \frac{\partial}{\partial x^{r}},\frac{\partial}{\partial y^{s}}) & \textrm{if $i=2r-1,k=2s$},\\
R(\frac{\partial}{\partial y^{r}},\frac{\partial}{\partial x^{s}},
 \frac{\partial}{\partial y^{r}},\frac{\partial}{\partial x^{s}}) & \textrm{if $i=2r,k=2s-1$},\\
R(\frac{\partial}{\partial y^{r}},\frac{\partial}{\partial
y^{s}},\frac{\partial}{\partial y^{r}},\frac{\partial}{\partial
y^{s}}) & \textrm{if $i=2r,k=2s$}.
\end{array} \right.
\end{displaymath}
Plugging  $\frac{\partial}{\partial x^{j}}=\frac{\partial}{\partial
z^{j}}+\frac{\partial}{\partial \overline{z}^{j}}$,
$\frac{\partial}{\partial
y^{j}}=\frac{1}{i}(\frac{\partial}{\partial
\overline{z}^{j}}-\frac{\partial}{\partial z^{j}})$ into the above
equalities we get
\begin{equation}\label{e:3.17}
R_{ikik}=R_{r\overline{s}r\overline{s}}+R_{s\overline{r}s\overline{r}}-R_{r\overline{s}s\overline{r}}
-R_{s\overline{r}r\overline{s}}
\end{equation}
if $(i, k)=(2r-1, 2s-1)$ or $(i,k)=(2r, 2s)$, and
$$
R_{ikik}=-(R_{r\overline{s}r\overline{s}}+R_{s\overline{r}s\overline{r}}+
R_{r\overline{s}s\overline{r}}+R_{s\overline{r}r\overline{s}})
$$
if $(i, k)=(2r-1, 2s)$ or $(i, k)=(2r, 2s-1)$. Note that
$$
g_{l\overline{d}}=g(\frac{\partial}{\partial
z^{l}},\frac{\partial}{\partial
\overline{z}^{d}})=g_{\overline{d}l}=g_{d\overline{l}}=g_{\overline{l}d}=\frac{\delta_{ld}}{2},\quad
g_{\overline{l}\overline{d}}=g_{ld}=0
$$
and that the nonzero components of the Riemannian curvature in the
complex local system $z^{1},\cdots,z^{n}$  are exactly
$R_{i\bar{j}k\bar{l}}$ and $R_{\bar ij\bar k l}$. Moreover,
$$
R_{i\overline{j}k\overline{l}}=-\frac{c}{2}(g_{i\overline{j}}g_{k\overline{l}}+g_{i\overline{l}}g_{\overline{j}k})
$$
on the K\"{a}hler manifolds of constant holomorphic sectional
curvature $c$ (by Proposition 7.6 of \cite[p. 169]{KoNo}). From
(\ref{e:3.17}) we derive
\begin{displaymath}
R_{ikik}= \left\{ \begin{array}{ll}
-\frac{c}{4}(\delta_{rs}-1) & \textrm{if $(i, k)=(2r-1, 2s-1)$ or $(i, k)=(2r, 2s)$}\\
\frac{c}{4}(3\delta_{rs}+1) & \textrm{if $(i, k)=(2r-1, 2s)$ or
$(i,k)=(2r, 2s-1)$}
\end{array} \right.
\end{displaymath}
This shows that
$$
R_{ikik}=\frac{c}{4}(3\delta_{ik'}+1)\quad\forall i\neq k.
$$
Plugging into (\ref{e:3.4}) yields
\begin{eqnarray*}
\frac{d}{dt}\ast\Omega&=&\Delta\ast\Omega+\ast\Omega\left\{Q(\lambda_{i},h_{jkl})
+\frac{c}{4}\sum_{k}\sum_{i\neq
k}\frac{\lambda_{i}(1-\lambda_{k}^{2})(1+3\delta_{ik'})}
{(1+\lambda_{k}^{2})(\lambda_{i}+\lambda_{i'})}\right\}\\
&=&\Delta\ast\Omega+\ast\Omega\left\{Q(\lambda_{i},h_{jkl}) +
c\sum_{k\;{\rm odd}}\frac{(1-\lambda_{k}^{2})^{2}}
{(1+\lambda_{k}^{2})^{2}}\right\}. \nonumber
\end{eqnarray*}
$\Box$\vspace{2mm}

As in the proof of \cite[\S3,Cor.4]{MeWa},  from  this and
Proposition~\ref{prop:3.4} we immediately get the following result.

\begin{proposition}\label{prop:3.9}
Under the assumptions of Proposition~\ref{prop:3.8},  for  any
$\Lambda\in [1, \Lambda_{0}(n))$ it holds that
\begin{eqnarray}\label{e:3.34}
(\frac{d}{dt}-\Delta)\ast\Omega\geq\delta_{\Lambda}\ast\Omega|II|^{2}
+ c\ast\Omega\sum_{k~{\rm odd}}\frac{(1-\lambda_{k}^{2})^{2}}
{(1+\lambda_{k}^{2})^{2}}
\end{eqnarray}
whenever $\frac{1}{\Lambda}\leq\lambda_{i}\leq\Lambda$ for
$i=1,\cdots,2n$. Here $|II|$ is the norm of the second fundamental
form of $\Sigma_{t}$.
\end{proposition}

{\bf From now on} we shall assume $c=0$. In this case we can improve the pinching condition.

\begin{proposition}\label{prop:3.10}
Under the assumptions of Proposition~\ref{prop:3.8}, if $c=0$ and
$\varphi$ is $\Lambda$-pinched with $\Lambda\in [1, \infty)$ then
$\varphi_{t}$ is still $\Lambda$-pinched on $ [0,T)$, i.e.
$$
\left.\begin{array}{ll}
\frac{1}{\Lambda}\leq\lambda_{i}(0)\leq\Lambda\\
\quad\forall i=1,\cdots,2n
\end{array}\right\}
\quad\Longrightarrow \left\{\begin{array}{ll}
\frac{1}{\Lambda}\leq\lambda_{i}(t)\leq\Lambda\quad\\
\forall i=1,\cdots,2n\quad\hbox{and}\quad \forall t\in [0,
T).\end{array}\right.
$$
Here  $[0, T)$ is the maximal existence interval of the mean
curvature flow, and $T>0$ or $T=\infty$.
\end{proposition}

\noindent{\bf Proof}.
 Since $\lambda_i$, $i=1,\cdots$, are
singular values of a linear symplectic map, we have
$\frac{1}{\lambda_i}\in\{\lambda_1,\cdots,\lambda_{2n}\}$ for
$i=1,\cdots, 2n$. (See Lemma 3 of \cite{MeWa}). So the question is
reduced to prove
$$
\left.\begin{array}{ll}
\lambda_{i}(0)\leq\Lambda\\
\quad\forall i=1,\cdots,2n
\end{array}\right\}
\quad\Longrightarrow \left\{\begin{array}{ll}
\lambda_{i}(t)\leq\Lambda\quad\\
\forall i=1,\cdots,2n\quad\hbox{and}\quad \forall t\in [0,
T).\end{array}\right.
$$
We shall use the method in \cite[Section 4]{TsWa} and \cite{Smo3}
to prove this.

 Let
$a^{j}, \;j=1,\cdots,n$ be as in Proposition~\ref{prop:3.3} with
$N=n$. Set
$$
e^{i}=\frac{1}{\sqrt{1+\lambda^{2}_{i}}}(a^{i},\lambda_{i}E(a^{i}))\quad\hbox{and}\quad
e^{2n+i}=\frac{1}{\sqrt{1+\lambda^{2}_{i}}}(J_{p}a^{i},-\lambda_{i}E(J_{p}a^{i}))
$$
for $i=1,\cdots,2n$.  Identifying the tangent space of $M\times
\widetilde M$  with $TM\oplus T\widetilde M$, let $\pi_{1}$ and
$\pi_{2}$ denote the projection onto the first and second factors in
the splitting. Then
\begin{eqnarray*}
&&\pi_{1}(e^{i})=\frac{a^{i}}{\sqrt{1+\lambda_{i}^{2}}},
~~\pi_{2}(e^{i})=\frac{\lambda_{i}E(a^{i})}{\sqrt{1+\lambda_{i}^{2}}}
;\\
&&\pi_{1}(e^{2n+i})=\frac{Ja^{i}}{\sqrt{1+\lambda_{i}^{2}}},
~~\pi_{2}(e^{2n+i})=\frac{-\lambda_{i}E(Ja^{i})}{\sqrt{1+\lambda_{i}^{2}}}
\end{eqnarray*}
for $i=1,\cdots,2n$.
 Let us  define the following parallel symmetric
two-tensor $S$ by
$$
S(X, Y)=\frac{\Lambda^{2}\langle\pi_{1}(X),\pi_{1}(Y)\rangle-
\langle\pi_{2}(X),\pi_{2}(Y)\rangle}{\Lambda^{2+\Xi}}
$$
for any $X,Y \in T(M\times \widetilde M)$, where $\Xi>0$ is a
parameter determined later.
 Then
 \begin{eqnarray*}
 &&S_{ij}:=S(e^{i},
e^{j})=\frac{(\Lambda^{2}-\lambda_{i}\lambda_{j})\delta_{ij}}{\Lambda^{2+\Xi}\cdot\sqrt{(1+\lambda_{i}^{2})(1+
\lambda_{j}^{2})}},\\
&&S_{r(2n+j)}:=S(e^{r},
e^{2n+j})=\frac{(\Lambda^{2}+\lambda_{r}\lambda_{j})\delta_{rj'}(-1)^{j+1}}{\Lambda^{2+\Xi}\cdot\sqrt{(1+\lambda_{r}^{2})(1+
\lambda_{j}^{2})}},\\
&&S_{(2n+i)(2n+j)}:=S(e^{2n+i},
e^{2n+j})=\frac{(\Lambda^{2}-\lambda_{i}\lambda_{j})\delta_{ij}}{\Lambda^{2+\Xi}\cdot\sqrt{(1+\lambda_{i}^{2})(1+
\lambda_{j}^{2})}}
\end{eqnarray*}
for $i,j,r=1,\cdots,2n$, and the matrix $S=(S_{kl})_{1\le k,l\le
4n}$ can be written in the block form
$$
\left(
     \begin{array}{ccccc}
       A & B  \\
       B^T & D
     \end{array}
   \right)
$$
 where $A=D={\rm
Diag}\Bigl(\frac{(\Lambda^{2}-\lambda^2_{1})}{\Lambda^{2+\Xi}\cdot(1+\lambda_{1}^{2})},\cdots,
\frac{(\Lambda^{2}-\lambda^2_{2n})}{\Lambda^{2+\Xi}\cdot(1+\lambda_{2n}^{2})}\Bigr)$.
So
\begin{equation}\label{e:3.19}
\hbox{$A$ is positive definite on $\Sigma_{t}$ if and only if
$\Lambda^{2}-\lambda_{i}^{2}>0, i=1,\cdots,2n$.}
\end{equation}

Obverse that $e^1,\cdots, e^{2n}$ forms an orthogonal basis for
the tangent space of $\Sigma_t$. As in \cite[Prop.3.2]{TsWa}, the
pullback of $S$ to $\Sigma_t$ satisfies the equation
\begin{equation} \label{e:3.20}
\begin{split}
(\frac{d}{dt}-\Delta)S_{ij}&= -h_{\alpha
li}H_{\alpha}S_{lj}-h_{\alpha jl}H_{\alpha}S_{li}+\mathcal
{R}_{kik\alpha}S_{\alpha j}+\mathcal {R}_{kjk\alpha}S_{\alpha
i} \\
 &+ h_{\alpha kl}h_{\alpha ki}S_{lj}+h_{\alpha kl}h_{\alpha
kj}S_{li}-2h_{\alpha ki}h_{\beta kj}S_{\alpha \beta}
 \end{split}
 \end{equation}
 for $i,j=1,\cdots,2n$,
where $\Delta$ is the rough Laplacian on 2-tensors over
$\Sigma_t$,
 $h_{ijk}=G(\nabla_{e^i}^{M\times\widetilde M}e^j, {\cal
J}e^k)$, and ${\cal R}_{kik\alpha}={\cal R}(e^k, e^i,e^k, e^\alpha)$
is the component of the curvature tensor ${\cal R}$ of
$(M\times\widetilde M, G)$ with ${\cal J}$ and $G=\pi_1^\ast g +
\pi_2^\ast  \tilde g$ as in Section 3.1.

 Consider the $2n\times 2n$ matrix  $(S_{ij}):=(S(e^i,
e^j)_{1\le i,j\le 2n}$. By (\ref{e:3.19}) we only need to prove
$$
(S_{ij})>0\;\hbox{at}\;t=0\Longrightarrow (S_{ij})>0\;\hbox{in}\;[0,
T).
$$
This can be directly derived from the following  analogue of
\cite[Lemma 4.1]{TsWa}. $\Box$\vspace{2mm}

\begin{proposition}\label{prop:3.11}
Let $x^{n+i}=y^i$, $i=1,\cdots,n$, and
$g_{ij}=g(\frac{\partial}{\partial x^i}, \frac{\partial}{\partial
x^j})$, $i,j=1,\cdots,2n$. For any given $\epsilon>0$, there
exists a parameter $\Xi>0$ such that the condition
$(T_{ij}):=(S_{ij})-\epsilon (g_{ij})>0$ is preserved along the
mean curvature flow.
\end{proposition}

\noindent{\textbf{Proof}}:
 Let $\alpha=2n+ \mu$ and $\beta=2n+
\nu$, $\mu, \nu=1,\cdots,2n$. As in \cite{TsWa}, (\ref{e:3.20})
yields
\begin{equation} \label{e:3.22}
\begin{split}
(\frac{d}{dt}-\Delta)T_{ij}&= -h_{\alpha
li}H_{\alpha}T_{lj}-h_{\alpha jl}H_{\alpha}T_{li}+\mathcal
{R}_{kik\alpha}S_{\alpha j}+\mathcal {R}_{kjk\alpha}S_{\alpha
i} \\
 &+ h_{\alpha kl}h_{\alpha ki}T_{lj}+h_{\alpha kl}h_{\alpha
kj}T_{li}+2\epsilon h_{\alpha ki}h_{\alpha kj}-2h_{\alpha
ki}h_{\beta kj}S_{\alpha \beta}.
 \end{split}
 \end{equation}
Let $N_{ij}$ denote the right hand side of (\ref{e:3.22}).
 A vector
$V=(V^1,\cdots, V^{2n})$ is called a null eigenvector $V$ of the
matrix $(T_{ij})$ if $\sum_{j}T_{ij}V^{j}=0\;\forall i$. By the
Hamilton's maximum principle \cite[Theorem 9.1]{Ha}, if we may prove
$$
\sum_{ij}N _{ij}V^{i}V^{j}\geq0
$$
for any null eigenvector $V$ of the matrix $(T_{ij})$, then the fact
that $(T_{ij})\ge 0$ at $t=0$ implies that $(T_{ij})\ge 0$ on $[0,
T)$, i.e. Proposition~\ref{prop:3.11} holds.

By a direct computation we only need to prove that at $t=0$
\begin{eqnarray}\label{e:3.23}
\sum_{ij}N _{ij}V^{i}V^{j}&=&\sum_{i,j,k,\alpha}\Bigl[2\epsilon
h_{\alpha ki}h_{\alpha kj}V^{i}V^{j}-2\sum_{\beta}h_{\alpha
ki}h_{\beta
kj}S_{\alpha\beta}V^{i}V^{j}\Bigr]\nonumber\\
&&+ 2\sum_{i,j,k,\alpha}\mathcal {R}_{kik\alpha}S_{\alpha
j}V^{i}V^{j}\geq 0
\end{eqnarray}
for any null eigenvector $V=(V^1,\cdots, V^{2n})$ of the matrix
$(T_{ij})$. It is easily estimated that
\begin{eqnarray*}
&&2\sum_{i,j,k,\alpha,\beta}h_{\alpha ki}h_{\beta
kj}S_{\alpha\beta}V^{i}V^{j}\\
&=&2\sum_{i,j,k,\mu,\nu}h_{2n+\mu,ki}h_{2n+\nu,kj}S_{2n+\mu,2n+\nu }V^{i}V^{j}\\
&=&2\sum_{i,j,k,\mu,\nu}\frac{h_{2n+\mu,ki}h_{2n+\nu,kj}(\Lambda^{2}
-\lambda_{\mu}\lambda_{\nu})\delta_{\mu\nu}V^{i}V^{j}}{\Lambda^{2+
\Xi}\cdot\sqrt{(1+\lambda_{\mu}^{2})(1+\lambda_{\nu}^{2})}}\\
&=&2\sum_{k}\sum_{\mu}\Bigl(\sum_{i,j}h_{2n+\mu,ki}h_{2n+
\mu,kj}V^{i}V^{j}\Bigr)\frac{\Lambda^{2}
-\lambda_{\mu}^{2}}{\Lambda^{2+\Xi}\cdot(1+\lambda_{\mu}^{2})}\\
&\leq&2\sum_{\mu}\sum_{k}\Bigl(\sum_{i,j}h_{2n+\mu,ki}h_{2n+\mu,kj}V^{i}V^{j}\Bigr)
\sum_{\nu}\frac{\Lambda^{2}
-\lambda_{\nu}^{2}}{\Lambda^{2+\Xi}\cdot(1+\lambda_{\nu}^{2})}\\
&\leq&\frac{4n}{\Lambda^{\Xi}}\sum_{i,j,k,\mu}h_{2n+\mu,ki}h_{2n+\mu,kj}V^{i}V^{j}.
\end{eqnarray*}
Here in the first inequality we used the facts\\
$\bullet$ $\sum_{i}(a_{i}b_{i})\leq (\sum_{i}a_{i})(\sum_{i}b_{i})$
for $a_{i}\geq0, b_{i}\geq0$, and\\
$\bullet$
$\sum_{i,j}h_{2n+\mu,ki}h_{2n+\mu,kj}V^{i}V^{j}=(\sum_{i}h_{2n+\mu,ki}V^i)^2
\geq0$,\\
 and the second one comes from the inequality
$$
\sum_{\nu}\frac{\Lambda^{2}
-\lambda_{\nu}^{2}}{\Lambda^{2+\Xi}\cdot(1+\lambda_{\nu}^{2})}\le
\sum_{\nu}\frac{\Lambda^{2} }{\Lambda^{2+\Xi}}\le
\frac{2n}{\Lambda^{\Xi}}.
$$
 So the first sum  in the right side of (\ref{e:3.23}) becomes
\begin{eqnarray*}
&&\sum_{i,j,k,\alpha}\Bigl[2\epsilon h_{\alpha ki}h_{\alpha
kj}V^{i}V^{j}-2\sum_{\beta}h_{\alpha ki}h_{\beta
kj}S_{\alpha\beta}V^{i}V^{j}\Bigr]\\
&&\geq\sum_{i,j,k,\mu}h_{2n+\mu,ki}h_{2n+\mu,kj}V^{i}V^{j}
\Bigl(2\epsilon-\frac{4n}{\Lambda^{\Xi}}\Bigr)
\end{eqnarray*}
because $\alpha=2n+\mu$ and $\beta=2n+\nu$, $\mu, \nu=1,\cdots,2n$.

For a given $\epsilon>0$ we can choose   $\Xi>0$ so large that
$\epsilon-\frac{2n}{\Lambda^{\Xi}}>0$. Then (\ref{e:3.23}) is proved
if we show
$$
\sum_{i,j,k,\alpha}\mathcal {R}_{kik\alpha}S_{\alpha
j}V^{i}V^{j}\geq 0
$$
for any null eigenvector $V$ of the matrix $(T_{ij})$. But this is
obvious because $(M\times\widetilde M, G)$ is flat and hence ${\cal
R}=0$. $\Box$\vspace{2mm}

 From Propositions~\ref{prop:3.9},~\ref{prop:3.10}  we immediately obtain
the following strengthen analogue of Proposition~\ref{prop:3.9}.

\begin{proposition}\label{prop:3.12}
For some $T>0$ let $[0, T)\ni t\to\Sigma_{t}$ be  the mean curvature
flow of the graph $\Sigma$ of a symplectomorphism
$\varphi:M\rightarrow\widetilde{M}$, where $M$ and $\widetilde{M}$
are K\"{a}hler-Einstein manifolds of constant holomorphic sectional
curvature $0$. Let $\ast\Omega(t)$ be the Jacobian of the projection
$\pi_{1}:\Sigma_{t}\rightarrow M$. For the constant $\Lambda_{0}(n)$
in (\ref{e:3.7}) and any $\Lambda\in [1, \Lambda_0(n))$, if
$\varphi$ is $\Lambda$-pinched initially, then $\ast\Omega$
satisfies
\begin{eqnarray*}
\left(\frac{d}{dt}-\Delta\right)\ast\Omega\geq\delta_{\Lambda}\ast\Omega|II|^{2}
\end{eqnarray*}
along the mean curvature flow, where $\delta_{\Lambda}$ is given in
(\ref{e:3.6}). In particular, $\min_{\Sigma_{t}}\ast\Omega$ is
nondecreasing as a function in $t$.
\end{proposition}

\section{Proofs of Theorems~\ref{th:1.1},~\ref{th:1.2} and \ref{th:1.3}}
\setcounter{equation}{0}

\subsection{Proofs of Theorems~\ref{th:1.1},~\ref{th:1.2}}
\setcounter{equation}{0}

Using Propositions~\ref{prop:3.5} and~\ref{prop:3.6} (resp.
Remark~\ref{rm:3.7}) and almost repeating the arguments in \S3.3,
\S3.4 of \cite{MeWa} we can complete the proof of
Theorem~\ref{th:1.1} (resp. Theorem~\ref{th:1.2}).

\subsection{Proof of Theorem~\ref{th:1.3}}\label{sec:5}
\setcounter{equation}{0}

\subsubsection{The long-time existence}
 Embedding
$M\times\widetilde M$  into some $\mathbb{R}^{N}$ isometrically,
as in \cite{MeWa}  the mean curvature flow equation   can be
written as $\frac{d}{dt}F(x,t)=H=\overline{H}+ V$ in terms of the
coordinate function $F(x,t)\in \mathbb{R}^{N}$, where $H\in
T_{\Sigma_t}(M\times\widetilde M)/T\Sigma_t$ and $\overline{H}\in
T_{\Sigma_t}\mathbb{R}^{N}/T\Sigma_t$ are  the mean curvature
vectors of $\Sigma_t$ in $M\times\widetilde M$ and
$\mathbb{R}^{N}$, respectively,  and $V=-II_{M}(e_a,e_a)$. Suppose
by a contradiction that there is a singularity at space time point
$(y_0,t_0)\in \mathbb{R}^{N}\times \mathbb{R}$. Let $d\mu_t$
denote the volume form of $\Sigma_t$, and let
$$
\rho_{(y_0,t_0)}(y,t)=\frac{1}{(4\pi(t_0-t))^n}\exp\left(\frac{-|y-y_0|^2}{4(t_0-t)}\right)
$$
be the backward heat kernel of $\rho_{(y_0,t_0)}$ at $(y_0,t_0)$.
Under our present assumptions, as in \cite[page 328]{MeWa}  we can
still use Proposition~\ref{prop:3.12} to derive the corresponding
inequality of \cite[page 328]{MeWa}, that is,
\begin{eqnarray*}
&&\frac{d}{dt}\int(1-\ast\Omega)\rho_{(y_0,t_0)}d\mu_t\\
&\leq& -\delta_{\Lambda}\int\ast\Omega\|II\|^2\rho_{(y_0,t_0)}d\mu_t
-\int(1-\ast\Omega)\rho_{(y_0,t_0)}\left\|\frac{F^{\perp}}{2(t_0-t)}+
\overline{H}+\frac{V}{2}\right\|^2d\mu_t\\
&&+\int(1-\ast\Omega)\rho_{(y_0,t_0)}\frac{\|V\|^2}{4}d\mu_t.
\end{eqnarray*}
 Then  the expected long-time existence  can be obtained by repeating the
remain arguments on the pages 328$\sim$330 of \cite{MeWa}.

\subsubsection{The convergence}

Let $\varphi: M\rightarrow\widetilde{M}$ be a $\Lambda$-pinched
symplectomorphism with $\Lambda\in (1, \Lambda_0(n))$. Take an
arbitrary $\Lambda_1\in (\Lambda, \Lambda_0(n))$.

\begin{lemma}\label{lem:5.1}
{\rm (\textbf{Djokovic inequality})}:
$$
\tan x\left\{\begin{array}{ll}
 > x+\frac{1}{3}x^{3},&\quad\hbox{if}\; \ 0<x<\frac{\pi}{2},\\
 <x+ f(\alpha)x^3,&\quad\hbox{if}\;0<x<\alpha<\frac{\pi}{2},
 \end{array}\right.
$$
where $f(\alpha)=\frac{\tan \alpha-\alpha}{\alpha^3}$, in particular
$f(\frac{\pi}{6})<\frac{4}{9}$.
\end{lemma}

The following lemma is key for us.

\begin{lemma}\label{lem:5.2}
For every $\Lambda_1\in [1, \Lambda_0(n))$ there exists a
$\widehat\Lambda_1>1$ such that for every $\Lambda\in (1,
\widehat\Lambda_1)$ we have $k,l>0$ to satisfy
\begin{eqnarray}
&&\frac{\pi}{2}\cdot
2^{nl}>\sqrt{\frac{\sqrt{21}-3}{2}}\cdot(\Lambda+\frac{1}{\Lambda})^{nl},\label{e:5.1}\\
&&\frac{l\delta_{\Lambda_1}}{10}\geq\frac{\tan(k(\frac{1}{2^{n}})^{l})}{k(\frac{1}{2^{n}})^{l}},\label{e:5.2}\\
&& \frac{\pi}{2}>k\cdot(\frac{1}{2^{n}})^{l}>
k\cdot(\frac{1}{(\Lambda+\frac{1}{\Lambda})^{n}})^{l}\geq
\sqrt{\frac{\sqrt{21}-3}{2}}.\label{e:5.3}
\end{eqnarray}
Moreover $\widehat\Lambda_1\ge
\left(2\exp\Bigl(\frac{0.141446\delta_{\Lambda_1}}{5n} \Bigr)+
2\exp\Bigl(\frac{0.141446\delta_{\Lambda_1}}{10n}
\Bigr)\sqrt{\exp\Bigl(\frac{0.141446\delta_{\Lambda_1}}{5n} \Bigr)
-1}-1\right)^{1/2}$.
\end{lemma}

Its proof will be given at the end of this section.

By the assumption of Theorem~\ref{th:1.3} we have $\Lambda_1\in
(\Lambda, \Lambda_0)$ such that $\Lambda<\widehat\Lambda_1$. Fix
this $\Lambda_1$ below. By Proposition~\ref{prop:3.12}  we have
\begin{equation}\label{e:5.4}
\frac{d}{dt}\ast\Omega\geq\Delta\ast\Omega+
\delta_{\Lambda_1}\cdot\ast\Omega\cdot |II|^{2}.
\end{equation}
From \cite[Section 7]{Wa2} we also know that
\begin{eqnarray}\label{e:5.5}
\frac{d}{dt}|II|^{2}&=&\Delta|II|^{2}-2|\nabla II|^{2}
+2\left[(\overline{\nabla}_{\partial_k}\overline{R})_{\underline{s}ijk}
+(\overline{\nabla}_{\partial_j}\overline{R})_{\underline{s}kik})\right]h_{sij}\nonumber\\
&&-4\overline{R}_{lijk}h_{slk}h_{sij}+8\overline{R}_{\underline{s}\ \underline{t}jk}h_{tik}h_{sij}\nonumber\\
&&-4\overline{R}_{lkik}h_{slj}h_{sij}+2\overline{R}_{\underline{s}k\underline{t}k}h_{tij}h_{sij}\nonumber\\
&&+2\sum_{s,t,i.m}(\sum_{k}(h_{sik}h_{tmk}-h_{smk}h_{tik}))^{2}\nonumber\\
&&+2\sum_{i,j,m,k}(\sum_{s}h_{sij}h_{smk})^{2},
\end{eqnarray}
where $\overline{R}$ is the curvature tensor and $\overline{\nabla}$
is the covariant derivative of the ambient space,
$\underline{s}=2n+s$.
 Now on one hand
\begin{eqnarray}\label{e:5.6}
&& 2\sum_{s, t, i, m}(\sum_{k}(h_{sik}h_{tmk}-h_{smk}h_{tik}))^{2}+2\sum_{i, j, m, k}(\sum_{s}h_{sij}h_{smk})^{2}\nonumber\\
&\leq & 4\sum_{s, t, i,
m}\Bigl[(\sum_{k}|h_{sik}|^2)(\sum_{k}|h_{tmk}|^2)+
(\sum_{k}|h_{smk}|^2)(\sum_{k}|h_{tik}|^2)\Bigr]\nonumber\\
&&\hspace{20mm}+2\sum_{i, j, m, k}(\sum_{s}h_{sij}^{2})(\sum_{s}h_{smk}^{2})\nonumber\\
&=& 8\sum_{s, i, k}h_{sik}^{2}\sum_{t, m, k}h_{tmk}^{2}+2(\sum_{ s, i, j}h_{sij}^{2})(\sum_{ s, m, k}h_{smk}^{2})\nonumber\\
&=&8|II|^4+2|II|^{4}=10|II|^4,
\end{eqnarray}
where the first inequality comes from
\begin{eqnarray*}
&&\bigl(\sum_{k}(h_{sik}h_{tmk}-h_{smk}h_{tik})\bigr)^{2}\\
&&\le
\bigl(\sum_{k}(|h_{sik}h_{tmk}|+ |h_{smk}h_{tik}|)\bigr)^{2}\\
&&\le
\bigl((\sum_{k}|h_{sik}|^2)^{\frac{1}{2}}(\sum_{k}|h_{tmk}|^2)^{\frac{1}{2}}+
(\sum_{k}|h_{smk}|^2)^{\frac{1}{2}}(\sum_{k}|h_{tik}|^2)^{\frac{1}{2}}\bigr)^{2}\\
&&\le 2\bigl[(\sum_{k}|h_{sik}|^2)(\sum_{k}|h_{tmk}|^2)+
(\sum_{k}|h_{smk}|^2)(\sum_{k}|h_{tik}|^2)\bigr].
\end{eqnarray*}
This and (\ref{e:5.5})-(\ref{e:5.6}) lead to
\begin{equation}\label{e:5.7}
\frac{d}{dt}|II|^{2}\leq\Delta|II|^{2}-2|\nabla II|^{2}+ 10
|II|^{4}.
\end{equation}

We hope to prove that $\max_{\Sigma_t}|II|^{2}\rightarrow 0$ as
$t\rightarrow\infty$. To this goal, for positive numbers $k, l, s$
determined later
 let us compute the evolution equation of
$\frac{|II|^2}{[\sin(k(\ast\Omega)^l)]^s}$ as follows:
\begin{eqnarray*}
&&\frac{d}{dt}\left(\frac{|II|^2}{[\sin(k(\ast\Omega)^l)]^s}\right)\\
&=&\frac{1}{[\sin(k(\ast\Omega)^l)]^s}\frac{d|II|^2}{dt}-
\frac{s\cdot k\cdot l
(\ast\Omega)^{l-1}|II|^2\cos(k(\ast\Omega)^l)}{[\sin(k(\ast\Omega)^l)]^{s+1}}\frac{d\ast\Omega}{dt},\\
&&\Delta\left(\frac{|II|^2}{[\sin(k(\ast\Omega)^l)]^s}\right)\\
&=&\frac{\Delta|II|^2}{[\sin(k(\ast\Omega)^l)]^s}-
\frac{s\cdot k\cdot l\cdot|II|^2\cdot(\ast\Omega)^{l-1}\cdot\cos(k(\ast\Omega)^l)\cdot\Delta\ast\Omega}{[\sin(k(\ast\Omega)^l)]^{s+1}}\\
&&-\frac{2s\cdot k\cdot l\cdot\nabla|II|^2\cdot(\ast\Omega)^{l-1}\cdot\cos(k(\ast\Omega)^l)\cdot\nabla\ast\Omega}{[\sin(k(\ast\Omega)^l)]^{s+1}}\\
&&+\frac{s\cdot k^2\cdot l^2\cdot |II|^2\cdot (\ast\Omega)^{2l-2} \cdot \sin(k(\ast\Omega)^l)\cdot |\nabla\ast\Omega|^{2}}{[\sin(k(\ast\Omega)^l)]^{s+1}}\\
&&-\frac{s\cdot k \cdot l \cdot (l-1) \cdot|II|^2\cdot (\ast\Omega)^{l-2} \cdot \cos(k(\ast\Omega)^l)\cdot |\nabla\ast\Omega|^{2}}{[\sin(k(\ast\Omega)^l)]^{s+1}}\\
&&+\frac{s\cdot (s+1)\cdot k^2\cdot l^2\cdot |II|^2\cdot
(\ast\Omega)^{2l-2} \cdot (\cos(k(\ast\Omega)^l))^2\cdot
|\nabla\ast\Omega|^{2}}{[\sin(k(\ast\Omega)^l)]^{s+2}}
\end{eqnarray*}
and hence
\begin{eqnarray}\label{e:5.8}
&&(\frac{d}{dt}-\Delta)\left(\frac{|II|^2}{[\sin(k(\ast\Omega)^l)]^s}\right)\nonumber\\
&=&\frac{1}{[\sin(k(\ast\Omega)^l)]^s}(\frac{d}{dt}-\Delta)|II|^2\nonumber\\
&&-\frac{s\cdot k\cdot l\cdot|II|^2\cdot(\ast\Omega)^{l-1}
\cdot\cos(k(\ast\Omega)^l)}{[\sin(k(\ast\Omega)^l)]^{s+1}}(\frac{d}{dt}-\Delta)\ast\Omega\nonumber\\
&&+\frac{4\cdot s\cdot k\cdot l\cdot
|II|\cdot\nabla|II|\cdot(\ast\Omega)^{l-1}
\cdot\cos(k(\ast\Omega)^l)\cdot\nabla\ast\Omega}{[\sin(k(\ast\Omega)^l)]^{s+1}}\nonumber\\
&&-\frac{s\cdot k^2\cdot l^2\cdot |II|^2\cdot (\ast\Omega)^{2l-2}
\cdot
\sin(k(\ast\Omega)^l)\cdot |\nabla\ast\Omega|^{2}}{[\sin(k(\ast\Omega)^l)]^{s+1}}\nonumber\\
&&+\frac{s\cdot k \cdot l \cdot (l-1) \cdot|II|^2\cdot
(\ast\Omega)^{l-2} \cdot
\cos(k(\ast\Omega)^l)\cdot |\nabla\ast\Omega|^{2}}{[\sin(k(\ast\Omega)^l)]^{s+1}}\nonumber\\
&&-\frac{s\cdot (s+1)\cdot k^2\cdot l^2\cdot |II|^2\cdot
(\ast\Omega)^{2l-2} \cdot
(\cos(k(\ast\Omega)^l))^2\cdot |\nabla\ast\Omega|^{2}}{[\sin(k(\ast\Omega)^l)]^{s+2}}\nonumber\\
&\leq&\frac{-2|\nabla II|^{2}+ 10 |II|^{4}}{[\sin(k(\ast\Omega)^l)]^s}\nonumber\\
&&-\frac{s\cdot k\cdot
l\cdot\delta_{\Lambda_1}\cdot|II|^4\cdot(\ast\Omega)^{l}\cdot\cos(k(\ast\Omega)^l)}{[\sin(k(\ast\Omega)^l)]^{s+1}}
\hspace{40mm}\hbox{(by (\ref{e:5.4}))}\nonumber\\
&&+\frac{4\cdot s\cdot k\cdot l\cdot
|II|\cdot\nabla|II|\cdot(\ast\Omega)^{l-1}\cdot
\cos(k(\ast\Omega)^l)\cdot\nabla\ast\Omega}{[\sin(k(\ast\Omega)^l)]^{s+1}}\nonumber\\
&&+(\hbox{the\ last\ three\ terms})\nonumber\\
&=&\frac{-2|\nabla II|^{2}}{[\sin(k(\ast\Omega)^l)]^s}+
10[\sin(k(\ast\Omega)^l)]^s\left(\frac{|II|^{2}}{[\sin(k(\ast\Omega)^l)]^s}\right)^2
\nonumber\\
&&-s\cdot k\cdot l\cdot \delta_{\Lambda_1} \cdot (\ast\Omega)^l
\cdot \cos(k(\ast\Omega)^l)\cdot
 [\sin(k(\ast\Omega)^l)]^{s-1}\left(\frac{|II|^{2}}{[\sin(k(\ast\Omega)^l)]^s}\right)^2\nonumber\\
&&+\frac{4\cdot s\cdot k\cdot l\cdot
|II|\cdot\nabla|II|\cdot(\ast\Omega)^{l-1}\cdot
\cos(k(\ast\Omega)^l)\cdot\nabla\ast\Omega}{[\sin(k(\ast\Omega)^l)]^{s+1}}\\
&&+(\hbox{the\ last\ three\ terms})\nonumber.
\end{eqnarray}

Note that the Cauchy-Schwarz inequality implies
\begin{eqnarray*}
|\nabla|II||^{2}&=&\sum_{i=1}^{2n}\left(\nabla_{i}\sqrt{\sum_{j,k,l=1}^{2n}h_{jkl}^{2}}\right)^{2}
=\sum_{i=1}^{2n}\left(\frac{2\sum_{j,k,l}h_{jkl}\partial_{i}h_{jkl}}{2|II|}\right)^{2}\nonumber\\
&\leq &\sum_{i=1}^{2n}\left(\sum_{j,k,l}\frac{h_{jkl}^{2}}{|II|^{2}}\sum_{j,k,l}(\partial_{i}h_{jkl})^{2}\right)\nonumber\\
&\leq & \sum_{i,j,k,l}(\partial_{i}h_{jkl})^2=|\nabla II|^{2}.
\end{eqnarray*}
The term in (\ref{e:5.8}) becomes
\begin{eqnarray*}
&&\frac{4\cdot s\cdot k\cdot l\cdot
|II|\cdot\nabla|II|\cdot(\ast\Omega)^{l-1}\cdot\cos(k(\ast\Omega)^l)
\cdot\nabla\ast\Omega}{[\sin(k(\ast\Omega)^l)]^{s+1}}\\
&\leq&\frac{4\cdot s\cdot k\cdot l\cdot
|II|\cdot|\nabla|II||\cdot(\ast\Omega)^{l-1}\cdot\cos(k(\ast\Omega)^l)
\cdot|\nabla\ast\Omega|}{[\sin(k(\ast\Omega)^l)]^{s+1}}\\
&\leq&\frac{4\cdot s\cdot k\cdot l\cdot |II|\cdot|\nabla
II|\cdot(\ast\Omega)^{l-1}\cdot\cos(k(\ast\Omega)^l)
\cdot|\nabla\ast\Omega|}{[\sin(k(\ast\Omega)^l)]^{s+1}}\\
&=&\frac{4}{[\sin(k(\ast\Omega)^l)]^{s}}\left[\frac{ s\cdot k\cdot
l\cdot |II|\cdot(\ast\Omega)^{l-1}
\cdot\cos(k(\ast\Omega)^l)\cdot|\nabla\ast\Omega|}{\sin(k(\ast\Omega)^l)}\right]\cdot|\nabla II|\\
&\leq&\frac{2}{[\sin(k(\ast\Omega)^l)]^{s}}\left[ \frac{s^2\cdot
 k^2\cdot l^2\cdot |II|^2\cdot(\ast\Omega)^{2l-2}\cdot(\cos(k(\ast\omega)^l))^2\cdot|\nabla\ast\Omega|^2}{(\sin(k(\ast\omega)^l))^2}
  +|\nabla II|^2\right]\\
&=&\frac{2|\nabla II|^{2}}{[\sin(k(\ast\Omega)^l)]^s}+\frac{2\cdot
s^2 \cdot k^2\cdot l^2\cdot
|II|^2\cdot(\ast\Omega)^{2l-2}\cdot(\cos(k(\ast\Omega)^l))^2\cdot|\nabla\ast\Omega|^2}{(\sin(k(\ast\Omega)^l))^{s+2}}.
\end{eqnarray*}
Hence we arrive at
\begin{eqnarray*}
&&(\frac{d}{dt}-\Delta)\left(\frac{|II|^2}{[\sin(k(\ast\Omega)^l)]^s}\right)\\
&\leq&
10[\sin(k(\ast\Omega)^l)]^s\left(\frac{|II|^{2}}{[\sin(k(\ast\Omega)^l)]^s}\right)^2
\\
&&-s\cdot k\cdot l\cdot \delta_{\Lambda_1} \cdot (\ast\Omega)^l
\cdot \cos(k(\ast\Omega)^l)\cdot
[\sin(k(\ast\Omega)^l)]^{s-1}\left(\frac{|II|^{2}}{[\sin(k(\ast\Omega)^l)]^s}\right)^2\\
&&+\frac{2\cdot s^2\cdot k^2\cdot l^2\cdot
 |II|^2\cdot(\ast\Omega)^{2l-2}\cdot(\cos(k(\ast\Omega)^l))^2\cdot|\nabla\ast\Omega|^2}{(\sin(k(\ast\Omega)^l))^{s+2}}\\
&&-\frac{s\cdot k^2\cdot l^2\cdot |II|^2\cdot
(\ast\Omega)^{2l-2} \cdot \sin(k(\ast\Omega)^l)\cdot |\nabla\ast\Omega|^{2}}{[\sin(k(\ast\Omega)^l)]^{s+1}}\\
&&+\frac{s\cdot k \cdot l \cdot (l-1) \cdot|II|^2\cdot
(\ast\Omega)^{l-2}
\cdot \cos(k(\ast\Omega)^l)\cdot |\nabla\ast\Omega|^{2}}{[\sin(k(\ast\Omega)^l)]^{s+1}}\\
&&-\frac{s\cdot (s+1)\cdot k^2\cdot l^2\cdot |II|^2\cdot
(\ast\Omega)^{2l-2} \cdot
 (\cos(k(\ast\Omega)^l))^2\cdot |\nabla\ast\Omega|^{2}}{[\sin(k(\ast\Omega)^l)]^{s+2}}\\
&=&\left(\frac{|II|^{2}}{[\sin(k(\ast\Omega)^l)]^s}\right)^2\cdot[\sin(k(\ast\Omega)^l)]^{s-1}\cdot
\bigl[10\cdot\sin(k(\ast\Omega)^l)-\\
&&\hspace{60mm}s\cdot k\cdot l\cdot \delta_{\Lambda_1} \cdot(\ast\Omega)^l \cdot \cos(k(\ast\Omega)^l)\bigr]\\
&&+\frac{  s\cdot k\cdot
l\cdot(\ast\Omega)^{l-2}|II|^{2}|\nabla\ast\Omega|^{2}}{[\sin(k(\ast\Omega)^l)]^{s+2}}[(s-1)\cdot
k\cdot l\cdot(\ast\Omega)^{l}\cdot(\cos(k(\ast\Omega)^l))^2\\
&&-k\cdot
l\cdot(\ast\Omega)^{l}\cdot(\sin(k(\ast\Omega)^l))^2+(l-1)\cos(k(\ast\Omega)^l)\sin(k(\ast\Omega)^l)].
\end{eqnarray*}
Take $s=1$ we obtain
\begin{eqnarray}\label{e:5.9}
&&(\frac{d}{dt}-\Delta)\left(\frac{|II|^2}{\sin(k(\ast\Omega)^l)}\right)\\
&\leq&\left(\frac{|II|^{2}}{\sin(k(\ast\Omega)^l)}\right)^2\cdot\bigl[10\cdot\sin(k(\ast\Omega)^l)-
k\cdot l\cdot \delta_{\Lambda_1}\cdot(\ast\Omega)^l \cdot \cos(k(\ast\Omega)^l)\bigr]\nonumber\\
&&+\frac{k\cdot
l\cdot(\ast\Omega)^{l-2}|II|^{2}|\nabla\ast\Omega|^{2}}{[\sin(k(\ast\Omega)^l)]^{3}}\bigl[-k\cdot
l\cdot(\ast\Omega)^{l}\cdot
(\sin(k(\ast\Omega)^l))^2\nonumber\\
&&\hspace{60mm}+(l-1)\cos(k(\ast\Omega)^l)\sin(k(\ast\Omega)^l)\bigr].\nonumber
\end{eqnarray}

\begin{claim}\label{cl:5.3}
If the positive numbers $k, l$ satisfy (\ref{e:5.1})-(\ref{e:5.3})
in Lemma~\ref{lem:5.2}, then
$$
10\sin(k(\ast\Omega)^l)-k\cdot
l\cdot\delta_{\Lambda_1}\cdot(\ast\Omega)^{l}\cdot\cos(k(\ast\Omega)^l)<
0
$$
and
$$
(l-1)\cdot\cos(k(\ast\Omega)^l) -k\cdot
l\cdot(\ast\Omega)^{l}\sin(k(\ast\Omega)^l)<0,
$$
that is
\begin{equation}\label{e:5.10}
\frac{l-1}{l\cdot k(\ast\Omega)^l}<\tan(k(\ast\Omega)^l)<\frac{
l\cdot\delta_{\Lambda_1}\cdot k\cdot(\ast\Omega)^{l}}{10}
\end{equation}
for any $\ast\Omega\in [\frac{1}{2^n}, \frac{1}{\Lambda+
\frac{1}{\Lambda}}]$ with $1<\Lambda<\widehat\Lambda_1$.
\end{claim}

We put off its proof. Then (\ref{e:5.9}) becomes
\begin{eqnarray*}
&&(\frac{d}{dt}-\Delta)\left(\frac{|II|^2}{\sin(k(\ast\Omega)^l)}\right)\\
&\leq& \left(\frac{|II|^{2}}{\sin(k(\ast\Omega)^l)}\right)^2
[10\cdot\sin(k(\ast\Omega)^l)-
 k\cdot l\cdot \delta_{\Lambda_1} \cdot(\ast\Omega)^l \cdot \cos(k(\ast\Omega)^l)].
\end{eqnarray*}
Let $g=\frac{|II|^2}{\sin(k(\ast\Omega)^l)}$ and
$$
K_{1}:=\max_{\ast\Omega\in[\frac{1}{(\Lambda+\frac{1}{\Lambda})^{n}},
\frac{1}{2^{n}}]}\left[10\cdot\sin(k(\ast\Omega)^l)- k\cdot l\cdot
\delta_{\Lambda_1} \cdot(\ast\Omega)^l \cdot
\cos(k(\ast\Omega)^l)\right].
$$
By Claim~\ref{cl:5.3}, $K_{1}< 0$ and
\begin{eqnarray}\label{e:5.11}
(\frac{d}{dt}-\Delta)g\leq K_{1}\cdot g^{2}.
\end{eqnarray}

Consider the initial value problem
\begin{equation}\label{e:5.12}
\frac{d}{dt}y= K_{1}\cdot y^{2}\quad\hbox{and}\quad y(0)=\rm
max_{\Sigma_{0}}g.
\end{equation}
The unique solution of it is given by
$y(t)=\frac{y(0)}{1-y(0)K_{1}t}$.
 By (\ref{e:5.11})-(\ref{e:5.12}) the comparison principle
for parabolic equations yields
\begin{equation}\label{e:5.13}
g=\frac{|II|^2}{\sin(k(\ast\Omega)^l)}\le  y(t)\quad\forall t>0.
\end{equation}
Since (\ref{e:5.3}) implies that the function
$$
\Bigl[\frac{1}{(\Lambda+\frac{1}{\Lambda})^{n}},\frac{1}{2^{n}}\Bigr]\ni
\ast\Omega\to \sin(k(\ast\Omega)^l)
$$
is bounded away from zero, we derive
\begin{equation}\label{e:5.14}
\max_{\Sigma_{t}}|II|^{2}\leq
\sin(k(\frac{1}{2^{n}})^l)\cdot\frac{y(0)}{1-y(0)K_{1}t} \rightarrow
0,\ \ t\rightarrow\infty.
\end{equation}
The desired claim is proved.  So up to proofs of Lemma~\ref{lem:5.2}
and Claim~\ref{cl:5.3}, we have proved that  the flow converges to a
totally geodesic Lagrangian submanifold at infinity.\vspace{2mm}

\noindent{\bf Proof of Claim~\ref{cl:5.3}}.  Fix the positive
numbers $k, l$ satisfying (\ref{e:5.1})-(\ref{e:5.3}) in
Lemma~\ref{lem:5.2}. By (\ref{e:5.3}) we have
$$
\frac{\pi}{2}>k\cdot(\frac{1}{2^{n}})^{l}\geq
k\cdot(\ast\Omega)^{l}\geq
k\cdot(\frac{1}{(\Lambda+\frac{1}{\Lambda})^{n}})^{l}\geq
\sqrt{\frac{\sqrt{21}-3}{2}}
$$
because $\ast\Omega\in [\frac{1}{(\Lambda+\frac{1}{\Lambda})^{n}},
\frac{1}{2^{n}}]$. Note that
$$
\sqrt{\frac{\sqrt{21}-3}{2}}=\inf\left\{x(x+\frac{1}{3}x^{3})\geq
1\,\biggm|\, 0<x<\pi/2\right\}\approx 0.8895436175241
$$
sits between $\frac{\pi}{3.5317}$ and $\frac{\pi}{3.5316}$. By
Lemma~\ref{lem:5.1} (the Djokovic inequality) we get
$$
k\cdot(\ast\Omega)^{l}(\tan(k(\ast\Omega)^l))>k\cdot(\ast\Omega)^{l}(k\cdot(\ast\Omega)^{l}
+\frac{1}{3}(k\cdot(\ast\Omega)^{l})^{3})\geq 1>\frac{l-1}{l},
$$
that is, the first inequality in (\ref{e:5.10}). Similarly, the
second inequality in (\ref{e:5.10}) follows from (\ref{e:5.2}).
Claim~\ref{cl:5.3} is proved. $\Box$\vspace{2mm}

\noindent{\bf Proof of Lemma~\ref{lem:5.2}}.
 For conveniences we set
$\tau:=\tau(\Lambda)=\Lambda+\frac{1}{\Lambda}$, which is larger
than $2$ because $\Lambda>1$.
 Since
$\frac{\pi}{2}>\sqrt{\frac{\sqrt{21}-3}{2}}$ we may fix a small
$\epsilon>0$ such that
$$
\frac{\pi}{2}>\frac{\pi}{2}-\epsilon>\sqrt{\frac{\sqrt{21}-3}{2}}.
$$
Set $\alpha=\frac{\pi}{2}-\epsilon$. Then (\ref{e:5.1}) holds for
any
\begin{equation}\label{e:5.15}
l\le\frac{\ln\big(\alpha/\sqrt{\frac{\sqrt{21}-3}{2}}\bigr)}{n\ln\frac{\tau}{2}}.
\end{equation}
More precisely, such a $l$ satisfies
$$
\alpha\cdot 2^{nl}\ge\sqrt{\frac{\sqrt{21}-3}{2}}\cdot\tau^{nl}.
$$
Hence we can always take $k=k_l>0$ such that
$$
\sqrt{\frac{\sqrt{21}-3}{2}}\cdot\tau^{nl}\le k\le \alpha\cdot
2^{nl}
$$
or equivalently
\begin{equation}\label{e:5.16}
\frac{\pi}{2}>\alpha\ge k\cdot(\frac{1}{2^{n}})^{l}>
k\cdot(\frac{1}{(\Lambda+\frac{1}{\Lambda})^{n}})^{l}\geq
\sqrt{\frac{\sqrt{21}-3}{2}}.
\end{equation}

 By the Djokovic inequality
$$
\frac{\tan(k(\frac{1}{2^{n}})^{l})}{k(\frac{1}{2^{n}})^{l}}\le 1+
f(\alpha)\bigl({k(\frac{1}{2^{n}})^{l}}\bigr)^2
$$
if $k\cdot(\frac{1}{2^{n}})^{l}\le \alpha$. So (\ref{e:5.2}) holds
if $k>0$ and $l>0$ are chosen to satisfy
\begin{equation}\label{e:5.17}
\frac{l\delta_{\Lambda_1}}{10}\geq 1+ f(\alpha)\alpha^2\ge 1+
f(\alpha)\bigl({k(\frac{1}{2^{n}})^{l}}\bigr)^2
\end{equation}
or equivalently
\begin{equation}\label{e:5.18}
l\ge \frac{10}{\delta_{\Lambda_1}}\cdot (1+ f(\alpha)\alpha^2).
\end{equation}
Hence we can take $l>0$ to satisfy (\ref{e:5.15}) and (\ref{e:5.18})
if
\begin{equation}\label{e:5.19}
\frac{\ln\big(\alpha/\sqrt{\frac{\sqrt{21}-3}{2}}\bigr)}{n\ln\frac{\tau}{2}}\ge
\frac{10}{\delta_{\Lambda_1}}\cdot (1+ f(\alpha)\alpha^2).
\end{equation}
Since the function
$$
(1, \infty)\to\R,\;\Lambda\mapsto \Lambda+ \frac{1}{\Lambda}
$$
is strictly increasing,  $\log\frac{\tau}{2}\to 0^{+}$  as
$\Lambda\to 1^{+}$. Hence for a given
$$
\frac{\pi}{2}>\alpha>\sqrt{\frac{\sqrt{21}-3}{2}},
$$
 there exists the largest $\Lambda_1^{(\alpha)}>1$ such that
(\ref{e:5.19}) holds for $\tau=\tau_\alpha=\Lambda_1^{(\alpha)}+
1/\Lambda_1^{(\alpha)}$, i.e.
\begin{equation}\label{e:5.20}
g(\alpha):=\frac{\alpha\ln\big(\alpha/\sqrt{\frac{\sqrt{21}-3}{2}}\bigr)}{\tan\alpha}\ge
\frac{10n}{\delta_{\Lambda_1}}\cdot \ln\frac{\tau_\alpha}{2}.
\end{equation}
 Of course, (\ref{e:5.20}) also holds
for for every $\tau= \Lambda+ \frac{1}{\Lambda}$ with $\Lambda\in
(1, \Lambda_1^{(\alpha)})$. Then
$$
\widehat\Lambda_1=\sup\biggl\{\Lambda_1^{(\alpha)}\,\biggm|\,
\sqrt{\frac{\sqrt{21}-3}{2}}<\alpha<\frac{\pi}{2} \quad\hbox{and
(\ref{e:5.19}) holds for $\tau=\tau_\alpha$}\biggr\}
$$
satisfies the desired condition. In Appendix~\ref{app:A} we shall
prove

\begin{claim}\label{cl:5.4}
There exists a unique $\alpha_0\in (\sqrt{\frac{\sqrt{21}-3}{2}},
\frac{\pi}{2})$ such that
$$
g(\alpha_0)=\sup\biggl\{g(\alpha)\,\biggm|\,\sqrt{\frac{\sqrt{21}-3}{2}}<\alpha<\frac{\pi}{2}
\biggr\}.
$$
Moreover $\alpha_0\approx 1.238756$ and $g(\alpha_0)\approx
0.141446$.
\end{claim}

Hence $\widehat\Lambda_1\ge\Lambda_1^{(\alpha_0)}$, where
$\Lambda_1^{(\alpha_0)}$ is determined by
$$
g(\alpha_0)=\frac{10n}{\delta_{\Lambda_1}}\cdot
\biggl[\ln\Bigl(\Lambda_1^{(\alpha_0)}+
\frac{1}{\Lambda_1^{(\alpha_0)}}\Bigr)-\ln 2\biggr],
$$
or more precisely
\begin{eqnarray*}
\Lambda_1^{(\alpha_0)}&=&\left(2\exp\Bigl(\frac{g(\alpha_0)\delta_{\Lambda_1}}{5n}
\Bigr)+ 2\exp\Bigl(\frac{g(\alpha_0)\delta_{\Lambda_1}}{10n}
\Bigr)\sqrt{\exp\Bigl(\frac{g(\alpha_0)\delta_{\Lambda_1}}{5n}
\Bigr) -1}-1\right)^{1/2}\\
&\approx &\left(2\exp\Bigl(\frac{0.141446\delta_{\Lambda_1}}{5n}
\Bigr)+ 2\exp\Bigl(\frac{0.141446\delta_{\Lambda_1}}{10n}
\Bigr)\sqrt{\exp\Bigl(\frac{0.141446\delta_{\Lambda_1}}{5n} \Bigr)
-1}-1\right)^{1/2}.
\end{eqnarray*}
 This
completes the proof of Lemma~\ref{lem:5.2}. $\Box$\vspace{2mm}

In summary the proof of Theorem~\ref{th:1.3} is complete.

~~~~~~~~~~~~~~~~~~~~~~~~~~~~~~~~~~~~~~~~~~~~~~~~~~~~~~~~~~~~~~~~
~~~~~~~~~~~~~~~~~~~~~~~~~~~~~~~~~~~~~$\square$\vspace{2mm}

\section{A concluding remark}\label{sec:6}
\setcounter{equation}{0}

Carefully checking the proofs of Theorems~\ref{th:1.1},~\ref{th:1.2}
we find that our real $2n$-dimensional compact K\"ahler-Einstein
manifolds  $(M, \omega, J, g)$ all satisfy the following three
conditions ({\bf A}), ({\bf B}) and ({\bf C}):

 \noindent{\bf (A)} The curvature tensor $R$
is constant on subbundle
$$
\{(X, JX, Y, JY)\,|\, g(X, Y)=0,\, g(X, JY)=0,\,g(X, X)=1=g(Y,Y)\}.
$$
 In other words, for any $p, q\in M$ and any unit orthogonal bases  of $(T_pM, J_p,
g_p)$ and $(T_qM, J_q, g_q)$, $\{a_1, \cdots, a_{2n}\}$ and $\{a'_1,
\cdots, a'_{2n}\}$ with $a_{2k}=J_pa_{2k-1}$ and
$a'_{2k}=J_qa'_{2k-1}$, $k=1,\cdots, n$, it holds that
$$R(a_i, a_k,
a_i, a_k)= R(a'_i, a'_k, a'_i, a'_k)\quad\forall 1\le i, k\le 2n.
$$

If $(M, \omega, J, g)$ is also homogeneous, this is equivalent to
the following weaker

 \noindent{\bf (A')} For any $p\in M$ and any unit orthogonal bases  of $(T_pM, J_p,
g_p)$, $\{a_1, \cdots, a_{2n}\}$ and $\{a'_1, \cdots, a'_{2n}\}$
with $a_{2k}=J_pa_{2k-1}$ and $a'_{2k}=J_qa'_{2k-1}$, $k=1,\cdots,
n$, it holds that $R(a_i, a_k, a_i, a_k)= R(a'_i, a'_k, a'_i, a'_k)$
for all $1\le i, k\le 2n$.

 \noindent{\bf (B)} ${\rm Re}(R(X,\overline{Y},X,\overline{Y}))\leq
 0$ for any $X,Y\in T^{(1,0)}M$.

 \noindent{\bf (C)} The  holomorphic sectional
curvature is positive, i.e. $\exists\, c_0>0$ such that
$$
R(u, Ju, u, Ju)=-4R\Bigr(\frac{u-\sqrt{-1}Ju}{2},
\frac{u+\sqrt{-1}Ju}{2}, \frac{u-\sqrt{-1}Ju}{2},
\frac{u+\sqrt{-1}Ju}{2}\Bigl)\ge c_0
$$
for any unit vector $u\in TM$.

By Propositions~\ref{prop:2.1},~\ref{prop:2.2} and
Corollaries~\ref{cor:2.4},~\ref{cor:2.5}, the manifolds $(G(n,
n+m;\C), h)$, and $(G^{\rm II}(n, 2n), h_{\rm II})$ and $(G^{\rm
III}(n, 2n), h_{\rm III})$ satisfy these conditions. On the other
hand,
 from (\ref{e:2.14}) we see that $(G^{\rm IV}(1, n+1), h_{\rm IV})$
 does not satisfy the condition ({\bf B}) though the condition ({\bf C}) holds for it.
Actually, in addition to irreducible Hermitian symmetric spaces of
compact type, there also exist countably K\"ahler C-spaces
associated with a complex simple Lie algebra of classical type that
have positive holomorphic sectional curvature.

We may obtain the following theorem, which generalizes
Theorems~\ref{th:1.1},~\ref{th:1.2}, but partially contains
~\ref{th:1.3}.

\begin{theorem}\label{th:6.1}
Let $(M, \omega, J, g)$ and $(\widetilde{M}, \tilde\omega, \tilde
J, \tilde g)$ be two real $2n$-dimensional compact
K\"{a}hler-Einstein manifolds satisfying the above conditions
({\bf A}) and ({\bf B}). Then for any $\Lambda$-pinched
symplectomorphism $\varphi:(M, \omega)\to (\widetilde M,
\widetilde\omega)$  with $\Lambda\in [1,
\Lambda_1(n)]\setminus\{\infty\}$, where $\Lambda_1(n)$ is given
by (\ref{e:1.1}), the following conclusions hold:
\begin{description}
\item[(i)]  The mean curvature flow $\Sigma_t$ of the graph of $\varphi$ in
$M\times\widetilde{M}$ exists smoothly for all $t>0$.

\item[(ii)] $\Sigma_{t}$ is the graph of a symplectomorphism
$\varphi_{t}$ for each $t>0$, and $\varphi_{t}$ is
$\Lambda'_n$-pinched along the mean curvature flow, where
$\Lambda'_n$ is defined by (\ref{e:1.2}).

\item[(iii)] If $\Lambda<\widehat\Lambda_1$ for some $\Lambda_1\in
(\Lambda, \Lambda_1(n)]\setminus\{\infty\}$, where
$\widehat\Lambda_1>1$ is a constant determined by $\Lambda_1$
and\, $n$ (see Lemma~\ref{lem:5.2}), then
 the flow converges to a Lagrangian submanifold of $M\times\widetilde{M}$ as
$t\rightarrow\infty$.

\item[(iv)] The flow converges to a totally geodesic Lagrangian
submanifold of $M\times\widetilde{M}$ and $\varphi_{t}$ converges
smoothly to a biholomorphic isometry from $M$ to $\widetilde{M}$
as $t\rightarrow\infty$ provided additionally that $(M, \omega, J,
g)$ and $(\widetilde{M}, \tilde\omega, \tilde J, \tilde g)$
satisfy the condition ({\bf C}). Consequently, the
symplectomorphism $\varphi: M\rightarrow\widetilde{M}$
 is symplectically isotopic to a biholomorphic isometry.
\end{description}
\end{theorem}

In order to prove it  we start with two simple lemmas.

\begin{lemma}\label{lem:6.2}
Let $R$ be the curvature tensor of a K\"ahler manifold $(M,g, J,
\omega)$ of real dimension $2N$.  For any local holomorphic
coordinate system $(z^1,\cdots, z^n)$ on it, let
$R_{r\overline{s}r\overline{s}}=R\bigr(\frac{\partial}{\partial
z^{r}}, \frac{\partial}{\partial \bar z^{s}},
\frac{\partial}{\partial z^{r}}, \frac{\partial}{\partial \bar
z^{s}}\bigl)$ and $z^s=x^s+ \sqrt{-1}y^s$, $s=1,\cdots, n$.  Then
\begin{eqnarray*}
R\Bigl(\frac{\partial}{\partial x^{s}},\frac{\partial}{\partial
y^{r}},\frac{\partial}{\partial x^{s}},\frac{\partial}{\partial
y^{r}}\Bigr)-R\Bigl(\frac{\partial}{\partial
x^{s}},\frac{\partial}{\partial x^{r}},\frac{\partial}{\partial
x^{s}},\frac{\partial}{\partial x^{r}}\Bigr)=-4{\rm
Re}(R_{r\overline{s}r\overline{s}})\quad\forall r,s.
\end{eqnarray*}
In particular, we have
$$
R\Bigl(\frac{\partial}{\partial x^{s}},\frac{\partial}{\partial
y^{s}},\frac{\partial}{\partial x^{s}},\frac{\partial}{\partial
y^{s}}\Bigr)=-4R_{s\overline{s}s\overline{s}}\;\quad\forall s,
$$
that is, the holomorphic sectional curvature in the direction
$\frac{\partial}{\partial x^{s}}$ is given by
$$
H\bigl(\frac{\partial}{\partial x^s}\bigr)=- \frac{4R_{s\bar ss\bar
s}}{[g(\frac{\partial}{\partial x^s}, \frac{\partial}{\partial
x^s})]^2}.
$$
\end{lemma}

\noindent{\bf Proof}. Since the only possible non-vanishing terms of
the curvature components are of the form $R_{i\bar jk\bar l}$ and
those obtained from the universal symmetries of the curvature
tensor, it is not hard to prove that
\begin{eqnarray}\label{e:6.1}
&&R\Bigl(\frac{\partial}{\partial x^{s}},\frac{\partial}{\partial
y^{r}},\frac{\partial}{\partial x^{s}},\frac{\partial}{\partial
y^{r}}\Bigr)\nonumber\\
&=&R\biggr(\frac{\partial}{\partial z^{s}}+ \frac{\partial}{\partial
\bar z^{s}}, \sqrt{-1}\Bigl(\frac{\partial}{\partial z^{r}}-
\frac{\partial}{\partial \bar z^{r}}\Bigl), \frac{\partial}{\partial
z^{s}}+ \frac{\partial}{\partial \bar z^{s}},
\sqrt{-1}\Bigl(\frac{\partial}{\partial z^{r}}-
\frac{\partial}{\partial
\bar z^{r}}\Bigr)\biggl)\nonumber\\
&=&-(R_{r\overline{s}r\overline{s}}+R_{s\overline{r}s\overline{r}}+
R_{r\overline{s}s\overline{r}}+R_{s\overline{r}r\overline{s}})
\end{eqnarray}
and
\begin{eqnarray}\label{e:6.2}
&&R\Bigl(\frac{\partial}{\partial x^{s}},\frac{\partial}{\partial
x^{r}},\frac{\partial}{\partial x^{s}},\frac{\partial}{\partial
x^{r}}\Bigr)\nonumber\\
&=&R\Bigr(\frac{\partial}{\partial z^{s}}+ \frac{\partial}{\partial
\bar z^{s}}, \frac{\partial}{\partial z^{r}}+
\frac{\partial}{\partial \bar z^{r}}, \frac{\partial}{\partial
z^{s}}+ \frac{\partial}{\partial \bar z^{s}},
\frac{\partial}{\partial z^{r}}+ \frac{\partial}{\partial
\bar z^{r}}\Bigr)\nonumber\\
&=&R_{r\overline{s}r\overline{s}}+R_{s\overline{r}s\overline{r}}-R_{r\overline{s}s\overline{r}}
-R_{s\overline{r}r\overline{s}}.
\end{eqnarray}
Note that
$R_{r\overline{s}r\overline{s}}=R_{\overline{s}r\overline{s}r}=\overline{R_{s\overline{r}s\overline{r}}}$.
It follows from this and (\ref{e:6.1})-(\ref{e:6.2}) that
\begin{eqnarray*}
&&R\Bigl(\frac{\partial}{\partial x^{s}},\frac{\partial}{\partial
y^{r}},\frac{\partial}{\partial x^{s}},\frac{\partial}{\partial
y^{r}}\Bigr)-R\Bigl(\frac{\partial}{\partial
x^{s}},\frac{\partial}{\partial x^{r}},\frac{\partial}{\partial
x^{s}},\frac{\partial}{\partial
x^{r}}\Bigr)\nonumber\\
&=&-2R_{r\overline{s}r\overline{s}}-
2R_{s\overline{r}s\overline{r}}=-4{\rm
Re}(R_{r\overline{s}r\overline{s}}).
\end{eqnarray*}
The second equality may be derived from (\ref{e:6.1}) directly.
Lemma~\ref{lem:6.2} is proved. $\Box$\vspace{2mm}

\begin{lemma}\label{lem:6.3}
Under the assumptions of Lemma~\ref{lem:6.2}, if $(M,g, J, \omega)$
also satisfies the condition ({\bf B}),  then ${\rm
Re}(R_{r\overline{s}r\overline{s}})\leq 0$ for all $1\leq r,s\leq
n$.
\end{lemma}

\noindent{\bf Proof.} Set
$X=\sum_{i=1}^{n}u_{i}\frac{\partial}{\partial z^{i}}$ and
$Y=\sum_{j=1}^{n}v_{j}\frac{\partial}{\partial z^{j}}$ with $u_{i},
v_{j}\in\C$. Then
\begin{eqnarray*}
R(X,\overline{Y},X,\overline{Y})=\sum_{i,j,k,l=1}^{n}R\Bigl(u_{i}\frac{\partial}{\partial
z^{i}}, \overline{v}_{j}\frac{\partial}{\partial \overline{z}^{j}},
u_{k}\frac{\partial}{\partial z^{k}} ,
\overline{v}_{l}\frac{\partial}{\partial \overline{z}^{l}} \Bigr)
=\sum_{i,j,k,l=1}^{n}u_{i}u_{k}\overline{v}_{j}\overline{v}_{l}R_{i\overline{j}k\overline{l}}
\end{eqnarray*}
and
\begin{eqnarray*}
R(Y,\overline{X},Y,\overline{X})=\sum_{i,j,k,l=1}^{n}R\Bigl(v_{j}\frac{\partial}{\partial
z^{j}}, \overline{u}_{i}\frac{\partial}{\partial \overline{z}^{i}},
v_{l}\frac{\partial}{\partial z^{l}} ,
 \overline{u}_{k}\frac{\partial}{\partial \overline{z}^{k}}\Bigr)
=\sum_{i,j,k,l=1}^{n}v_{j}v_{l}\overline{u}_{i}\overline{u}_{k}R_{j\overline{i}l\overline{k}}.
\end{eqnarray*}
Since
$R_{j\overline{i}l\overline{k}}=\overline{R_{i\overline{j}k\overline{l}}}$
we get
\begin{eqnarray*}
R(X,\overline{Y},X,\overline{Y})+R(Y,\overline{X},Y,\overline{X})
&=&\sum_{i,j,k,l=1}^{n}\Bigl(u_{i}u_{k}\overline{v}_{j}\overline{v}_{l}R_{i\overline{j}k\overline{l}}
+\overline{u_{i}u_{k}\overline{v}_{j}\overline{v}_{l}R_{i\overline{j}k\overline{l}}}\,\Bigr)\\
&=&R(X,\overline{Y},X,\overline{Y})+\overline{R(X,\overline{Y},X,\overline{Y})}\\
&=& 2{\rm Re}(R(X,\overline{Y},X,\overline{Y})).
\end{eqnarray*}
Taking $X=\frac{\partial}{\partial z^{r}}$,
$Y=\frac{\partial}{\partial z^{s}}$, the desired results are
obtained. $\Box$\vspace{2mm}

The following proposition implies Theorem~\ref{th:6.1}(i) and
(ii).

\begin{proposition}\label{prop:6.4}
Let $(M, \omega, J, g)$ be a real $2n$-dimensional compact
K\"{a}hler-Einstein manifold satisfying the conditions ({\bf A})
and ({\bf B}). Then for  any  symplectomorphism $\varphi:
M\rightarrow \widetilde{M}$ it holds that
\begin{eqnarray}\label{e:6.3-}
\frac{d}{dt}\ast\Omega\geq\Delta\ast\Omega+\ast\Omega\cdot
Q(\lambda_{i},h_{jkl}),
\end{eqnarray}
 along the mean curvature flow $\Sigma_t$ of the graph $\Sigma$ of
 $\varphi$. Furthermore, if $\varphi$ is $\Lambda$-pinched for
 some $\Lambda\in (1, \Lambda_1(n))$,
  then the symplectomorphism $\varphi_t:M\rightarrow
\widetilde{M}$, whose graph is $\Sigma_t$,
 is $\Lambda'_n$-pinched and
\begin{eqnarray}\label{e:6.3}
\frac{d}{dt}\ast\Omega\geq\Delta\ast\Omega+\delta_{\Lambda}\cdot\ast\Omega|II|^{2}
\end{eqnarray}
along the mean curvature flow. In particular,
$\min_{\Sigma_t}\ast\Omega$ is nondecreasing as a function in $t$.
\end{proposition}

\noindent{\bf Proof}. By the condition ({\bf A}),
$R_{ikik}=\widetilde{R}_{ikik}\;\forall i,k$. Hence the second term
in the big bracket of (\ref{e:3.4}) can be written as follows (
omitting $|_p$ in $\frac{\partial}{\partial x^{t}}|_p$ and
$\frac{\partial}{\partial y^{t}}|_p$),
\begin{eqnarray}\label{e:6.4}
&&\sum_{k}\sum_{i\neq
k}\frac{\lambda_{i}(R_{ikik}-\lambda_{k}^{2}\widetilde{R}_{ikik})}
{(1+\lambda_{k}^{2})(\lambda_{i}+\lambda_{i'})}\\
&=&\sum_{k}\sum_{i\neq
k}\frac{\lambda_{i}(1-\lambda_{k}^{2})R_{ikik}}
{(1+\lambda_{k}^{2})(\lambda_{i}+\lambda_{i'})}\nonumber\\
&=&\sum_{k=2r-1,i=2s-1,r\neq
s}\frac{\lambda_{2s-1}(1-\lambda_{2r-1}^{2})R(\frac{\partial}{\partial
x^{s}},\frac{\partial}{\partial x^{r}},\frac{\partial}{\partial
x^{s}},\frac{\partial}{\partial
x^{r}})}{(1+\lambda_{2r-1}^{2})(\lambda_{2s-1}+\lambda_{2s})}\nonumber\\
&&+\sum_{k=2r-1,i=2s}\frac{\lambda_{2s}(1-\lambda_{2r-1}^{2})R(\frac{\partial}{\partial
y^{s}},\frac{\partial}{\partial x^{r}},\frac{\partial}{\partial
y^{s}},\frac{\partial}{\partial
x^{r}})}{(1+\lambda_{2r-1}^{2})(\lambda_{2s-1}+\lambda_{2s})}\nonumber\\
&&+\sum_{k=2r,i=2s-1}\frac{\lambda_{2s-1}(1-\lambda_{2r}^{2})R(\frac{\partial}{\partial
x^{s}},\frac{\partial}{\partial y^{r}},\frac{\partial}{\partial
x^{s}},\frac{\partial}{\partial
y^{r}})}{(1+\lambda_{2r}^{2})(\lambda_{2s-1}+\lambda_{2s})}\nonumber\\
&&+\sum_{k=2r,i=2s,r\neq
s}\frac{\lambda_{2s}(1-\lambda_{2r}^{2})R(\frac{\partial}{\partial
y^{s}},\frac{\partial}{\partial y^{r}},\frac{\partial}{\partial
y^{s}},\frac{\partial}{\partial
y^{r}})}{(1+\lambda_{2r}^{2})(\lambda_{2s-1}+\lambda_{2s})}\nonumber\\
&=&\sum_{r\neq s}\frac{R(\frac{\partial}{\partial
x^{s}},\frac{\partial}{\partial x^{r}},\frac{\partial}{\partial
x^{s}},\frac{\partial}{\partial
x^{r}})}{(\lambda_{2s-1}+\lambda_{2s})}
\left[\frac{\lambda_{2s-1}(1-\lambda_{2r-1}^{2})}{(1+\lambda_{2r-1}^{2})}
+\frac{\lambda_{2s}(1-\lambda_{2r}^{2})}{(1+\lambda_{2r}^{2})}\right]\nonumber\\
&&+\sum_{r,s}\frac{R(\frac{\partial}{\partial
x^{s}},\frac{\partial}{\partial y^{r}},\frac{\partial}{\partial
x^{s}},\frac{\partial}{\partial
y^{r}})(\lambda_{2r}^{2}-1)(\lambda_{2s}-\lambda_{2s-1})}{(\lambda_{2s-1}+\lambda_{2s})(1+\lambda_{2r}^{2})}\nonumber
\end{eqnarray}
\begin{eqnarray*}
&=&\sum_{r,s}\frac{(\lambda_{2r}^{2}-1)(\lambda_{2s}-\lambda_{2s-1})}{(\lambda_{2s-1}+\lambda_{2s})(1+\lambda_{2r}^{2})}
\left[R(\frac{\partial}{\partial x^{s}},\frac{\partial}{\partial
y^{r}},\frac{\partial}{\partial x^{s}},\frac{\partial}{\partial
y^{r}})-R(\frac{\partial}{\partial x^{s}},\frac{\partial}{\partial
x^{r}},\frac{\partial}{\partial x^{s}},\frac{\partial}{\partial
x^{r}})\right]\nonumber\\
&=&\sum_{r,s}\frac{(\lambda_{2r}^{2}-1)(\lambda_{2s}^{2}-1)}{(1+\lambda_{2s}^{2})(1+\lambda_{2r}^{2})}
\left[R(\frac{\partial}{\partial x^{s}},\frac{\partial}{\partial
y^{r}},\frac{\partial}{\partial x^{s}},\frac{\partial}{\partial
y^{r}})-R(\frac{\partial}{\partial x^{s}},\frac{\partial}{\partial
x^{r}},\frac{\partial}{\partial x^{s}},\frac{\partial}{\partial
x^{r}})\right]\nonumber\\
&=&\sum_{r\neq
s}\frac{(\lambda_{2r}^{2}-1)(\lambda_{2s}^{2}-1)}{(1+\lambda_{2s}^{2})(1+\lambda_{2r}^{2})}
\left[R(\frac{\partial}{\partial x^{s}},\frac{\partial}{\partial
y^{r}},\frac{\partial}{\partial x^{s}},\frac{\partial}{\partial
y^{r}})-R(\frac{\partial}{\partial x^{s}},\frac{\partial}{\partial
x^{r}},\frac{\partial}{\partial x^{s}},\frac{\partial}{\partial
x^{r}})\right]\nonumber\\
&&+\sum_{r=s}\frac{(\lambda_{2r}^{2}-1)(\lambda_{2s}^{2}-1)}{(1+\lambda_{2s}^{2})(1+\lambda_{2r}^{2})}
\left[R(\frac{\partial}{\partial x^{s}},\frac{\partial}{\partial
y^{r}},\frac{\partial}{\partial x^{s}},\frac{\partial}{\partial
y^{r}})-R(\frac{\partial}{\partial x^{s}},\frac{\partial}{\partial
x^{r}},\frac{\partial}{\partial x^{s}},\frac{\partial}{\partial
x^{r}})\right]\nonumber\\
&=&\sum_{r\neq
s}\frac{(\lambda_{2r}^{2}-1)(\lambda_{2s}^{2}-1)}{(1+\lambda_{2s}^{2})(1+\lambda_{2r}^{2})}
\left[-4{\rm Re}(R_{r\overline{s}r\overline{s}})\right]
+\sum_{r=s}\frac{(\lambda_{2r}^{2}-1)(\lambda_{2s}^{2}-1)}{(1+\lambda_{2s}^{2})(1+\lambda_{2r}^{2})}
\left[-4{\rm Re}(R_{s\overline{s}s\overline{s}})\right]         \nonumber\\
&\ge & 0\nonumber
\end{eqnarray*}
because of Lemmas~\ref{lem:6.2},~\ref{lem:6.3} and our choice that
$\lambda_{2i-1}\leq 1\leq\lambda_{2i}, i=1,\cdots,n$. This leads
to (\ref{e:6.3-}).

Now if $\varphi$ is $\Lambda$-pinched, then
$\frac{1}{\Lambda}\le\lambda_i(0)\le\Lambda$ for $i=1,\cdots,2n$.
Since $\Lambda_1(n)<\Lambda_0(n)$ in the case
$\Lambda_0(n)<\infty$, by Proposition~\ref{prop:3.4} we get
$Q(\lambda_{i}(0),h_{jkl})\geq\delta_{\Lambda}\sum_{ijk}h_{jkl}^{2}$
and hence
$$
\left(\frac{d}{dt}-\triangle\right)\ast\Omega\ge
0\quad\hbox{at}\quad t=0.
$$
Note that Lemma 5 of \cite{MeWa} implies that
$\frac{1}{2^{n}}-\epsilon(n, \Lambda)\le\ast\Omega$ at $t=0$,
where $\epsilon(n,
\Lambda)=\frac{1}{2^{n}}-\frac{1}{(\Lambda+\frac{1}{\Lambda})^{n}}$.
Then repeating the proof of Proposition 4 and Corollary 5 in
\cite{MeWa} we may get (\ref{e:6.3}).
 $\Box$\vspace{2mm}

Using this proposition we may prove the long-time existence in
Theorem~\ref{th:6.1} (i) as in \cite[\S 3.3]{MeWa} (or that of
Theorem~\ref{th:1.1}).
\

\noindent{\bf The proof of Theorem~\ref{th:6.1}(iii)}.  The idea
is similar to that of Theorem~\ref{th:1.3}. All arguments from the
beginning of Section 4.2.2 to (\ref{e:5.6}) in the proof of
convergence in Theorem~\ref{th:1.3} are still valid. Then  there
exists a positive number $K_2$ depending on the manifolds $M$ and
$\widetilde{M}$ such that
\begin{eqnarray*}
\sum_{s,i,j}\left(\sum_{k}[(\overline{\nabla}_{\partial_k}\overline{R})_{\underline{s}ijk}
+(\overline{\nabla}_{\partial_j}\overline{R})_{\underline{s}kik})]\right)^{2}\leq
K_2
\end{eqnarray*}
and hence
\begin{eqnarray*}
\sum_{s,i,j,k}2\left[(\overline{\nabla}_{\partial_k}\overline{R})_{\underline{s}ijk}
+(\overline{\nabla}_{\partial_j}\overline{R})_{\underline{s}kik})\right]h_{sij}\leq
K_2 +|II|^{2}.
\end{eqnarray*}
As there it follows from the boundedness of the curvature that
\begin{equation}\label{e:6.5}
\frac{d}{dt}|II|^{2}\leq\Delta|II|^{2}-2|\nabla II|^{2}+ 10
|II|^{4}+K_1 |II|^{2}+K_2,
\end{equation}
 where  $K_1$ is a nonnegative constant that depends on the
dimensions of $M$ and $\widetilde{M}$. With the same proof we may
get the corresponding result of (\ref{e:5.9}), i.e.
\begin{eqnarray}\label{e:6.6}
&&(\frac{d}{dt}-\Delta)\left(\frac{|II|^2}{\sin(k(\ast\Omega)^l)}\right)\\
&\leq&\left(\frac{|II|^{2}}{\sin(k(\ast\Omega)^l)}\right)^2\cdot\bigl[10\cdot\sin(k(\ast\Omega)^l)-
k\cdot l\cdot \delta_{\Lambda_1}\cdot(\ast\Omega)^l \cdot \cos(k(\ast\Omega)^l)\bigr]\nonumber\\
&&+\frac{k\cdot
l\cdot(\ast\Omega)^{l-2}|II|^{2}|\nabla\ast\Omega|^{2}}{[\sin(k(\ast\Omega)^l)]^{3}}\bigl[-k\cdot
l\cdot(\ast\Omega)^{l}\cdot
(\sin(k(\ast\Omega)^l))^2\nonumber\\
&&\hspace{60mm}+(l-1)\cos(k(\ast\Omega)^l)\sin(k(\ast\Omega)^l)\bigr]\nonumber\\
&&+K_1
\frac{|II|^{2}}{\sin(k(\ast\Omega)^l)}+\frac{K_2}{\sin(k(\ast\Omega)^l)}.\nonumber
\end{eqnarray}
By Claim~\ref{cl:5.3}, it follows from (\ref{e:6.6}) that
\begin{eqnarray*}
&&(\frac{d}{dt}-\Delta)\left(\frac{|II|^2}{\sin(k(\ast\Omega)^l)}\right)\\
&\leq& \left(\frac{|II|^{2}}{\sin(k(\ast\Omega)^l)}\right)^2
[10\cdot\sin(k(\ast\Omega)^l)-
 k\cdot l\cdot \delta_{\Lambda_1} \cdot(\ast\Omega)^l \cdot \cos(k(\ast\Omega)^l)]\\
 &&+K_1 \frac{|II|^{2}}{\sin(k(\ast\Omega)^l)}+\frac{K_2}{\sin(k(\ast\Omega)^l)}.
\end{eqnarray*}
Let $g=\frac{|II|^2}{\sin(k(\ast\Omega)^l)}$,
$K_{4}:=\max{\frac{K_2}{\sin(k(\ast\Omega)^l)}}
=\frac{K_2}{\sin(k(\frac{1}{(\Lambda+\frac{1}{\Lambda})^{n}})^l)}$
and
$$
K_{3}:=\max_{\ast\Omega\in[\frac{1}{(\Lambda+\frac{1}{\Lambda})^{n}},
\frac{1}{2^{n}}]}\left[10\cdot\sin(k(\ast\Omega)^l)- k\cdot l\cdot
\delta_{\Lambda_1} \cdot(\ast\Omega)^l \cdot
\cos(k(\ast\Omega)^l)\right].
$$
By Claim~\ref{cl:5.3}, $K_{3}< 0$ and
\begin{eqnarray}\label{e:6.7}
(\frac{d}{dt}-\Delta)g\leq K_{3}\cdot g^{2}+K_1\cdot g+K_4.
\end{eqnarray}

Consider the initial value problem
\begin{equation}\label{e:6.8}
\frac{d}{dt}y= K_{3}\cdot y^{2}+K_1\cdot y +K_4\quad\hbox{and}\quad
y(0)=\rm max_{\Sigma_{0}}g.
\end{equation}

If $y(0)>\frac{-K_{1}-\sqrt{K_{1}^{2}-4K_{3}K_{4}}}{2K_{3}}$, the
unique solution of (\ref{e:6.8}) is given by
\begin{eqnarray*}
y(t)=\frac{(K_{1}+\sqrt{K_{1}^{2}-4K_{3}K_{4}})\cdot\exp(\sqrt{K_{1}^{2}
-4K_{3}K_{4}}t+K_{5})-K_{1}+\sqrt{K_{1}^{2}-4K_{3}K_{4}}}
{-2K_{3}\cdot[\exp(\sqrt{K_{1}^{2}-4K_{3}K_{4}}t+K_{5})-1]},
\end{eqnarray*}
where
$K_5=\ln\frac{2K_{3}y(0)+K_{1}-\sqrt{K_{1}^{2}-4K_{3}K_{4}}}{2K_{3}y(0)+K_{1}+\sqrt{K_{1}^{2}-4K_{3}K_{4}}}$.
Clearly, $y(t)\rightarrow
\frac{K_{1}+\sqrt{K_{1}^{2}-4K_{3}K_{4}}}{-2K_{3}}$ as
$t\rightarrow\infty$.

If $y(0)=\frac{-K_{1}-\sqrt{K_{1}^{2}-4K_{3}K_{4}}}{2K_{3}}$, then
$y(t)\equiv\frac{-K_{1}-\sqrt{K_{1}^{2}-4K_{3}K_{4}}}{2K_{3}}$.

If $y(0)<\frac{-K_{1}-\sqrt{K_{1}^{2}-4K_{3}K_{4}}}{2K_{3}}$, then
there exists a $T>0$
 such that on $[0, T]$ we have $y(t)-\frac{-K_{1}-\sqrt{K_{1}^{2}-4K_{3}K_{4}}}{2K_{3}}\leq 0$, and therefore
$$
 \left(y(t)+\frac{K_{1}}{2K_{3}}\right)^{2}=-\exp(\sqrt{K_{1}^{2}-4K_{3}K_{4}}t+K_{5})+\frac{K_{1}^{2}-4K_{3}K_{4}}{4K_{3}^{2}}\geq 0
$$
where
$K_{5}=\ln\left(-y(0)^{2}-\frac{K_{1}y(0)}{K_{3}}-\frac{K_{4}}{K_{3}}\right)$.
It follows that
$$
T=\frac{\ln(\frac{K_{1}^{2}-4K_{3}K_{4}}{4K_{3}^{2}})-K_{5}}{\sqrt{K_{1}^{2}-4K_{3}K_{4}}}\geq
0\quad\hbox{and}\quad
y(T)=-\frac{K_{1}}{2K_{3}}<\frac{-K_{1}-\sqrt{K_{1}^{2}-4K_{3}K_{4}}}{2K_{3}}.
$$
Hence we can continue this procedure and get
$$
 \left(y(t)+\frac{K_{1}}{2K_{3}}\right)^{2}=-\exp(\sqrt{K_{1}^{2}-4K_{3}K_{4}}t+K_{5})+\frac{K_{1}^{2}-4K_{3}K_{4}}{4K_{3}^{2}}\geq 0
$$
 for all time $t\ge 0$. From this we derive
\begin{eqnarray*}
y(t)&=&-\frac{K_{1}}{2K_{3}}+\sqrt{-\exp(\sqrt{K_{1}^{2}-4K_{3}K_{4}}t+K_{5})+\frac{K_{1}^{2}-4K_{3}K_{4}}{4K_{3}^{2}}}\\
&\leq &
-\frac{K_{1}}{2K_{3}}+\sqrt{\frac{K_{1}^{2}-4K_{3}K_{4}}{4K_{3}^{2}}}=\frac{-K_{1}-\sqrt{K_{1}^{2}-4K_{3}K_{4}}}{2K_{3}}
\end{eqnarray*}
if $-\frac{K_{1}}{2K_{3}}\leq
y(0)<\frac{-K_{1}-\sqrt{K_{1}^{2}-4K_{3}K_{4}}}{2K_{3}}$,  and
\begin{eqnarray*}
y(t)&=& -\frac{K_{1}}{2K_{3}}-\sqrt{-\exp(\sqrt{K_{1}^{2}-4K_{3}K_{4}}t+K_{5})+\frac{K_{1}^{2}-4K_{3}K_{4}}{4K_{3}^{2}}}\\
&\leq & -\frac{K_{1}}{2K_{3}}
\end{eqnarray*}
if $0 \leq y(0)<-\frac{K_{1}}{2K_{3}}$.

 By (\ref{e:6.7})-(\ref{e:6.8}) the comparison principle
for parabolic equations yields
\begin{equation}\label{e:6.9}
g=\frac{|II|^2}{\sin(k(\ast\Omega)^l)}\le  y(t)\quad\forall t>0.
\end{equation}
Since (\ref{e:5.3}) implies that the function
$$
\Bigl[\frac{1}{(\Lambda+\frac{1}{\Lambda})^{n}},\frac{1}{2^{n}}\Bigr]\ni
\ast\Omega\to \sin(k(\ast\Omega)^l)
$$
is bounded away from zero, we derive
\begin{equation}\label{e:6.10}
\max_{\Sigma_{t}}|II|^{2}\leq \sin(k(\frac{1}{2^{n}})^l)\cdot
y(t)\leq\sin(k(\frac{1}{2^{n}})^l)\cdot L,
\end{equation}
where $L=\frac{-K_{1}-\sqrt{K_{1}^{2}-4K_{3}K_{4}}}{2K_{3}}$ if
$y(0)\geq-\frac{K_{1}}{2K_{3}}$, and $L=-\frac{K_{1}}{2K_{3}}$ if
$0\leq y(0)<-\frac{K_{1}}{2K_{3}}$. Hence $|II|^2$ is uniformly
bounded. Namely, we have proved that  the flow converges to a
Lagrangian submanifold at infinity provided that the flow exists for
all the time. ({\it Note}: Different from the case of tori we cannot
prove $\max_{\Sigma_{t}}|II|^{2}\to 0$ as $t\to\infty$, and hence
cannot assert that the limit submanifold is  totally geodesic.)

\

\noindent{\bf The proof of Theorem~\ref{th:6.1}(iv)}.  The idea is
similar to that of Theorem~\ref{th:1.1}. In the present case we have
the following

\begin{proposition}\label{prop:6.5}
Under the assumptions of Proposition~\ref{prop:6.4}, suppose further
that $(M, \omega, J, g)$ also satisfies the condition ({\bf C}).
 Then along the mean curvature flow a
similar inequality to that of Proposition~\ref{prop:3.5}  holds,
i.e.
\begin{eqnarray}\label{e:6.11}
\frac{d}{dt}\ast\Omega\geq\Delta\ast\Omega+\delta_{\Lambda}\cdot\ast\Omega|II|^{2}+
c_0\cdot\ast\Omega\sum_{k~{\rm odd}}\frac{(1-\lambda_{k}^{2})^{2}}
{(1+\lambda_{k}^{2})^{2}}.
\end{eqnarray}
\end{proposition}

\noindent{\bf Proof}. Under the further assumption, by (\ref{e:6.4})
we have
\begin{eqnarray*}
&&\sum_{k}\sum_{i\neq
k}\frac{\lambda_{i}(R_{ikik}-\lambda_{k}^{2}\widetilde{R}_{ikik})}
{(1+\lambda_{k}^{2})(\lambda_{i}+\lambda_{i'})}\\
&=&\sum_{r\neq
s}\frac{(\lambda_{2r}^{2}-1)(\lambda_{2s}^{2}-1)}{(1+\lambda_{2s}^{2})(1+\lambda_{2r}^{2})}
\left[-4{\rm Re}(R_{r\overline{s}r\overline{s}})\right]
+\sum_{r=s}\frac{(\lambda_{2r}^{2}-1)(\lambda_{2s}^{2}-1)}{(1+\lambda_{2s}^{2})(1+\lambda_{2r}^{2})}
\left[-4{\rm Re}(R_{s\overline{s}s\overline{s}})\right]         \nonumber\\
&\ge &
c_0\sum_{r=s}\frac{(\lambda_{2r}^{2}-1)(\lambda_{2s}^{2}-1)}{(1+\lambda_{2s}^{2})(1+\lambda_{2r}^{2})}.
\nonumber
\end{eqnarray*}
This and Propositions~\ref{prop:3.3} and ~\ref{prop:3.4} give
(\ref{e:6.11}). $\Box$\vspace{2mm}

As in \cite{MeWa}, using this we may prove that
$\lambda_{i}\rightarrow 1$ and $\max_{\Sigma_{t}}|II|^{2}\rightarrow 0$ as
$t\rightarrow\infty$, and hence that the flow converges to a totally
geodesic Lagrangian submanifold of $M\times\widetilde{M}$ as
$t\rightarrow\infty$ and that $\varphi_{t}$ converges smoothly to a
biholomorphic isometry $\varphi_{\infty}:M\rightarrow
\widetilde{M}$. Theorem~\ref{th:6.1} is proved. $\Box$\vspace{2mm}

  A theorem by Matsushima and Borel-Remmert claimed that
every compact homogeneous K\"ahler manifold is the K\"ahler product
of a flat complex torus (known as the {\it Albanese} torus of $(M,
J)$) and a K\"ahler $C$-space (cf. \cite[Theorem 8.97]{Be}). As a
consequence, a compact homogeneous K\"ahler manifold admits a
K\"ahler-Einstein structure if and only if it is a complex torus or
is simply-connected. If we restrict the manifolds in
Theorem~\ref{th:6.1} to homogeneous K\"ahler-Einstein manifolds,
then Theorem~\ref{th:6.1} has sense only for simply-connected case
(because the better result has been obtained for complex tori).

\appendix
\section{Appendix:\quad Proof of Claim~\ref{cl:5.4} }\label{app:A}\setcounter{equation}{0}

For the function $g(\alpha)$ in (\ref{e:5.20}),
 a direct computation yields
\begin{eqnarray}
&&g'(\alpha)=\frac{1}{(\sin \alpha)^{2}}\biggl[\sin
\alpha\cdot\cos\alpha\cdot\ln\big(\alpha/\sqrt{\frac{\sqrt{21}-3}{2}}\bigr)\nonumber\\
&&\hspace{40mm}+\sin
\alpha\cdot\cos\alpha-\alpha\ln\big(\alpha/\sqrt{\frac{\sqrt{21}-3}{2}}\bigr)\biggr],\nonumber\\
&& g''(\alpha)=\frac{1}{(\sin \alpha)^{3}}\biggl[\frac{(\sin
\alpha)^{2}\cdot\cos\alpha}{\alpha}
+2\alpha\cos\alpha\cdot\ln\big(\alpha/\sqrt{\frac{\sqrt{21}-3}{2}}\bigr)\nonumber\\
&&\hspace{40mm}-2\sin \alpha -2\sin
\alpha\cdot\ln\big(\alpha/\sqrt{\frac{\sqrt{21}-3}{2}}\bigr)\biggr].\label{e:A.3}
\end{eqnarray}
Clearly,
$\lim_{\alpha\rightarrow\frac{\pi}{2}}g(\alpha)=0=g(\sqrt{\frac{\sqrt{21}-3}{2}})$,
and $g(\alpha)>0$ on $(\sqrt{\frac{\sqrt{21}-3}{2}},\frac{\pi}{2})$.
Moreover,
$$
g'(\frac{\pi}{2})=-\frac{\pi}{2}\ln\big(\pi/2\sqrt{\frac{\sqrt{21}-3}{2}}\bigr)<0\quad\hbox{and}\quad
g'(\sqrt{\frac{\sqrt{21}-3}{2}})=1/\tan\sqrt{\frac{\sqrt{21}-3}{2}}>0.
$$
(Note that $\sqrt{\frac{\sqrt{21}-3}{2}}\approx 0.8895436175241$
sits between $\frac{\pi}{3.5317}$ and $\frac{\pi}{3.5316}$).
 Hence $g(\alpha)$ attains its maximum at some
point $\alpha_{0}\in (\sqrt{\frac{\sqrt{21}-3}{2}},\frac{\pi}{2})$
with $g'(\alpha_0)=0$. Since any zero $\alpha$ of $g'$ in
$(\sqrt{\frac{\sqrt{21}-3}{2}},\frac{\pi}{2})$  satisfies the
following equation
\begin{eqnarray}\label{e:A.4}
\sin
\alpha\cdot\cos\alpha\cdot\ln\big(\alpha/\sqrt{\frac{\sqrt{21}-3}{2}}\bigr)
+\sin
\alpha\cdot\cos\alpha-\alpha\ln\big(\alpha/\sqrt{\frac{\sqrt{21}-3}{2}}\bigr)=0,
\end{eqnarray}
plugging (\ref{e:A.4}) into (\ref{e:A.3}) we get
\begin{eqnarray}\label{e:A.5}
g''(\alpha) &=&\frac{1}{(\sin \alpha)^{3}}\biggl[\frac{(\sin
\alpha)^{2}\cdot\cos\alpha}{\alpha}
+2\alpha\cos\alpha\cdot\ln\big(\alpha/\sqrt{\frac{\sqrt{21}-3}{2}}\bigr)\nonumber\\
&&\hspace{40mm}-2\sin \alpha
-2\sin \alpha\cdot\ln\big(\alpha/\sqrt{\frac{\sqrt{21}-3}{2}}\bigr)\biggr]\nonumber\\
&=&\frac{1}{(\sin \alpha)^{3}}\biggl[\frac{(\sin
\alpha)^{2}\cdot\cos\alpha}{\alpha}
+2\cos\alpha\cdot(\sin \alpha\cdot\cos\alpha)(1+\ln\big(\alpha/\sqrt{\frac{\sqrt{21}-3}{2}}\bigr))\nonumber\\
&&-2\sin \alpha
-2\sin \alpha\cdot\ln\big(\alpha/\sqrt{\frac{\sqrt{21}-3}{2}}\bigr)\biggr]\nonumber\\
&=&\frac{1}{(\sin \alpha)^{3}}\biggl[\frac{(\sin
\alpha)^{2}\cdot\cos\alpha}{\alpha}-
2(\sin \alpha)^{3}-2(\sin \alpha)^{3}\ln\big(\alpha/\sqrt{\frac{\sqrt{21}-3}{2}}\bigr)\biggr]\nonumber\\
&=&\frac{1}{\sin \alpha}\biggl[\frac{\cos\alpha}{\alpha}-2\sin
\alpha- 2\sin
\alpha\cdot\ln\big(\alpha/\sqrt{\frac{\sqrt{21}-3}{2}}\bigr)\biggr].
\end{eqnarray}
Observe that the function $u(\alpha)=\frac{\cos\alpha}{\alpha}$ is
decreasing on $(\sqrt{\frac{\sqrt{21}-3}{2}},\frac{\pi}{2})$ because
of $u'(\alpha)=-\frac{\sin
\alpha}{\alpha}-\frac{\cos\alpha}{\alpha^{2}}<0$. From
$\sqrt{\frac{\sqrt{21}-3}{2}}\approx 0.8895436175241$ we derive
\begin{eqnarray}\label{e:A.6}
&&\frac{\cos\alpha}{\alpha}-2\sin \alpha-2\sin \alpha\cdot\ln\big(\alpha/\sqrt{\frac{\sqrt{21}-3}{2}}\bigr)\nonumber\\
&<&\frac{\cos\alpha}{\alpha}-2\sin \alpha
\leq\frac{\cos\sqrt{\frac{\sqrt{21}-3}{2}}}
{\sqrt{\frac{\sqrt{21}-3}{2}}} -2\sin
\sqrt{\frac{\sqrt{21}-3}{2}}<0.
\end{eqnarray}
That is,  $g''(\alpha)<0$ for any zero $\alpha$ of $g'$ in
$(\sqrt{\frac{\sqrt{21}-3}{2}},\frac{\pi}{2})$. It follows that
each zero $\alpha$ of of $g'$ in
$(\sqrt{\frac{\sqrt{21}-3}{2}},\frac{\pi}{2})$ is a local maximum
point of $g$. This implies that $g'$ has a unique zero $\alpha_0$
in $(\sqrt{\frac{\sqrt{21}-3}{2}},\frac{\pi}{2})$ and that
$$
g(\alpha_{0})=(\cos
\alpha_{0})^{2}(1+\ln\big(\alpha_{0}/\sqrt{\frac{\sqrt{21}-3}{2}}\bigr))
=\frac{\alpha_{0}(\cos \alpha_{0})^{2}}{\alpha_{0}-\sin
\alpha_{0}\cdot\cos \alpha_{0}}
$$
is the maximum of $g$ in
$(\sqrt{\frac{\sqrt{21}-3}{2}},\frac{\pi}{2})$. We can compute
$\alpha_0\approx 1.238756$ and $g(\alpha_0)\approx 0.141446$.


\begin{thebibliography}{L3}

\bibitem{Ab} M. Abreu,  Topology of symplectomorphism groups of $S^2\times
S^2$. {\it Invent. Math.}, {\bf 131}(1998), no. 1, 1--23.

\bibitem{AbMc} M. Abreu and D. McDuff,  Topology of symplectomorphism groups of
rational ruled surfaces. {\it J. Amer. Math. Soc.}, {\bf 13}(2000),
no. 4, 971--1009.

\bibitem{AnGr} S. Anjos and G. Granja,  Homotopy decomposition of a
group of symplectomorphisms of $S^2\times S^2$. {\it Topology}, {\bf
43}(2004), no. 3, 599--618.

\bibitem{Ban} A. Banyaga, {\it The structure of classical
diffeomorphism groups}. Mathematics and its Applications, 400. {\it
Kluwer Academic Publishers Group, Dordrecht}, 1997.

\bibitem{Be} A. L. Besse, {\it Einstein Manifolds}, Spinger-Verlag,
1987.

\bibitem{Bo} A. Borel, On the curvature tensor of the Hermitian symmetric manifolds,
{\it Annals of Mathematics}, Second Series, {\bf 71}(1960), no.3,
508-521.


\bibitem{CaVe} E. Calabi and E. Vesentini, On Compact, Locally Symmetric Kahler Manifolds,
{\it Annals of Mathematics}, Second Series, {\bf 71}(1960), no.3,
472-507.


\bibitem{Ch} L.S. Charlap,  {\it Bieberbach groups and flat
manifolds}. Universitext. Springer-Verlag, New York, 1986.



\bibitem{Gr} M. Gromov, Pseudo holomorphic curves
in symplectic manifolds, {\it Invent.Math.}, {\bf 82}(1985),307-347.


\bibitem{Ha} R. Hamilton, Three-manifolds with positive Ricci curvature. {\it J.
Differential Geometry}, {\bf 17}(1982), no.2, 255-306.




\bibitem{He} S. Helgason, {\it Differential Geometry, Lie groups and Symmetric Spaces
 }, Academic Press, New York,  1978.

\bibitem{HoZe} H. Hofer and E. Zehnder, {\it Symplectic
invariants and Hamiltonian dynamics}. Birkh\"uuser Verlag, Basel,
1994

\bibitem{KoNo} S. Kobayashi and K. Nomizu, {\it Foundations of
Differential Geometry}, Vol.II, Interscience Publishers, 1969.


\bibitem{LeO}  H\^{o}ng-V\^{a}n L\^{e} and K. Ono, Parameterized Gromov-Witten invariants
and topology of symplectomorphism groups. {\it Groups of
diffeomorphisms}, 51--75, Adv. Stud. Pure Math., 52, Math. Soc.
Japan, Tokyo, 2008.

\bibitem{Le} K. Leichtweiss, Zur Riemannschen Geometrie in
Grassmannschen Mannigfaltigkeiten, {\it Math. Z.}, {\bf 76}(1961),
334-336.



\bibitem{LiWat} Jiayong Li and J. A. Watts, The orienation-preserving diffeomorphism
group of $S^2$ deforms to $SO(3)$ smoothly.  {\it Transformation Groups}, {\bf 16}(2011),
no.2, 537-553.

\bibitem{Lu1} Qi-keng Lu, The elliptic geometry of extended spaces,
{\it Chinese Math.}, {\bf 4}(1963), 54-69. (translation of Acta
Math. Sin., {\bf 13}(1963), 49-62, by the Am. Math.Soc.).

\bibitem{Lu2} Qi-keng Lu, {\it The classical manifolds and the classical domains},
Shanghai Scientific $\&$ Technical Publishers, 1963.

\bibitem{Lu3} Qi-keng Lu, {\it The New Results of the classical manifolds and the classical domains},
Shanghai Scientific $\&$ Technical Publishers, 1997.

\bibitem{Mc} D. McDuff, A survey of the topological properties of
symplectomorphism groups, In: {\it Topology, geometry and quantum
field theory}, 173-193, London math. Soc. Lecture Note Ser., {\bf
308}, Cambridge Univ. Press, 2004.


\bibitem{McSa} D. McDuff and D. Salamon, {\it Introduction to
Symplectic Topology}, Clarendon Press Oxford, 1995.



\bibitem{MeWa} I. Medos and M.-T.Wang, Deforerming Symplectomorphisms of complex
projective spaces by mean curvature flow. {\it J. Diff. Geom.},
{\bf 87}(2011) 309-341.

\bibitem{Si} L. Simon, Asymptotics for a class of nonlinear evolution equations,
 with applications to geometric problems. {\it Ann. of Math.}, {\bf (2) 118}(1983),
 no. 3, 525-571.

\bibitem{Mok} Ngaiming Mok, {\it Metric Rigidity Theorems on Hermitian
Locally Symmetric Manifolds}, Series in Pure Mathematics Vol.6,
World Scientic Publishing Co. Pte. Ltd, 1989.

\bibitem{Po} L. Polterovich, {\it The geometry of the group of symplectic diffeomorphisms}.
Lectures in Mathematics ETH Z\"uich. Birkh\"auser Verlag, Basel,
2001.

\bibitem{Sei} P. Seidel,  On the group of symplectic automorphisms of $\bold
C{\rm P}^m\times\bold C{\rm P}^n$. {\it Northern California
Symplectic Geometry Seminar}, 237--250, Amer. Math. Soc. Transl.
Ser. 2, 196, {\it Amer. Math. Soc., Providence}, RI, 1999.

\bibitem{Sma} S. Smale, Diffeomorphisms of the 2-spheres, {\it Proc.
Amer. Math. Soc.}, {\bf 10}(1959),621-626.

\bibitem{Smo1} K.Smoczyk, A canonical way to deform a Lagrangian submanifolds.
preprint, dg-ga/9605005.

\bibitem{Smo2} K.Smoczyk,  Angle theorems for the
Lagrangian mean curvature flow. {\it Math. Z.}, {\bf 240}(2002), no.
4, 849-883.

\bibitem{Smo3} K.Smoczyk, Longtime existence of the Lagrangian mean curvature flow.
{\it Calc.Var.}, {\bf 20}(2004),25-46.


\bibitem{TsWa} M. -P. Tsui and M.-T. Wang, Mean curvature flows and homotopy of
maps between spheres. {\it Comm. Pure App. Math.},  {\bf 57}(2004),
no. 8, 1110-1126.

\bibitem{W} Hsien-Chung Wang, Closed Manifolds with Homogeneous Complex Structure,
{\it American Journal of Mathematics}, {\bf 76}(1954), no.1,  1-32.




\bibitem{Wa1}  M.-T. Wang,  Mean curvature flow in higher codimension.  math.DG (math.AP). math.DG/0204054.

\bibitem{Wa2}  M.-T. Wang, Mean curvature flow of surface in Einstein Four Manifolds.
{\it J.Differential Geom.},{\bf 57(2)}(2001), 301-338.

\bibitem{Wa3}  M.-T. Wang, Long-time existence and convergence of graphic mean curvature flow
in arbitrary codimension. {\it Invent. Math.}, {\bf 148}(2002), no.
3, 525--543.

\bibitem{Wa4}  M.-T. Wang, Deforming area preserving diffeomorphism of surfaces by mean
curvature flow. {\it Math. Res. Lett.}, {\bf 8}(2001), no. 5-6,
651--661.

\bibitem{Wa5}  M.-T. Wang, A convergence result of the Lagrangian mean curvature flow.
{\it Proceedings of the third International Congress of Chinese
Mathematicians.}  arXiv:math/0508354v1.



\bibitem{Wo1} Yung-Chow, Wong, Differential geometry of Grassmann manifolds, {\it Proc.
N.A.S.}, {\bf 57}(1967), 589¨C594.

\bibitem{Wo2} Yung-Chow,Wong, Sectional Curvatures of Grassmann Manifolds, {\it Proc.
N.A.S.}, {\bf 60}(1968), 75-79.

\end{thebibliography}
\end{document}